\documentclass[12pt,twoside]{report}

\usepackage{amsfonts}       %
\usepackage[leqno]{amsmath} 
\usepackage{amssymb}        %
\usepackage{amstext}        %

\usepackage{amsthm}         
\usepackage{amscd}          
\usepackage{makeidx}        
\makeindex

\usepackage{latexsym}
\usepackage[dvips]{epsfig}
\usepackage[all,dvips]{xy}  
\usepackage{multicol}
\usepackage{enumerate}      

\usepackage{color}

\newtheorem{meinzaehler}{ist quatsch}[chapter]
\newtheorem{thm}[meinzaehler]{Theorem}
\newtheorem{Def}[meinzaehler]{Definition}
\newtheorem{prop}[meinzaehler]{Proposition}
\newtheorem{lemma}[meinzaehler]{Lemma}

\newtheorem{remark}[meinzaehler]{Remark}
\newtheorem{cor}[meinzaehler]{Corollary}
\newtheorem{expl}[meinzaehler]{Example}
\newtheorem{Tabelle}[meinzaehler]{Table}

\newtheorem{conj}[meinzaehler]{Conjecture}

\newcommand{\Numerierung}
    {\refstepcounter{meinzaehler}{\bf(\arabic{chapter}.\arabic{meinzaehler}) }
    }
\newcommand{\meinchapter}[1]{\refstepcounter{chapter}\section*{\arabic{chapter} $\ $ #1}
    \addcontentsline{toc}{chapter}{\arabic{chapter}$\ $ #1}
    \markboth{\scshape \arabic{chapter}. #1}{}
}

\newcommand{\Absatz}{\par\bigskip\par}
\newcommand{\absatz}{\par\medskip\par}
\newcommand{\klabsatz}{\par\smallskip\par}
\newcommand{\neueZeile}{\\}
\newcommand{\Seitenumbruch}{\clearpage}

\newcommand{\N}{\mathbb{N}}
\newcommand{\Z}{\mathbb{Z}}

\newcommand{\C}{\mathbb{C}}
\newcommand{\F}{\mathbb{F}}
\newcommand{\G}{\mathbb{G}}
\newcommand{\NN}{\mathbb{N}}
\newcommand{\ZZ}{\mathbb{Z}}

\newcommand{\GG}{\mathbb{G}}

\newcommand{\CCC}{\mathcal{C}}
\newcommand{\OOO}{\mathcal{O}}

\newcommand{\LLL}{\mathcal{L}}
\newcommand{\NNN}{\mathcal{N}}

\newcommand{\pPp}{\mathfrak{p}}

\newcommand{\mMm}{\mathfrak{m}}

\newcommand{\eps}{\varepsilon}
\newcommand{\pHi}{\varphi}

\newcommand{\tran}{{}^t}
\newcommand{\STO}{O^{st}}

\newcommand{\Cent}{\text{Cent}}  

\newcommand{\GL}{\text{GL}}
\newcommand{\SL}{\text{SL}}
\newcommand{\ORTH}{\text{O}}
\newcommand{\SO}{\text{SO}}
\newcommand{\Sp}{\text{Sp}}
\newcommand{\SP}{\text{Sp}}
\newcommand{\GSp}{\text{GSp}}
\newcommand{\PGL}{\text{PGL}}
\newcommand{\Spin}{\text{Spin}}
\newcommand{\GSPIN}{\text{GSpin}}
\newcommand{\rank}{\text{rank}}

\newcommand{\ov}{\overline}

\newenvironment{engeListe}%
{%
        \begin{list}{}{\setlength{\parsep}{0pt}\setlength{\itemsep}{0pt}
                    \setlength{\labelsep}{8pt}\setlength{\labelwidth}{15pt}
                      }
}{\end{list}}

\definecolor{hellgrau}{gray}{0.9}
\definecolor{mittelgrau}{gray}{0.7}
\definecolor{mittelgrau2}{gray}{0.705}
\definecolor{dunkelgrau}{gray}{0.5}
\definecolor{dunkelgrau2}{gray}{0.51}

\begin{document}

\pagestyle{headings}

\markboth{}{}
\title{Remarks on the Fundamental Lemma for stable twisted Endoscopy of Classical Groups }
\author{Joachim Ballmann\footnote{{supported by Deutsche Forschungsgemeinschaft}}
    \\ Rainer Weissauer
    \\ $\text{Uwe Weselmann}^1$}
\date{\today}
\maketitle

\CompileMatrices
\thispagestyle{empty}
\Seitenumbruch 

\markboth{}{}
\tableofcontents\thispagestyle{empty}        \Seitenumbruch
\pagenumbering{roman}\setcounter{page}{1}
\section*{Introduction}
\addcontentsline{toc}{chapter}{Introduction}
\Absatz
These notes go back to the beginning of the century as the fundamental work of Ngo \cite{Ngo} on the fundamental lemma was not known.
There is some overlap with the paper of Waldspurger \cite{Waldtwist}, which has been written later.
We reproduce the paper  here for historical reasons and to get it available for the public as a reference.

\Absatz
The original aim of these notes was to prove a fundamental lemma for the stable lift from $H=\SP_4$ to
$\tilde G=\widetilde{\PGL}_5$ over a local non archimedean field $F$ with residue characteristic $\ne 2$. Here
$\widetilde{\PGL}_5=\PGL_5\rtimes \langle \Theta\rangle$
is generated by its normal subgroup $\PGL_5$ of index $2$ and
the involution $\Theta:g\mapsto J\cdot \tran g^{-1}\cdot J^{-1}$, where $J$ is the antidiagonal matrix
with entries $(1,-1,1,-1,1)$.
\Absatz
We will (Cor.\ref{BC2TheoremCor}) prove that if the semisimple elements
 $\gamma\Theta\in \widetilde{\PGL}_5(F)$ and $\eta\in \SP_4(F)$ match (see \ref{DefMatching} for
 a definition of matching) then the corresponding stable orbital integrals (see \ref{DefOrbint})
 for the unit elements in the Hecke algebra match:
 \begin{gather}\label{Fule45Gleichung}
\STO_{\gamma\Theta}(1,\widetilde{\PGL}_{5})=\STO_\eta(1, \Sp_{4}).
\end{gather}
This theorem will have important applications in the theory of automorphic representations of
the group $\GSp_4$ over a number field and for the $l$--adic Galois-representations on the corresponding
Shimura varieties \cite{W1}, \cite{W2}, \cite{Uwe}.
\Absatz
In analyzing (\ref{Fule45Gleichung}) using the Kazhdan-trick (lemma \ref{Kazhdan-Lemma} below)
we recognized that all essential computations had already been done by Flicker in \cite{FGL4},
where the corresponding fundamental lemma  for the lift from $\GSp_4\simeq \GSPIN_5$ to
$\widetilde{\GL_4}\times \GG_m$ has been proved. This phenomenon seems to be known to the
experts \cite{HalesSp4SO5}.
\Absatz
More generally one can discuss the fundamental lemma for the stable lift from $H$ to a classical
split group with outer automorphism $\tilde G$, where $H$ is the stable endoscopic group for
$\tilde G$. This fundamental lemma  describes a relationship between ordinary stable orbital integrals on the
endoscopic group $H$ and $\Theta$-twisted orbital integrals on $\tilde G$.
We will discuss the following lifts in detail:\klabsatz
$$\begin{array} {ll}
H & \tilde G\\
\SP_{2n} & \PGL_{2n+1}\rtimes\langle g\mapsto J\tran g^{-1} J^{-1}\rangle \\
\GSPIN_{2n+1}\quad &
(\GL_{2n}\times \GG_m)\rtimes\langle (g,a)\mapsto (J\tran g^{-1} J^{-1},\det g\cdot a)\rangle\\
\SP_{2n} & \widetilde{SO_{2n+2}}\simeq \ORTH_{2n+2}.
\end{array}$$
In each case we will reduce the fundamental lemma using the Kazhdan trick and a lot of observations in
linear algebra to a statement
which we call the $BC$-conjecture and which seems to be proven only for $n=1,2$:
\absatz
{\bf Conjecture:}{\em  If the regular topologically unipotent and algebraically semisimple
elements $ u\in \SO_{2n+1}(\OOO_F)$ and
 $v\in \Sp_{2n}(\OOO_F)$ are  $BC$-matching (see \ref{BCMatching}) then}
\begin{gather}\tag{$\text{BC}_n$}
    \STO_u(1,\SO_{2n+1})=\STO_v(1,\Sp_{2n}).
\end{gather}
\absatz
Thus the ($\Theta$-twisted) fundamental lemmas for the three series of endoscopy will be reduced
to a fundamental lemma like statement for ordinary (i.e. untwisted) stable orbital integrals on
the groups $\SO_{2n+1}$ of type $B_n$ resp. $\SP_{2n}$ of type $C_n$. An outline of the proofs will
be given in chapter \ref{KAPITEL Anfang}.

\thispagestyle{empty}     \Seitenumbruch

\pagenumbering{arabic} \setcounter{page}{1}
\meinchapter{Stable endoscopy and matching}
\Numerierung {\bf Notations.} In this paper we will denote by $F$ a $p$-adic field with ring of integers $\OOO_F$,
prime ideal $\pPp$ and uniformizing element $\varpi=\varpi_F$. The residue
field of characteristic $p$ is denoted $\kappa=\kappa_F=\OOO_F/\pPp$. By $\bar F$ we denote an
algebraic closure of $F$.
In the whole paper we will assume that $p\ne 2$. Only in this chapter $R$ denotes an arbitrary
integral domain.
\Absatz
\Numerierung {\bf Split Groups with automorphism.}\label{splitmitauto}
Let $G/R$ be a connected reductive split group scheme. We fix some "splitting" i.e.
a tripel $(B,T,\{X_\alpha\}_{\alpha\in \Delta})$ where $T$ denotes a maximal split-torus
inside a rational Borel $B$, $\Delta=\Delta_G=\Delta(G,B,T)\subset \Phi(G,T)\subset X^*(T)$ the set
 of simple roots inside the system of roots and
the $X_\alpha$ are a system (nailing)  of isomorphisms between the additive group scheme $\G_a$
and the unipotent
root subgroups $B_\alpha$. Here $X^*(T)=Hom(T,\GG_m)$ denotes the character module of $T$, while
$X_*(T)=Hom(\GG_m,T)$ will denote the cocharacter module of $T$.
Let  $\Theta\in Aut(G)$ be an automorphism of $G$ which fixes the
splitting, i.e. stabilizes $B$ and $T$ and permutes the $X_\alpha$. We  assume $\Theta$
to be of finite order $l$. We denote by
$$\tilde G=G\rtimes \langle \Theta\rangle$$
the (nonconnected) semidirect product of $G$ with $\Theta$.
 $\Theta$ acts on the (co)character module via
$X_*(T)\ni \alpha^\vee \mapsto \Theta\circ\alpha^\vee$ resp. $X^*(T)\ni \alpha\mapsto \alpha\circ
\Theta^{-1}$.
\Absatz
\Numerierung{\bf The dual group.}
Let $\hat G=\hat G(\C)$ be the dual group of $G$. By definition $\hat G$ has a tripel
$(\hat B,\hat T,\{\hat X_{\hat\alpha}\})$ such that we have identifications
$X^*(\hat T)=X_*(T), \enskip X_*(\hat T)= X^*(T)$ which identifies the (simple) roots
$\hat\alpha\in X^*(\hat T)$ with the (simple) coroots $a^\vee \in X_*(T)$, and the (simple) coroots
$\hat\alpha^\vee\in X_*(\hat T)$ with the (simple) roots $\alpha \in X^*(T)$. There exists a
unique automorphism $\hat\Theta$ of $\hat G$ which stabilizes $(\hat B,\hat T,\{\hat X_{\hat\alpha}\})$
and induces on $(X_*(\hat T),X^*(\hat T))$ the same automorphism as $\Theta$ on $(X^*(T),X_*(T))$.
\Absatz
\Numerierung{\bf The $\Theta$-invariant subgroup in $\hat G$.}
Let $\hat H=(\hat G^{\hat\Theta})^\circ$ be the connected component of the subgroup of $\hat\Theta$-fixed
elements in $\hat G$. It is a reductive split group with triple
$(\hat B_H,\hat T_H,\{X_\beta\}_{\beta\in\Delta_{\hat H}})$, where $\hat B_H=\hat B^{\hat\Theta}$,
$\hat T_H=\hat T^{\hat\Theta}$ and the nailing will be explained soon. We have the inclusion of
cocharacter modules
$X_*(\hat T_H)   
=X_*(\hat T)^{\hat\Theta}\subset X_*(\hat T)$ and a projection for the
character module $$P_\Theta:X^*(\hat T)\twoheadrightarrow (X^*(\hat T)_{\hat\Theta})_{free}=X^*(\hat T_H),$$
where $(X^*(\hat T)_{\hat\Theta})_{free}$ denotes the maximal free quotient of the coinvariant
module $X^*(\hat T)_{\hat\Theta}$.
For a $\Z[\Theta]$-module $X$ we define a map
\begin{gather*} S_\Theta:\quad X\to X^\Theta,\qquad x\mapsto \sum_{i=0}^{ord_x(\Theta)-1}\Theta^{i}(x)
\end{gather*}
where $ord_x(\Theta)=\min\{i>0\mid \Theta^{i}(x)=0\} $ is length of the orbit $\langle\Theta\rangle(x)$,
 which may vary on $X$.
\klabsatz
 For the roots $\Phi$ and coroots $\Phi^\vee$ of a given root datum
 $(X^*,X_*,\Phi,\Phi^\vee)$ we have to introduce a modified map $S'_\Theta$ by
 \begin{align*}
 S'_\Theta(\alpha)\quad &= \quad c(\alpha)\cdot S_\Theta(\alpha)\qquad\text{ where}\\
 c(\alpha)\quad &=\quad\frac{2}{\langle\alpha^\vee, S_\Theta(\alpha)\rangle}
 \end{align*}
resp. by the formula where the roles of $\alpha$ and $\alpha^\vee$ are exchanged. For all simple root systems
with automorphisms which are not of type $A_{2n}$ we have $\langle \alpha^\vee,\Theta^{i}(\alpha)\rangle=0$
for $i=1,\ldots, ord_\alpha(\Theta)-1$ which implies $c(\alpha)=1$ i.e.
$S'_\Theta(\alpha)=  S_\Theta(\alpha)$.
We furthermore introduce the subset of short-middle roots and the dual concept of long-middle coroots:
\begin{align*} \Phi(\hat G, \hat T)^{sm}\quad &=\quad
    \left\{\alpha\in \Phi(\hat G,\hat T)\left\vert \frac{1}{2}\cdot
        P_\Theta(\alpha)\notin  P_\Theta(\Phi(\hat G,\hat T))\right\}\right.\\
         \Phi^\vee(\hat G,\hat T)^{lm}\quad &= \quad \left\{ \alpha^\vee\left\vert \alpha\in
          \Phi(\hat G, \hat T)^{sm}\right\}\right.
\end{align*}
\begin{prop} \label{CoWurzelbeschreibungbeiInvarianten}
With the above notations we have
\begin{align} \label{PThetaGleichung}
\Phi(\hat H,\hat T_H)\quad &= \quad P_\Theta(\Phi(\hat G, \hat T)^{sm}) \qquad\qquad\qquad\text{for the roots}
\\ \label{SThetaGleichung}
  \Phi^\vee(\hat H,\hat T_H)\quad &= \quad S'_\Theta( \Phi^\vee(\hat G,\hat T)^{lm})
   \qquad\qquad\quad \text{ for the coroots}
\\  \notag \Delta^\vee_{\hat H}=\Delta^\vee(\hat H,\hat B_H,\hat T_H)\quad &=
        \quad S'_\Theta(\Delta^\vee_{\hat G})
       \qquad\qquad\quad\;\text{for the simple coroots}
\\ \notag \Delta_{\hat H}=\Delta(\hat H,\hat B_H,\hat T_H)\quad &= \quad P_\Theta(\Delta_{\hat G})
        \qquad\qquad\qquad\text{for the simple roots}
\end{align}
\end{prop}
Proof: This follows from \cite[8.1]{St}, which is restated in Theorem
\ref{Steinberg: G_hoch_tTheta ist reduktiv} below.
\qed
\begin{Def}[stable $\Theta$-endoscopic group]
In the above situation a connected reductive split group scheme $H/R$ will be called
a stable $\Theta$-endoscopic group for $(G,\Theta)$ resp. $\tilde G$  if its dual group is together with the
splitting isomorphic to the above $(\hat H,\hat B_H,\hat T_H,\{X_\beta\}_{\beta\in\Delta_{\hat H}})$.
\end{Def}
Remarks: Since $H$ is unique up to
 isomorphism (up to unique isomorphism if we consider $H$ together with a splitting)
 we can call $H$ {the} stable $\Theta$-endoscopic group
for $(G,\Theta)$. For a maximal split-torus $T_H\subset H$ we have:
\begin{align}\label{RelationX_*}
 X_*(T_H)\quad &= \quad X_*(T)_\Theta \qquad\text{ for the cocharacter-module}\\ \notag
 \label{RelationX^*} X^*(T_H)\quad &= \quad X^*(T)^\Theta \qquad\text{ for the character-module}
\end{align}
\Absatz\Numerierung
  To get examples we use the following {\bf notations}:
  \klabsatz\quad $diag(a_1,\ldots,a_n)\in \GL_n$ denotes the diagonal
  matrix $(\delta_{i,j}\cdot a_i)_{ij}$ and
  \klabsatz\quad $antidiag(a_1,\ldots,a_n)\in \GL_n$ the antidiagonal matrix
  $(\delta_{i,n+1-j}\cdot a_i)_{ij}$ with $a_1$ in the upper right corner.
  We introduce the following matrix
  \begin{gather*}
    J=J_n=(\delta_{i,n+1-j}(-1)^{i-1})_{1\le i,j\le n}=antidiag(1,-1,\ldots,(-1)^{n-1})\in \GL_n(R).
  \end{gather*}
  and its modification $J'_{2n}=antidiag(1,-1,1,\ldots,(-1)^{n-1},(-1)^{n-1},\ldots,1,-1,1)$.
  Since $\tran J_n=(-1)^{n-1}\cdot J_n$  and $J'_{2n}$ is symmetric we can define the
  \begin{align*} \text{standard symplectic group} \quad\qquad\Sp_{2n}\quad &= \quad \Sp(J_{2n})\\
    \text{standard split odd orthogonal  group}\qquad \SO_{2n+1}\quad &= \quad \SO(J_{2n+1}).\\
    \text{standard split even orthogonal  group}\quad\qquad \SO_{2n}\quad &= \quad \SO(J'_{2n}).
  \end{align*}
  We consider the groups $\GL_n,\SL_n,\PGL_n,\SP_{2n},\SO_n$ with the splittings consisting
  of the diagonal torus, the Borel consisting of upper triangular matrices and the standard nailing.
  We remark that the following map defines an involution of $\GL_n,\SL_n$ and $\PGL_n$:
  \begin{gather*} \Theta=\Theta_n: g\mapsto J_n\cdot \tran g^{-1}\cdot J_n^{-1}.
  \end{gather*}
\begin{expl}[$A_{2n}\leftrightarrow C_n$]\label{BspPGL2n+1Sp2n}
   \em
   \begin{gather*}
   \begin{matrix}
   \qquad G=\PGL_{2n+1}, &\Theta=\Theta_{2n+1}&\quad\text{ has dual }\quad&
   \hat G=\SL_{2n+1}(\C),\,& \hat\Theta=\Theta_{2n+1}\\
   &&&\bigcup&\\
   \qquad\quad H=\SP_{2n}&&\quad\text{ has dual }\quad&\hat H=\SO_{2n+1}(\C)&\\
   \end{matrix}
   \end{gather*}
\end{expl}
\Absatz
\begin{expl}[$A_{2n-1}\leftrightarrow B_n$]\label{BspGL2nGl1GSpin2n+1}
  \em
  The group $G=\GL_{2n}\times \G_m$ has the automorphism
  \begin{gather*} \Theta:(g,a)\mapsto (\Theta_{2n}(g),\det(g)\cdot a)
  \end{gather*}
  which is an involution since $\det(\Theta_{2n}(g))=\det g^{-1}$.
  The dual $\hat \Theta \in Aut(\hat G)$ satisfies
  \begin{align*}
     \hat\Theta(g,b)&=(\Theta_{2n}(g)\cdot b, b), \text{\qquad so that we get}
  \end{align*}
  \begin{align*}
     &\begin{matrix}
         \qquad\qquad G=\GL_{2n}\times \G_m, &\Theta\quad&\quad\text{ has dual }\qquad&
         \hat G=\GL_{2n}(\C)\times \C^\times,\,& \hat\Theta
       \\    &&&\bigcup&\\
         \qquad\qquad\quad H=\GSPIN_{2n+1}&&\quad\text{ has dual }\qquad&\hat H=\GSp_{2n}(\C).&
       \\
      \end{matrix}
  \end{align*}
  Recall that $\GSPIN_{2n+1}$ can be realized as the quotient
  $\left(\G_m\times \Spin_{2n+1}\right)/\mu_2$, where $\mu_2\simeq\{\pm 1\}$
  is embedded diagonally, so that we get an exact sequence
  $$1\;\rightarrow\; \Spin_{2n+1}\;\rightarrow\;\ \GSPIN_{2n+1}\;\xrightarrow{ \mu }\; \G_m\;\rightarrow \;1,$$
  where the "multiplier" map $\mu$ is induced by the projection to the $\G_m$ factor followed by squaring.
   Thus the derived group of $\GSPIN_{2n+1}$ is $\Spin_{2n+1}$, i.e.
  a connected, split and simply connected group.
\end{expl}
\Absatz
\begin{expl}[$D_{n+1}\leftrightarrow C_n$]\label{BspSO2n+2Sp2n}
  \em We furthermore have the situation:
  \begin{gather*}
    \begin{matrix}
      \qquad G=\SO_{2n+2}, &\Theta=int(s)&\quad\text{ has dual }\quad&
      \hat G=\SO_{2n+2}(\C),\,& \hat\Theta=int(\hat s)\\
      &&&\bigcup&\\
      \qquad\quad H=\SP_{2n}&&\quad\text{ has dual }\quad&\hat H=\SO_{2n+1}(\C)&\\
    \end{matrix}
  \end{gather*}
  where $s\in \ORTH_{2n+2}$ denotes the reflection which interchanges the standard basis
  vectors $e_{n+1}$ and $e_{n+2}$ and fixes all other basis elements $e_i$. The
  $\hat s$ is of the same shape.
\end{expl}
\Absatz
\Numerierung {\bf Matching elements}\label{DefMatching}
  Since each semisimple $\Theta$-conjugacy resp. conjugacy class in $G(\bar F)$ resp. $H(\bar F)$
  meets $T(\bar F)$ resp. $T_H(\bar F)$ we have canonical bijections
  \begin{align*} i_G:\quad G(\bar F)_{ss}/\Theta-conj \qquad &\simeq\quad T(\bar F)_\Theta /(W_G)^\Theta\\
    i_H:\quad H(\bar F)_{ss}/conj\qquad\quad &\simeq\quad T_H(\bar F)/W_H
  \end{align*} where $W_G=Norm(T,G)/T$ and $W_H=Norm(T_H,H)/T_H$ denote the Weyl-groups.
  We further have an isomorphism
  \begin{align*} N_{KS}:\quad T(\bar F)_\Theta\simeq(X_*(T)\otimes {\bar F}^\times)_\Theta
    \simeq X_*(T)_\Theta\otimes{\bar F}^\times=X_*(T_H)\otimes {\bar F}^\times \simeq T_H(\bar F)
  \end{align*} and observe $W_H\simeq (W_G)^\Theta$. Therefore we may define:\klabsatz
  Two ($\Theta$-)semisimple elements
  $\gamma\Theta\in G(F)\Theta$ and $h\in H(F)$ are called {\it matching} if their ($\Theta$-)stable conjugacy
  classes in $G(\bar F)$ resp. $H(\bar F)$ correspond to each other via the isomorphism
  $i_H^{-1}\circ N_{KS}\circ i_G$.
\Absatz\Numerierung\label{BCMatching}{\bf $BC$-matching}:
  We have an isomorphism between the standard diagonal tori:
  \begin{align*}
    i_{BC}: T_{\SO_{2n+1}}&\to T_{\Sp_{2n}}, \quad
    \\
    diag(t_1,\ldots,t_n,1,t_n^{-1},\ldots,t_1^{-1})
    &\mapsto diag(t_1,\ldots,t_n,t_n^{-1},\ldots,t_1^{-1}).
  \end{align*}
  We observe that $i_{BC}$ induces an isomorphism of Weylgroups:
  \begin{gather*}
    W_{\SO_{2n+1}}\quad\simeq\quad S_n\ltimes \{\pm1\}^n\quad\simeq W_{\Sp_{2n}}
  \end{gather*}
  Two semisimple elements $h\in \SO_{2n+1}(F)$ and $b\in \Sp_{2n}(F)$ are called {\it $BC$-matching}
  if their stable conjugacy classes correspond to each other under the isomorphism
  $i_{\Sp_{2n}}^{-1}\circ i_{BC} \circ i_{\SO_{2n+1}}$.
\Absatz
\begin{expl} \label{BspMatchingPGL2n+1Sp2n}
  \em
  In example \ref{BspPGL2n+1Sp2n} above the norm map $N_{KS}: T\to T_\Theta \simeq T_H$ is given by
  \begin{align} \label{NormmapPGL2n+1Sp2n}
    \gamma\quad=\quad diag(t_1,t_2,\ldots,t_n,t_{n+1},t_{n+2},\ldots,t_{2n+1})
    \qquad &\in T\subset\PGL_{2n+1}\\
    \notag\mapsto h=diag(t_1/t_{2n+1},t_2/t_{2n},\ldots,t_n/t_{n+2},t_{n+2}/t_n,\ldots,t_{2n+1}/t_1)&\in
    T_H\subset\Sp_{2n}.
  \end{align}
\end{expl}
Proof: We identify $X_*(T)\simeq \Z^{2n+1}/\Z \simeq \bigoplus_{i=1}^{2n+1} \Z f_i/\Z\sum_{i=1}^{2n+1} f_i$
  and $X^*(T)\simeq\{\sum_{i=1}^{2n+1} \alpha_i e_i\mid \sum_{i=1}^{2n+1} \alpha_i=0\}$, such that
  \begin{align*}
    f_i:\quad t\;\mapsto\; & diag(1,\ldots,1,\underset{i}{t},1,\ldots,1)\in T\qquad\text{ and }\\
    e_i-e_j:\quad T\;\ni\; & diag(t_1,\ldots,t_{2n+1})\;\mapsto\; t_i/t_j.
  \end{align*}
  Similarly we identify $T_H\simeq \GG_m^n$ via $diag(t_1,\ldots,t_n,t_n^{-1},\ldots,t_1^{-1})\mapsto
  (t_i)_{1\le i\le n}$ and write $X_*(T_H)\simeq X_*(\GG_m^n)\simeq \bigoplus_{i=1}^n \Z f'_i$ resp.
  $X^*(T_H)\simeq X^*(T_H)\simeq \bigoplus_{i=1}^n \Z e'_i$. The involution $\Theta$ acts as
  \begin{align*} \Theta(f_i)=-f_{2n+2-i},\qquad \Theta(e_i-e_j)=e_{2n+2-j}-e_{2n+2-i}
  \end{align*}
  Now it is clear that we have an identification $P_\Theta: X_*(T)_\Theta \simeq X_*(T_H)$ given by
  $P_\Theta(f_i)=f'_i$ and $P_\Theta(f_{2n+2-i})=-f'_i$ for $1\le i\le n$. This forces $P_\Theta(f_{n+1})
  =P_\Theta(-\sum_{i\ne n+1}f_i)=P_\Theta(\Theta(\sum_{i=1}^n f_i)-\sum_{i=1}^n f_i)=0$. Dual to this
  we have an injection $\iota: X^*(T_H)\simeq X^*(T)^\Theta\subset X^*(T)$ such that
  $\iota(e'_i)=e_i-e_{2n+2-i}$. It is furthermore clear that this $P_\Theta$ induces the map
  (\ref{NormmapPGL2n+1Sp2n})
  The claim now follows if show that $P_\Theta$ and $\iota$ correspond to the natural maps
  on the side of the dual groups which are characterized by the equations of Proposition
  \ref{CoWurzelbeschreibungbeiInvarianten}. Especially we have to check the duals of the relations
  (\ref{PThetaGleichung}) and (\ref{SThetaGleichung}) for our explicitly given maps
  $P_\Theta$ and $\iota$, namely the equations
  \begin{align} \label{PThetaGleichungendo}
       \Phi^\vee( H, T_H)\quad &= \quad  P_\Theta(\Phi^\vee( G, T)^{sm})\qquad\qquad\text{for the coroots}
     \\ \label{SThetaGleichungendo}
       \iota(\Phi(H, T_H))\quad &= \quad S'_\Theta( \Phi(G,T)^{lm})
       \qquad\qquad\quad \text{ for the roots.}
  \end{align}
 \absatz
  But we have
  \begin{align}\label{Cowurzelgleichung1}
    P_\Theta(\pm(f_i-f_j))\quad &= \quad \pm(f'_i-f'_j)\quad=\quad P_\Theta(\pm(f_{2n+2-j}-f_{2n+2-i})),\\
    \label{Cowurzelgleichung2}
    P_\Theta(\pm(f_i-f_{2n+2-j}))\quad &=\quad \pm(f'_i+f'_j), \\
    \notag
    P_\Theta(\pm(f_i-f_{n+1}))
    \quad &=\quad \pm f'_i \quad =\quad P_\Theta(\pm(f_{n+1}-f_{2n+2-i}))
  \end{align}
  where $1\le i,j \le n$ in all three equations, but where additionally $i\ne j$ in (\ref{Cowurzelgleichung1})
  while $i=j$ is allowed in (\ref{Cowurzelgleichung2}). Nevertheless $P_\Theta(\pm(f_i-f_{2n+2-i}))=2
  \cdot f'_i$ is not a member of the right hand side of (\ref{PThetaGleichungendo}) since
  $f'_i=P_\Theta(f_i-f_{n+1})$. By the well known description of $\Phi^\vee(\Sp_{2n},T_H)$ we
  get the equality (\ref{PThetaGleichungendo}).
 \klabsatz
  Similarly we get
  \begin{align*}
     S_\Theta(\pm(e_i-e_j))\quad &= \quad \iota(\pm(e'_i-e'_j))\quad=\quad S_\Theta(\pm(e_{2n+2-j}-e_{2n+2-i})),\\
     S_\Theta(\pm(e_i-e_{2n+2-j})) \quad &=\quad \iota(\pm(e'_i+e'_j)), \\
     S_\Theta(\pm(e_i-e_{n+1}))\quad &=\quad \pm(e_i-e_{2n+2-i})
     \quad=\quad S_\Theta(\pm(e_{n+1}-e_{2n+2-i}))\\
     \quad &=\quad\pm\iota(e'_i)\quad=\quad S_\Theta(\pm(e_i-e_{2n+2-i}))
  \end{align*}
  where $1\le i\ne j\le n$ in the first two equations and $1\le i \le n$ in the last two. Since
  $\Phi(\Sp_{2n},T_H)=\{e'_i-e'_j\mid 1\le i\ne j\le n\}\cup \{e'_i+e'_j\mid 1\le i\ne j\le n\}\cup
  \{2\cdot e'_i\mid 1\le i\le n\}$ and since $\langle f_i-f_{n+1},S_\Theta(e_i-e_{n+1})\rangle=1$
  but $\langle f_i-f_j,S_\Theta(e_i-e_j)\rangle=2$ for all $i,j\ne n+1, 1\le i,j \le 2n+1$
  we get the claim
  (\ref{SThetaGleichungendo}) in view of the definition of $S'_{\Theta}$.
 \klabsatz
  Finally it is clear, that $S'_\Theta$ maps the set of simple roots
  $\{e_i-e_{i+1}, e_n-e_{n+1},e_{n+1}-e_{n+2},e_{2n+1-i}-e_{2n+2-i}\mid 1\le i\le n-1\}$ of $G$ to
  the set of simple roots $\iota(\{e'_i-e'_{i+1},2\cdot e'_n\mid 1\le i\le n-1\})$ of $H$
  and that $P_\Theta$ maps the set of simple coroots
  $\{f_i-f_{i+1}, f_n-f_{n+1},f_{n+1}-f_{n+2},f_{2n+1-i}-f_{2n+2-i}\mid 1\le i\le n-1\}$
  of $G$ to the set of simple coroots $\{f'_i-f'_{i+1},f'_n\mid 1\le i\le n-1\}$.
\qed

\begin{expl}
  \label{BspMatchingGL2nGSpin2n+1}
  \em
  In example \ref{BspGL2nGl1GSpin2n+1} above  we consider additionally the projection
  $pr_{ad}: \GSPIN_{2n+1}\to \SO_{2n+1}=\Spin_{2n+1}/\{\pm 1\}$. Then the composite map
  $pr_{ad}\circ N_{KS}: T\to T_{ad}\subset \SO_{2n+1}$ is given by
  \begin{align} \label{NormmapGL2nGSPIN2n+1}
    \gamma\quad=\quad (diag(t_1,t_2,\ldots,t_n,t_{n+1},\ldots,t_{2n}),t_0)
    \qquad &\in T\subset\GL_{2n}\times \GG_m\\
    \notag\mapsto h=diag(t_1/t_{2n},\ldots,t_n/t_{n+1},1,t_{n+1}/t_n,\ldots,t_{2n}/t_1)&\in
  T_{ad}\subset\SO_{2n+1}.
  \end{align}
\end{expl}
Proof: We consider the following basis $(e_i)_{0\le i\le 2n}$ of $X^*(T)$:
  \begin{align*}
    e_i:\quad T\;\ni\; &(diag(t_1,\ldots,t_{2n}),t_0)\quad\mapsto \quad t_i.
  \end{align*}
  Let $(f_i)_{0\le i\le 2n}$ be the dual basis of $X_*(T)$.
  We furthermore identify $T_{ad}\simeq \GG_m^n$ via $diag(t_1,\ldots,t_n,1,t_n^{-1},\ldots,t_1^{-1})\mapsto
  (t_i)_{1\le i\le n}$ and write $X_*(T_{ad})\simeq X_*(\GG_m^n)\simeq \bigoplus_{i=1}^n \Z \tilde f_i$ resp.
  $X^*(T_H)\simeq X^*(T_H)\simeq \bigoplus_{i=1}^n \Z \tilde e_i$. The involution $\Theta$ acts via
  \begin{align*} \Theta(e_i)\quad &=\quad -e_{2n+1-i},\qquad\quad \Theta(f_i)\quad =\quad -f_{2n+1-i}+f_0\qquad
    \text{for }1\le i\le 2n\\
    \Theta(e_0)\quad&=\quad e_0 +\sum_{i=1}^{2n}e_i,\qquad\quad \Theta(f_0)\quad=\quad f_0
  \end{align*}
  Now it is clear that $X^*(T_H)=X^*(T)^\Theta$ has as basis $ (e'_0,e'_1,\ldots, e'_n)$ where
  $e'_i=e_i-e_{2n+1-i}$ for $1\le i\le n$ and $e'_0=e_0+\sum_{i=n+1}^{2n} e_i$.
  Let $(f'_0,f'_1,\ldots,f'_n)$ be the dual basis of $X_*(T_H)=X_*(T)_\Theta$.
  Then the projection map $P_\Theta:X_*(T)\to X_*(T)_\Theta$ satisfies:
  $P_\Theta(f_0)=f'_0,\quad P_\Theta(f_i)=f'_i,\quad P_\Theta(f_{2n+1-i})= -f'_i + f'_0$ for $1\le i\le n$.
 \klabsatz
  For the (co)root systems we get:
  \begin{align*} \Phi(\GSPIN_{2n+1}&,T_H)\quad = \quad P_\Theta(\Phi(\GL_{2n}\times \GG_m,T))\\
    = \quad &\{\pm e'_i \pm e'_j\mid 1\le i<j\le n\}\cup \{e'_i\mid 1\le i\le n\}\\
    \Phi^\vee(\GSPIN_{2n+1}&,T_H)\quad =\quad S'_\Theta(\Phi^\vee(\GL_{2n}\times \GG_m,T)),\\
    = \quad &\{\pm(f'_i-f'_j)\mid 1\le i<j\le n\}\cup
    \{\pm(f'_i+f'_j-f'_0)\mid 1\le i\le j\le n\}.
  \end{align*}
  The cocharacter group of the center of $\GSPIN_{2n+1}$ is then recognized as $\Z f'_0$. Thus we may
  define a surjection $pr_{ad}:X_*(T_H)\to X_*(T_{ad})$ by $f'_0\mapsto 0$ and $f'_i\mapsto \tilde f_i$
  for $1\le i\le n$. Dually one has the injection $\iota:X^*(T_{ad})\hookrightarrow X^*(T_H), \tilde e_i
  \mapsto e'_i$. Now it is clear that  $pr_{ad}\circ P_\Theta$ induces the map
  (\ref{NormmapGL2nGSPIN2n+1}) and it remains to show that we have the following relations analogous
  to (\ref{PThetaGleichungendo}) and (\ref{SThetaGleichungendo}):
  \begin{align} \label{prPThetaGleichungendo}
    \Phi^\vee( \SO_{2n+1}, T_{ad})\quad &= \quad  pr_{ad}\circ P_\Theta(\Phi^\vee( G, T)^{sm})\qquad\qquad\text{for the coroots}
    \\ \label{prSThetaGleichungendo}
    \iota(\Phi(\SO_{2n+1}, T_{ad}))\quad &= \quad S'_\Theta( \Phi(G,T)^{lm})
    \qquad\qquad\qquad\quad \text{ for the roots.}
  \end{align}
  \absatz
  But this follows immediately from the above description of $P_\Theta(\Phi(\GL_{2n}\times \GG_m,T))$ and
  $S'_\Theta(\Phi^\vee(\GL_{2n}\times \GG_m,T))$ in view of the very simple shape of $pr_{ad}$ and $\iota$
  and the knowledge of the (co)root system of $\SO_{2n+1}$.
  The relation for the bases of the (co)root systems is checked in a similar way.
\qed
\begin{expl}\label{BspMatchingmurelation}
  \em
   In example \ref{BspGL2nGl1GSpin2n+1} above we now analyze the relation between the multiplier map
   $\mu$ and matching. We claim:
   If $(h,a)\in \GL_{2n}(F)\times F^\times$ and $\eta\in \GSPIN_{2n+1}(F)$ match then we have:
   \begin{align*}
     \mu(\eta)\quad=\quad \det(h)\cdot a^{2}.
   \end{align*}
\end{expl}
Proof: In the notations of \ref{BspMatchingGL2nGSpin2n+1} the element $e'=2e'_0+\sum_{i=1}^n e'_i=
   2e_0+\sum_{i=1}^{2n} e_i\in X^*(T_H)=X^*(T)^\Theta$ corresponds to the character $(h,a)\mapsto \det(h)\cdot a^2$.
   Since $e'$ is orthogonal to the coroots $\Phi^\vee(\GSPIN_{2n+1},T_H)$ it has to correspond to a multiple of
   the multiplier $\mu$. Now it is easy to see that it corresponds in fact to $\mu$.
\qed
\begin{expl}
  \label{BspMatchingSO2n+2Sp2n}
  \em
  In example \ref{BspSO2n+2Sp2n} above the norm map $N_{KS}: T\to T_\Theta \simeq T_H$ is given by
  \begin{align} \label{NormmapSO2n+2Sp2n}
    \gamma\quad=\quad diag(t_1,t_2,\ldots,t_n,t_{n+1},t_{n+1}^{-1},\ldots,t_1^{-1})
    \qquad &\in T\subset\SO_{2n+2}\\
    \notag\mapsto h=diag(t_1,\ldots,t_n,t_n^{-1},\ldots,t_1^{-1})&\in
    T_H\subset\Sp_{2n}.
  \end{align}
\end{expl}

Proof: We consider the following basis $(e_i)_{1\le i\le n+1}$ of $X^*(T)$:
  \begin{align*}
    e_i:\quad T\ni & diag(t_1,\ldots,t_{n+1},t_{n+1}^{-1},\ldots,t_1^{-1})\quad\mapsto \quad t_i.
  \end{align*}
  Let $(f_i)_{1\le i\le n+1}$ be the dual basis of $X_*(T)$.
  The involution $\Theta$ acts via
  \begin{align*} \Theta(e_i)\quad &=\quad e_i,\qquad\qquad\quad \Theta(f_i)\quad =\quad f_i\qquad
    \text{for }1\le i\le n\\
    \Theta(e_{n+1})\quad&=\quad -e_{n+1},\qquad \Theta(f_{n+1})\quad=\quad -f_{n+1}.
  \end{align*}
  We furthermore use the bases $(e'_i)_{1\le i\le n}$ of $X^*(T_H)$ and $(f'_i)_{1\le i\le n}$ of
  $X_*(T_H)$ from example
  \ref{BspMatchingPGL2n+1Sp2n}.
  It is clear that we have an isomorphism
  $\iota:X^*(T_H)\simeq X^*(T)^\Theta$ given by $e'_i\mapsto e_i$ and a dual isomorphism
  $P_\Theta: (X_*(T)_\Theta)_{free}\simeq X_*(T_H)$ induced by the dual map $P_\Theta: f_{n+1}\mapsto 0$
  and $P_\Theta(f_i)\mapsto f'_i$ for $1\le i\le n$. It is clear that this $P_\Theta$ induces the
  map (\ref{NormmapSO2n+2Sp2n}).
 \klabsatz
  Recall the (co)root systems of $\SO_{2n+2}$
  \begin{align*} \Phi(\SO_{2n+2},T)\quad &= \quad \{\pm e_i \pm e_j\mid 1\le i<j\le n+1\}\\
    \Phi^\vee(\SO_{2n+2},T)\quad &= \quad \{\pm f_i\pm f_j\mid 1\le i<j\le n+1\}.
  \end{align*}
  We have
  $P_\Theta(f_i\pm f_j)= f'_i\pm f'_j,\quad S'_\Theta(e_i\pm e_j)=\iota(e'_i\pm e'_j)$ for $1\le i<j\le n$
  and $P_\Theta(f_i\pm f_{n+1})=f'_i,\quad S'_\Theta(e_i\pm e_{n+1})=e_i\pm e_{n+1}+\Theta(e_i\pm e_{n+1})=
  \iota(2\cdot e'_i)$ for $1\le i\le n$. The relations
  (\ref{PThetaGleichungendo}) and (\ref{SThetaGleichungendo}) now follow from the knowledge
  of the (co)root system of $\SP_{2n}$.
  It is clear that $S'_\Theta$ maps the simple roots $e_i-e_{i+1}, e_n\pm e_{n+1}$ of $\SO_{2n+2}$
  to the simple roots $e'_i-e'_{i+1},2\cdot e'_n$ of $\SP_{2n}$ and $P_\Theta$ maps the simple
  coroots $f_i-f_{i+1},f_n\pm f_{n+1}$ of $\SO_{2n+2}$ to the simple coroots $f'_i-f'_{i+1},f'_n$
  of $\SP_{2n}$.
\qed

     \Seitenumbruch

\meinchapter{Centralizers}\label{KAPITEL Anfang}

\begin{thm}[Steinberg]\label{Steinberg: G_hoch_tTheta ist reduktiv}
  Let $T$ be a $\Theta$--stable, maximal subtorus of $G$ and $t\in T(\bar F)$.
  The group $G^{t\Theta}$ of fixed point in $G$ under $\text{\em int}(t)\circ\Theta$ is reductive.
  The root system of the connected component $(G^{t\Theta})^\circ$ of $1$,
  viewed as a subsystem of $\Phi_\Theta=P_\Theta(\Phi(T,G))$,
  which might be identified with $\Phi_{\text{\em res}}=\{\alpha|_{T^{\Theta\circ}}\mid\alpha\in\Phi\}$,
  is given by
  \[\Phi(G^{t\Theta\circ},T^{\Theta\circ})
    \simeq\left\{P_\Theta(\alpha)\in \Phi_\Theta\ \Big|\  S_\Theta(\alpha)(t)=\left\{
    \begin{array}{rl}
      1 & \text{if }\frac{1}{2}P_\Theta(\alpha)\not\in\Phi_\Theta\\
     -1 & \text{if }\frac{1}{2}P_\Theta(\alpha)\in\Phi_\Theta
    \end{array}
   \right.\right\}.
  \]
\end{thm}
Proof: \cite[8.1]{St}
\qed
\Absatz\Numerierung\label{sTheta ist speziell auf NormCent}
  Define $N_\Theta(t)=\prod_{i=0}^{\text{ord}(\Theta)-1} \Theta^{i}(t)$.
  In case $\Phi(G)$ is of type $A_{2n}$ or one $\Theta$--Orbit of components of type $A_{2n}$
  the root system $\Phi(G^{t\Theta\circ},T^{\Theta\circ})$ in \ref{Steinberg: G_hoch_tTheta ist reduktiv}
  is a maximal reduced subsystem of
  \[\{P_\Theta(\alpha)\ \mid \ \alpha(N_\Theta(t))=1\}.
  \]
\Absatz\Numerierung\label{Extended Dynkin-diagram}
  For an irreducible root system $\Phi$ with basis $\Delta$ and $\theta\in \text{Aut}(\Phi,\Delta)$ we
  denote by $\tilde \alpha\in\Phi^-$ the negative root such that $-c(\tilde\alpha)S_\Theta(\tilde\alpha)
  =-S'_\Theta(\tilde\alpha)$ is the highest root in $S'_\Theta(\Phi^+)$ with respect to the basis
  $S'_\Theta(\Delta)$.
  \klabsatz
  To $(\Phi,\Delta,\Theta)$ we associate the extended Dynkin diagram
  \[ \begin{array}{rcl}
    \Delta_{ext}(\Phi,\Theta)
    & := & P_\Theta(\Delta)\cup \{c(\tilde\alpha)P_\Theta(\tilde\alpha)\}
   \end{array}
  \]

\begin{prop}[Dynkin]\label{Dynkins Satz}
  Let $G$ and $\Theta$ be as in \ref{splitmitauto}.
  \begin{itemize}
    \item[\bf (a)]
        For every $t\in T$ there exists $w\in W^{\Theta}$, such that $\Phi((G^{w(t)\Theta})^\circ)$
        has a basis in $\Delta_{ext}(\Phi,\Theta)$.
    \item[\bf (b)]
        Every proper subsystem of $\Delta_{ext}(\Phi,\Theta)$
        occurs as $\Phi\left((G^{t\Theta})^\circ\right)$ for some $t\in T$.
  \end{itemize}
\end{prop}
Proof: For a detailed proof we refer to \cite{Balldiss}[2.42]. We remark that in case $\Phi=A_{2n}$
  the extended Dynkin diagram $(P_\Theta(\Phi),\Delta_{ext}(\Phi,\Theta))$ looks like
  \begin{align*}
  \begin{xy}<8mm,0mm>:
  (.5,0)*{>},
  (6.5,0)*{>},
    (0,0)*{\bullet}\PATH ~={**\dir{=}}
  '(1,0)*{\circ} ~={**\dir{-}}
  '(2,0)*{\circ}
  '(3,0)*{\circ} ~={**\dir{.}}
  '(4,0)*{\circ} ~={**\dir{-}}
  '(5,0)*{\circ}
  '(6,0)*{\circ}~={**\dir{=}}
  '(7,0)*{\circ}
  \end{xy}
  \end{align*}
  Therefore (a) will follow in the case $G=\SL_{2n+1}$  from the fact, proven in lemma \ref{Zentralisatoren}
  that the groups $G^{t\Theta}$  are isomorphic to groups of the form
\[G^{t\Theta}=\left(
\begin{array}{c@{}c}
\begin{array}{c@{}c@{\,}c}
\cline{1-1}
\multicolumn{1}{|@{}c@{}|}{\colorbox{mittelgrau}{\parbox[c][5mm][c]{5mm}{\centerline{$\star$}}}}     &       &
\\
\cline{1-2}
        &\multicolumn{1}{|@{}c@{}|}{\colorbox{hellgrau}{$\times$}}       &
\\
\cline{2-2}     &       &\ddots
\\
\parbox[c][3mm][c]{3mm}{\ }       &       &
\\
\hline\hline
\parbox[c][3mm][c]{3mm}{\ }        &       &
\\
     &       & \vphantom{\ddots}
\\
     &       &
\\
\cline{1-1}
\multicolumn{1}{|@{}c@{}|}{\colorbox{mittelgrau2}{\parbox[c][5mm][c]{5mm}{\centerline{$\star$}}}} &       &
\\
\cline{1-1}
\end{array}
\begin{array}{c||c}
    &   \parbox[c][20mm][c]{2mm}{\;}
\\
\cline{1-2}
        \multicolumn{1}{|@{}c@{}}{\colorbox{dunkelgrau}{\parbox[c][3mm][c]{3mm}{\centerline{$*$}}}}
    &   \multicolumn{1}{@{}c@{}|}{\colorbox{dunkelgrau}{\parbox[c][3mm][c]{3mm}{\centerline{$*$}}}}
\\
        \multicolumn{1}{|@{}c@{}}{\colorbox{dunkelgrau}{\parbox[c][3mm][c]{3mm}{\centerline{$*$}}}}
    &   \multicolumn{1}{@{}c@{}|}{\colorbox{dunkelgrau}{\parbox[c][3mm][c]{3mm}{\centerline{$*$}}}}
\\
\cline{1-2}
    &   \parbox[c][20mm][c]{2mm}{\;}
\end{array}
\begin{array}{c@{\,}c@{}c}
\cline{3-3}
&       &\multicolumn{1}{|@{}c@{}|}{\colorbox{mittelgrau}{\parbox[c][5mm][c]{5mm}{\centerline{$\star$}}}}
\\
\cline{3-3}
\vphantom{\colorbox{hellgrau}{$\times$}}       &       &
\\
\vphantom{\ddots} &     &
\\
\parbox[c][3mm][c]{3mm}{\ }     &       &
\\
\hline\hline
\parbox[c][3mm][c]{3mm}{\ }     &       &
\\
\ddots &       &
\\
\cline{2-2}     &\multicolumn{1}{|@{}c@{}|}{\colorbox{hellgrau}{$\times$}}&
\\
\cline{2-3}
&       &\multicolumn{1}{|@{}c@{}|}{\colorbox{mittelgrau2}{\parbox[c][5mm][c]{5mm}{\centerline{$\star$}}}}
\\
\cline{3-3}
\end{array}
\end{array}
\right)
\begin{array}{l}
\simeq\text{Sp}(2m_{k})\times\text{Gl}(m_{k-1})\times\cdots\\
\phantom{\simeq}
\cdots\times\text{Gl}(m_1)\times\text{SO}(2m_0+1)
\end{array}
\]
 where $(m_k,m_{k-1},\ldots , m_1,m_0)$ runs through all partitions of $n$ with $m_k,m_0\ge 0$ and
$m_i\ge 1$ for $0<i<k$.
This means that the $G^{t\Theta}$ are of type
$C_{m_k}\times A_{m_{k-1}-1}\times\cdots\times A_{m_1-1}\times B_{m_0}$
and the $\Delta(G^{t\Theta})$ are subsets of $P_\Theta(\Delta)\cup\{-P_\Theta(\alpha^+)\}$,
where  $\alpha^+=\eps_1-\eps_{2n+1}$ is the highest root for $G$.
\qed

\Absatz\Numerierung\label{ErlaeuterungDynkindiagram}
In the following table we list all simple root systems $\Phi=\Phi(G)$, such that the
semisimple (simply connected split) group $G$ has an outer automorphism $\Theta$, together
with the root systems $\Phi(H)$ of the stable endoscopic groups $H$ of $(G,\theta)$.
The ordinary simple roots are
marked by a $\circ$, the additional root $\tilde\alpha$  by a $\bullet$. We get
six blocks, separated by double lines, which contain the following information:
\begin{align*}
\begin{array}{c|c|c|c}
\Phi(G) &\quad\Delta_{ext}(\Phi(G),id)\quad& \theta & \quad\Delta_{ext}(\Phi(G),\theta) \\
\hline \Phi(H)&&\theta=id &\quad \Delta_{ext}(\Phi(H),id)
\end{array}
\end{align*}
Here $\Delta_{ext}(\Phi(G),id)$ is arranged such that the $\theta$-orbits of roots are
in vertical order. We may think of $\Delta(\Phi(G),\theta)$ as a quotient diagram
of $\Delta(\Phi(G),id)$, but their seems to be no rule which describes the additional
root $c(\tilde\alpha)P_\Theta(\tilde\alpha)$. To obtain $\Delta_{ext}(\Phi(H),id)$ one observes at first that
 $\Phi(\hat G)\simeq\Phi(G)$, since $G$ is of type $ADE$,
 then remarks $\Delta(\Phi(\hat H),id)=\Delta(\Phi(\hat G),\theta)$, reversing the arrows in this
 diagram one gets the
 diagram of $\Delta(\Phi(H),id)$, which finally has to be extended to $\Delta_{ext}(\Phi(H),id)$.
 \Absatz
\Seitenumbruch
\begin{Tabelle}\em \qquad $\Delta_{ext}(\Phi(G),\theta)$ versus $\Delta_{ext}(\Phi(H),id)$
\begin{align*}
\begin{array}{c|c|c|c}
\Phi(G) &\Delta_{ext}(\Phi,id)& \theta & \Delta_{ext}(\Phi,\theta) \\
\hline\hline
\underset{n\ge 3}{A_{2n-1}} &
\begin{xy} <8mm,0mm>:
    (0,1)*{ }, (0,-1.5)*{ }, 
    (0,0)*{\bullet}\PATH
            ~={**\dir{-}}
    '(1, .5)*{\circ}
    '(2, .5)*{\circ} ~={**\dir{.}}
    '(3, .5)*{\circ} ~={**\dir{-}}
    '(4, .5)*{\circ}
    '(5,  0)*{\circ}
    '(4,-.5)*{\circ}
    '(3,-.5)*{\circ} ~={**\dir{.}}
    '(2,-.5)*{\circ} ~={**\dir{-}}
    '(1,-.5)*{\circ}
    '(0,  0)*{\circ}
\end{xy} &
\text{Ord}(\theta)=2
&\begin{xy} <8mm,0mm>:
    (0,1)*{ }, (0,-1)*{ }, 
    (4.5,0)*{<}, 
    (0,.5)*{\bullet}\PATH
            ~={**\dir{-}}
    '(1,0)*{\circ}
    '(2,0)*{\circ} ~={**\dir{.}}
    '(3,0)*{\circ} ~={**\dir{-}}
    '(4,0)*{\circ} ~={**\dir{=}}
    '(5,0)*{\circ} ~={**\dir{}}
    '(0,-.5)*{\circ} ~={**\dir{-}}
    '(1,0)*{\circ} ~={**\dir{}}
\end{xy}%
\\
\hline
\underset{n\ge3}{B_n} && \theta=id &
\begin{xy} <8mm,0mm>:
    (0,1)*{ }, (0,-1)*{ }, 
    (4.5,0)*{>},
    (0,.5)*{\bullet}\PATH
            ~={**\dir{-}}
    '(1,0)*{\circ}
    '(2,0)*{\circ} ~={**\dir{.}}
    '(3,0)*{\circ} ~={**\dir{-}}
    '(4,0)*{\circ} ~={**\dir{=}}
    '(5,0)*{\circ} ~={**\dir{}}
    '(0,-.5)*{\circ} ~={**\dir{-}}
    '(1,0)*{\circ}
\end{xy}%
\\
\hline\hline
A_3 & \begin{xy} <8mm,0mm>:
    (0,1)*{ }, (0,-1)*{ }, 
    (0,0)*{\bullet}\PATH
            ~={**\dir{-}}
    '(1, .5)*{\circ}
    '(2,  0)*{\circ}
    '(1,-.5)*{\circ}
    '(0,  0)*{\circ}
\end{xy} &
\text{Ord}(\theta)=2 &
\begin{xy} <8mm,0mm>:
    (.5,0)*{<}, (1.5,0)*{>},
    (0,1)*{ }, (0,-1)*{ }, 
    (0,0)*{\circ}\PATH ~={**\dir{=}}
    '(1,0)*{\circ}
    '(2,0)*{\bullet}
\end{xy}%
\\
\hline
B_2 && \theta=id &
\begin{xy} <8mm,0mm>:
    (0,1)*{ }, (0,-1)*{ }, 
    (0.5,0)*{>}, (1.5,0)*{<},
    (0,0)*{\circ}\PATH
            ~={**\dir{=}}
    '(1,0)*{\circ}
    '(2,0)*{\bullet}
\end{xy}%
\\
\hline\hline
\underset{n\ge 1}{A_{2n}} &
\begin{xy} <8mm,0mm>:
    (5.5,0)*{ }, (-.5,0)*{ }, (0,1)*{ }, (0,-1.5)*{ },
    (0,0)*{\bullet}\PATH
            ~={**\dir{-}}
    '(1, .5)*{\circ}
    '(2, .5)*{\circ} ~={**\dir{.}}
    '(3, .5)*{\circ} ~={**\dir{-}}
    '(4, .5)*{\circ}
    '(5, .5)*{\circ}
    '(5,-.5)*{\circ}
    '(4,-.5)*{\circ}
    '(3,-.5)*{\circ} ~={**\dir{.}}
    '(2,-.5)*{\circ} ~={**\dir{-}}
    '(1,-.5)*{\circ}
    '(0,  0)*{\circ}
\end{xy} &
\text{Ord}(\theta)=2 &
\begin{xy} <8mm,0mm>:
    (.5,0)*{>}, (4.5,0)*{>},    
    (4.5,-1)*{>},   
    (1.5,-1)*{\text{resp. if $n=1$}},
    (0,.5)*{},(0,-1.5)*{},
    (0,0)*{\bullet}\PATH
                                 ~={**\dir{=}}
    '(1,0)*{\circ}               ~={**\dir{-}}
    '(2,0)*{\circ}               ~={**\dir{.}}
    '(3,0)*{\circ}               ~={**\dir{-}}
    '(4,0)*{\circ}               ~={**\dir{=}}
    '(5,0)*{\circ}               ~={**\dir{ }}
    '(4,-1)*{\bullet}            ~={**\dir{=}}
    '(5,-1)*{\circ}              ~={**\dir{-}}
    '(4,-1)*{\circ}              ~={**\dir{ }}
    '(4,-1.1)*{\phantom{\circ}}  ~={**\dir{-}}
    '(5,-1.1)*{\phantom{\circ}}  ~={**\dir{ }}
    '(4,-.9)*{\phantom{\circ}}   ~={**\dir{-}}
    '(5,-.9)*{\phantom{\circ}}
\end{xy}
\\
\hline
\underset{n\ge 1}{C_n} && \theta=id  &
\begin{xy} <8mm,0mm>:
    (0,1)*{ }, (0,-1)*{ }, 
    (.5,0)*{>}, (4.5,0)*{<}, 
    (0,0)*{\bullet}\PATH ~={**\dir{=}}
    '(1,0)*{\circ} ~={**\dir{-}}
    '(2,0)*{\circ} ~={**\dir{.}}
    '(3,0)*{\circ} ~={**\dir{-}}
    '(4,0)*{\circ} ~={**\dir{=}}
    '(5,0)*{\circ}
\end{xy}%
\\
\hline\hline
\underset{n\ge 3}{D_{n+1}} &
\begin{xy} <8mm,0mm>:
    (0,1)*{ }, (0,-1)*{ }, 
    (0,.5)*{\bullet}\PATH
            ~={**\dir{-}}
    '(1,  0)*{\circ}
    '(2,  0)*{\circ} ~={**\dir{.}}
    '(3,  0)*{\circ} ~={**\dir{-}}
    '(4,  0)*{\circ}
    '(5, .5)*{\circ} ~={**\dir{}}
    '(5,-.5)*{\circ} ~={**\dir{-}}
    '(4,  0)*{\circ} ~={**\dir{}}
    '(0,-.5)*{\circ} ~={**\dir{-}}
    '(1,  0)*{\circ}
\end{xy} &
\text{Ord}(\theta)=2 &
\begin{xy} <8mm,0mm>:
    (0,1)*{ }, (0,-1)*{ }, 
    (.5,0)*{<}, (4.5,0)*{>}, 
    (0,0)*{\bullet}\PATH ~={**\dir{=}}
    '(1,0)*{\circ}       ~={**\dir{-}}
    '(2,0)*{\circ}       ~={**\dir{.}}
    '(3,0)*{\circ}       ~={**\dir{-}}
    '(4,0)*{\circ}       ~={**\dir{=}}
    '(5,0)*{\circ}
\end{xy}%
\\
\hline
\underset{n\ge 3}{C_n} & & \theta=id &
\begin{xy} <8mm,0mm>:
    (0,1)*{ }, (0,-1)*{ }, 
    (.5,0)*{>}, (4.5,0)*{<}, 
    (0,0)*{\bullet}\PATH
            ~={**\dir{=}}
    '(1,0)*{\circ} ~={**\dir{-}}
    '(2,0)*{\circ} ~={**\dir{.}}
    '(3,0)*{\circ} ~={**\dir{-}}
    '(4,0)*{\circ} ~={**\dir{=}}
    '(5,0)*{\circ}
\end{xy}%
 \\
\hline\hline
E_6 &
\begin{xy} <8mm,0mm>:
    (0,1)*{ }, (0,-1)*{ }, 
     (0,0)*{\bullet}\PATH
            ~={**\dir{-}}
    '(1, 0)*{\circ}
    '(2, 0)*{\circ}
    '(3, .5)*{\circ}
    '(4, .5)*{\circ} ~={**\dir{}}
    '(4,-.5)*{\circ} ~={**\dir{-}}
    '(3,-.5)*{\circ}
    '(2,0)*{\circ}
\end{xy} &
\text{Ord}(\theta)=2 &
\begin{xy} <8mm,0mm>:
    (0,1)*{ }, (0,-1)*{ }, 
    (1.5,0)*{>},
    (-1,0)*{}, (5,0)*{ }, 
    (0,0)*{\circ}\PATH
            ~={**\dir{-}}
    '(1,0)*{\circ} ~={**\dir{=}}
    '(2,0)*{\circ} ~={**\dir{-}}
    '(3,0)*{\circ}
    '(4,0)*{\bullet}
\end{xy}%
\\
\hline
F_4 && \theta=id &
\begin{xy} <8mm,0mm>:
    (0,1)*{ }, (0,-1)*{ }, 
    (1.5,0)*{<},
    (0,0)*{\circ}\PATH
            ~={**\dir{-}}
    '(1,0)*{\circ} ~={**\dir{=}}
    '(2,0)*{\circ} ~={**\dir{-}}
    '(3,0)*{\circ} ~={**\dir{-}}
    '(4,0)*{\bullet}
\end{xy}%
\\
\hline\hline
D_4 &
\begin{xy} <8mm,0mm>:
    (0,1)*{ }, (0,-1)*{ }, 
    (0,0)*{\bullet}\PATH
            ~={**\dir{-}}
    '(1,0)*{\circ}
    '(2,.5)*{\circ} ~={**\dir{}}
    '(2,-.5)*{\circ} ~={**\dir{-}}
    '(1,0)*{\circ}
    '(2,0)*{\circ}
\end{xy} &
   \text{Ord}(\theta)=3 &
\begin{xy} <8mm,0mm>:
    (0,1)*{ }, (0,-1)*{ }, 
    (.5,0)*{>},
    (0,0)*{\circ}\PATH
            ~={**\dir3{-}}
    '(1,0)*{\circ} ~={**\dir{-}}
    '(2,0)*{\bullet}
\end{xy}%
\\
\hline
G_2 && \theta=id &
\begin{xy} <8mm,0mm>:
    (0,1)*{ }, (0,-1)*{ }, 
    (0.5,0)*{<},
    (0,0)*{\circ}\PATH
            ~={**\dir3{-}}
    '(1,0)*{\circ} ~={**\dir{-}}
    '(2,0)*{\bullet}
\end{xy}%
\end{array}
\end{align*}
\end{Tabelle}

\Absatz\Numerierung\label{Motivation}{\bf Comparison of diagrams:} By construction the
 diagrams $\Delta_{ext}(\Phi(G),\theta)$ and $\Delta_{ext}(\Phi(H),id)$ are arranged in vertical order
 on the right hand side of each of the six blocks, and the corresponding spherical diagrams are obtained
 from each other by reversing the arrows. By inspection we see that the same statement holds for
 the extended diagrams with the exception that in case $A_{2n}\leftrightarrow C_n$ there is no
  reversal of the arrow joining the additional root with the standard diagram.
\absatz
 A similar phenomenon appears when we consider centralizers:
 If $s\theta\in G(\bar F)\theta$ and $\sigma\in H(\bar F)$ match they can be assumed to lie in
 the diagonal tori such that:
 $T(\bar F)\ni s\mapsto \sigma\in T_H(\bar F)=T_\Theta(\bar F)$. If we compute
 $\Phi(G^{s\theta})$ and $\Phi(H^\sigma)$ using \ref{Steinberg: G_hoch_tTheta ist reduktiv} and
 \ref{Dynkins Satz}(b) we see by inspection that they can be arranged in vertical order as subdiagrams
 (in the sense of of \ref{Dynkins Satz}(b)) of the diagrams
 $\Delta_{ext}(\Phi(G),\theta)$ and $\Delta_{ext}(\Phi(H),id)$ , but the
 oriented arrows get reversed with the exception mentioned above.
  If we write $(G^{s\Theta})_{ad}$ resp. $(H^\sigma)_{ad}$
 as product of simple groups, we therefore get common factors of type $A$ and $D$ but factors
 of type $B$ in general correspond to factors of type $C$ and vice versa and factors of type
 $G_2$ and $F_4$ come in with reversed arrows.
\Absatz\Numerierung\label{Outline}
  Our strategy to prove the fundamental lemma in the case of classical split groups
  now goes as follows: By the Kazhdan lemma
  \ref{Kazhdan-Lemma} the (twisted) stable orbital integral of $g\theta\in G(F)\theta$ in the group
  $\tilde G$ can be replaced by the ordinary stable orbital integral in the group $G^{s\theta}$
  of the topologically unipotent part $u$. Now $G^{s\theta}$ resp. $H^\sigma$ is isogenous to a product
  $G^{s\theta}_+\times G_*^{s\theta}$ resp. $H^{\sigma}_+\times H_*^{\sigma}$ ,
  such that $G^{s\theta}_+$ is of type $B$ or $C$, $H^\sigma_+$ is the other of these two types and
  $G_*^{s\theta}$ is isogenous to $H_*^\sigma$. Similarly the stable orbital integral of some
  $\gamma\in H(F)$, which matches with $g$, can be computed as the stable orbital integral of
  the topologically unipotent part $v$ in the group $H^\sigma$, where the residually semisimple
  part $\sigma$ of $\gamma$ matches with $s\theta$. Decomposing $u=(u_+,u_*)$ and $v=(v_+,v_*)$
  we get that $u_*$ and $v_*$ coincide up to stable conjugation und up to powering, so they have
  matching stable orbital integrals. The fundamental lemma for $g\theta$ and $\gamma$ will
  now follow if we can assume that the stable orbital integrals of $u_+$ and $v_+$ match,
  which is essentially the $BC$ conjecture \ref{BCVermutung}, stated already in the introduction.

\Absatz
\Numerierung\label{sc oder ad genuegt fuer Zentr}
In case $\Phi=A_{2n}$ consider the exact sequence:
\[ 1\longrightarrow \mu_{2n+1}\longrightarrow \SL_{2n+1}\longrightarrow \PGL_{2n+1}\longrightarrow 1. \]
Since the involution $\Theta$ acts as $-1$ on $\mu_{2n+1}\simeq\ZZ/(2n+1)\ZZ$
one gets the long exact sequence
\[1\longmapsto\underbrace{\kappa(G)^{\Theta}}_{=1}\longrightarrow G^{t\Theta}
    \longrightarrow G_{ad}^{t\Theta}\longrightarrow
    \underbrace{H^1(\langle\Theta\rangle,\kappa(G))}_{=1}\longrightarrow...
\]
for every $t\in\text{Inn}(G)$.
\begin{expl} \em {\bf The case $\Phi=A_{4}$}\\
Consider $s\in T(\bar F)$, where $T$ is the diagonal torus in $G=\PGL_5$.
There are seven possible types of groups $G^{s\Theta}$ nonisomorphic over $\bar F$.
($\zeta_n$ denotes a (fixed) primitive $n$--th root of unity in $\bar F$.)
\end{expl}
\begin{align}
\tag*{\bf 1. Case:}
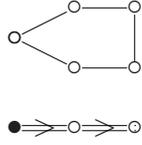
\begin{xy} <8mm,0mm>:
    (.5,-1.5)*{>}, (1.5,-1.5)*{>},
    (0,0)*{\circ}\PATH
                        ~={**\dir{-}}
    '(1, .5)*{\circ}
    '(2, .5)*{\circ}
    '(2,-.5)*{\circ}
    '(1,-.5)*{\circ}
    '(0,  0)*{\circ}    ~={**\dir{}}
    '(0,-1.5)*{\bullet} ~={**\dir{=}}
    '(1,-1.5)*{\circ}
    '(2,-1.5)*{\circ}
\end{xy}
\qquad\qquad
s=1\phantom{mmmmm}
\end{align}
$G^{s\Theta}=\text{SO}(5)=\{g\in\text{Sl}(5)\mid{}^tg\cdot J\cdot g =J\}$
is of type $\Phi(G^{s\Theta})=B_2$ and $\pi_1(\Phi(G^{s\Theta}))\simeq\ZZ/2$.
\begin{align}
\tag*{\bf 2. Case:}
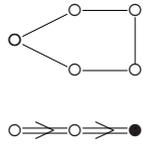
\begin{xy} <8mm,0mm>:
    (.5,-1.5)*{>}, (1.5,-1.5)*{>},
    (0,0)*{\circ}\PATH
                        ~={**\dir{-}}
    '(1, .5)*{\circ}
    '(2, .5)*{\circ}
    '(2,-.5)*{\circ}
    '(1,-.5)*{\circ}
    '(0,  0)*{\circ}    ~={**\dir{}}
    '(0,-1.5)*{\circ} ~={**\dir{=}}
    '(1,-1.5)*{\circ}
    '(2,-1.5)*{\bullet}
\end{xy}
\qquad\qquad
s=\left(\begin{smallmatrix}
        \zeta_4 &\cdot  &\cdot  &\cdot      &\cdot \\
        \cdot   &\zeta_4 &\cdot &\cdot      &\cdot \\
        \cdot   &\cdot  &1      &\cdot      &\cdot \\
        \cdot   &\cdot  &\cdot  &\zeta_4^{-1} &\cdot \\
        \cdot   &\cdot  &\cdot  &\cdot      &\zeta_4^{-1} \\
  \end{smallmatrix}\right)
\end{align}
$G^{s\Theta}=\text{Sp}(2\cdot2)=\{g\in\text{Sl}(5)\mid{}^tg\cdot \tilde J\cdot g =\tilde J\}$,
where $\tilde J=J\cdot Diag(-1,-1,1,1,1)$,
is of type $\Phi(G^{s\Theta})=C_2$ and $\pi_1(\Phi(G^{s\Theta}))=1$.
\begin{align}
\tag*{\bf 3. Case:}
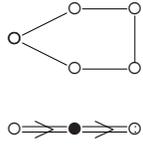
\begin{xy} <8mm,0mm>:
    (.5,-1.5)*{>}, (1.5,-1.5)*{>},
    (0,0)*{\circ}\PATH
                        ~={**\dir{-}}
    '(1, .5)*{\circ}
    '(2, .5)*{\circ}
    '(2,-.5)*{\circ}
    '(1,-.5)*{\circ}
    '(0,  0)*{\circ}    ~={**\dir{}}
    '(0,-1.5)*{\circ} ~={**\dir{=}}
    '(1,-1.5)*{\bullet}
    '(2,-1.5)*{\circ}
\end{xy}
\qquad\qquad
s=\left(\begin{smallmatrix}
        \zeta_4 &\cdot &\cdot &\cdot &\cdot \\
        \cdot &1     &\cdot &\cdot &\cdot \\
        \cdot &\cdot &1     &\cdot &\cdot \\
        \cdot &\cdot &\cdot &1    &\cdot \\
        \cdot &\cdot &\cdot &\cdot &\zeta_4^{-1} \\
  \end{smallmatrix}\right)
\end{align}
$G^{s\Theta}\simeq \text{SO}(3)\times \text{Sp}(2)$,
is of type $\Phi(G^{s\Theta})=B_1\times C_1=A_1^2$ and $\pi_1(\Phi(G^{s\Theta}))\simeq\ZZ/2$.
\begin{align}
\tag*{\bf 4. Case:}
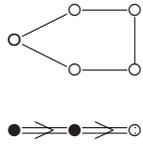
\begin{xy} <8mm,0mm>:
    (.5,-1.5)*{>}, (1.5,-1.5)*{>},
    (0,0)*{\circ}\PATH
                        ~={**\dir{-}}
    '(1, .5)*{\circ}
    '(2, .5)*{\circ}
    '(2,-.5)*{\circ}
    '(1,-.5)*{\circ}
    '(0,  0)*{\circ}    ~={**\dir{}}
    '(0,-1.5)*{\bullet} ~={**\dir{=}}
    '(1,-1.5)*{\bullet}
    '(2,-1.5)*{\circ}
\end{xy}
\qquad\qquad
s=\left(\begin{smallmatrix}
        \zeta_8 &\cdot &\cdot &\cdot &\cdot \\
        \cdot   &1     &\cdot &\cdot &\cdot \\
        \cdot   &\cdot &1     &\cdot &\cdot \\
        \cdot   &\cdot &\cdot &1     &\cdot \\
        \cdot   &\cdot &\cdot &\cdot &\zeta_8^{-1} \\
  \end{smallmatrix}\right)
\end{align}
$G^{s\Theta}\simeq\text{SO}(3)\times\text{Gl}(1)$
\begin{align}
\tag*{\bf 5. Case:}
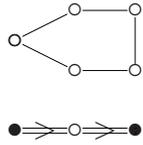
\begin{xy} <8mm,0mm>:
    (.5,-1.5)*{>}, (1.5,-1.5)*{>},
    (0,0)*{\circ}\PATH
                        ~={**\dir{-}}
    '(1, .5)*{\circ}
    '(2, .5)*{\circ}
    '(2,-.5)*{\circ}
    '(1,-.5)*{\circ}
    '(0,  0)*{\circ}    ~={**\dir{}}
    '(0,-1.5)*{\bullet} ~={**\dir{=}}
    '(1,-1.5)*{\circ}
    '(2,-1.5)*{\bullet}
\end{xy}
\qquad\qquad
s=\left(\begin{smallmatrix}
        \zeta_8 &\cdot  &\cdot  &\cdot      &\cdot \\
        \cdot   &\zeta_8 &\cdot &\cdot      &\cdot \\
        \cdot   &\cdot  &1      &\cdot      &\cdot \\
        \cdot   &\cdot  &\cdot  &\zeta_8^{-1} &\cdot \\
        \cdot   &\cdot  &\cdot  &\cdot      &\zeta_8^{-1} \\
  \end{smallmatrix}\right)
\end{align}
$G^{s\Theta}\simeq\text{Gl}(2)$
\begin{align}
\tag*{\bf 6. Case:}
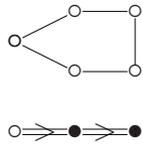
\begin{xy} <8mm,0mm>:
    (.5,-1.5)*{>}, (1.5,-1.5)*{>},
    (0,0)*{\circ}\PATH
                        ~={**\dir{-}}
    '(1, .5)*{\circ}
    '(2, .5)*{\circ}
    '(2,-.5)*{\circ}
    '(1,-.5)*{\circ}
    '(0,  0)*{\circ}    ~={**\dir{}}
    '(0,-1.5)*{\circ} ~={**\dir{=}}
    '(1,-1.5)*{\bullet}
    '(2,-1.5)*{\bullet}
\end{xy}
\qquad\qquad
s=\left(\begin{smallmatrix}
        \zeta_4 &\cdot  &\cdot  &\cdot      &\cdot \\
        \cdot   &\zeta_8 &\cdot &\cdot      &\cdot \\
        \cdot   &\cdot  &1      &\cdot      &\cdot \\
        \cdot   &\cdot  &\cdot  &\zeta_8^{-1} &\cdot \\
        \cdot   &\cdot  &\cdot  &\cdot      &\zeta_4^{-1} \\
  \end{smallmatrix}\right)
\end{align}
$G^{s\Theta}\simeq\text{Sp}(2)\times\text{Gl}(1)$
\begin{align}
\tag*{\bf 7. Case:}
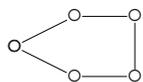
\begin{xy} <8mm,0mm>:
    (.5,-1.5)*{>}, (1.5,-1.5)*{>},
    (0,0)*{\circ}\PATH
                        ~={**\dir{-}}
    '(1, .5)*{\circ}
    '(2, .5)*{\circ}
    '(2,-.5)*{\circ}
    '(1,-.5)*{\circ}
    '(0,  0)*{\circ}    ~={**\dir{}}
    '(0,-1.5)*{\bullet} ~={**\dir{=}}
    '(1,-1.5)*{\bullet}
    '(2,-1.5)*{\bullet}
\end{xy}
\qquad\qquad
s=\left(\begin{smallmatrix}
        \zeta_6 &\cdot  &\cdot  &\cdot      &\cdot \\
        \cdot   &\zeta_{12} &\cdot &\cdot      &\cdot \\
        \cdot   &\cdot  &1      &\cdot      &\cdot \\
        \cdot   &\cdot  &\cdot  &\zeta_{12}^{-1} &\cdot \\
        \cdot   &\cdot  &\cdot  &\cdot      &\zeta_6^{-1} \\
  \end{smallmatrix}\right)
\end{align}
$G^{s\Theta}\simeq\text{Gl}(1)\times\text{Gl}(1)$

\Seitenumbruch
\meinchapter{Topological Jordan decomposition}
\begin{Def}\label{Die Kategorie tJZ}
  Let $\mathbf{tJZ}_p$ be the category,
  \index{$\mathbf{tJZ}_p$}
  whose objects are topological groups, such that the neighborhood filter of
  $1$ has a basis consisting of pro--$p$--groups, and whose morphisms are
  the continuous homomorphisms.
\end{Def}
\begin{Def}
An element $g$ of $G\in Ob(\mathbf{tJZ}_p)$ is called
\begin{itemize}
\item strongly compact, if  $g$ lies in a compact subgroup of $G$.
\index{stark kompakt}
\item topologically unipotent, if $\lim_{n\to\infty}g^{p^{n}}=1$.
\index{topologisch unipotent}
\index{residuell halbeinfach}
\item residually semisimple, if $g$ is of a finite order, which is prime to $p$.
\end{itemize}
\end{Def}
\Numerierung\label{Hales' abs. sesi. Beschreibung}
For an element $g$ it is equivalent to be topologically unipotent and to lie in a pro--$p$--subgroup
of $G$. In the definition above one can replace the sequence  $(p^n)$ by an arbitrary
sequence $(p^{n_k})_k$
satisfying $\lim_{k\to\infty}n_k=\infty$.
\begin{expl}
\label{G(F) in tJZ-Kategorie}
\em
Let $F$ be a $p$--adic field and $\tilde G$ an affine linear algebraic group.
Then $\tilde G(F)$ in an Object of $\mathbf{tJZ}_p$.
\end{expl}

An element  $g\in\tilde G(F)$ in the example \ref{G(F) in tJZ-Kategorie} is strongly compact,
iff the set $g^\Z$ is bounded inside $\tilde G(F)$.
A third equivalent formulation is,
that the eigenvalues of $\rho(g)$ are units in $\bar F^\times$ for one/for all faithful
representation(s) $\rho: G\rightarrow GL(V)$.

\absatz
Each pro--$p$--group $U$ has a unique structure as a topological $\Z_p$--module, which extends the
canonical structure as  $\Z$--module (comp. \cite[\S15.2]{Hasse}).
\begin{lemma}[topological Jordan decomposition]\label{allgemeine topologische Jordanzerlegung}
\index{topologische Jordanzerlegung}
Let $G$ be an object of $\mathbf{tJZ}_p$.
Every strongly compact $g\in G$ has a unique decomposition
\[g = g_{u}\cdot g_{s} = g_{s}\cdot g_{u}\;,\]
where $g_{u}\in G$ is topologically unipotent and $g_{s}\in G$ is residually semisimple.

\end{lemma}
Proof:
We consider the (abelian!) closure $\ov{\langle g^\Z\rangle}$ of the abelian group $\langle g^\Z\rangle$,
which is contained in a compact subgroup of $G$ and is therefore itself compact.
Since $G\in Ob(\mathbf{tJZ}_p)$ there exists an open pro--$p$--subgroup of $\ov{\langle g^\Z\rangle}$.
The set $U$ of all topologically unipotent
elements in $\ov{\langle g^\Z\rangle}$ contains this open subgroup, is a group since
$\ov{\langle g^\Z\rangle}$ is abelian, and is therefore an open pro--$p$--subgroup $U$.
The compactness of $\ov{\langle g^\Z\rangle}$ implies that $U$ has a finite index $N$ in it, which
has to be prime to $p$. Since $U$ is a $\Z_p$--module and $N\in \Z_p^\times$ there exists a
(topologically unipotent!) element
$g_u\in U$ such that $g^N=g_u^N$. Since $g_ug=gg_u$ the element $g_s=g\cdot g_u^{-1}$ satisfies
 $g_s^N=1$, i.e.
is residually semisimple, and we get $g=g_sg_u$ and
$g=g\cdot g_u\cdot g_u^{-1}=g_u\cdot g\cdot g_u^{-1}=g_ug_s$,
i.e the existence of the decomposition is proved.
\absatz
If $g=g'_ug'_s=g'_sg'_u$ is a second topological Jordan decomposition with $(g'_s)^{N'}=1$ we choose a
$p$-power $Q=p^m$ such that $Q\equiv 1\mod NN'$ and get $\lim_{\alpha\to\infty} g^{Q^\alpha}
=\lim_{\alpha\to\infty} (g'_u)^{Q^\alpha}\cdot (g'_s)^{Q^\alpha}
=\lim_{\alpha\to\infty} (g'_u)^{Q^\alpha}\cdot g'_s
=g'_s$, and by the same argument: $\lim_{\alpha\to\infty} g^{Q^\alpha}=g_s$, i.e.
$g'_s=g_s$ and therefore also $g_u=g'_u$ i.e. the uniqueness assertion.
\qed

\begin{cor}[Properties of the topological Jordan decomposition]\label{top Jordan Eigenschaften}
$\ $
\begin{engeListe}
\item[{\bf (1)}] Let $g\in G$ be strongly compact, $N\in\N$
    be prime to $p$, such that $g^N$ lies in some pro--$p$--group,
    and let $Q$ be a $p$--power with $Q\equiv 1$ (mod $N$).
    Then we have
\[  \lim_{m\to\infty}g^{Q^m}=g_s.
\]
\item[{\bf (2)}] We have $g_u\in G^{g_s}$ and $G^g=\Cent(g_u, G^{g_s})$.
\item[{\bf (3)}] Residually semisimple elements are semisimple.
\item[{\bf (4)}] Let $m$ be prime to $p$ and $u$ be topologically unipotent. Then there exists a unique
topologically unipotent $u_1$ such that \[  u_1^m=u.\]
\item[{\bf (5)}]  The topological Jordan decomposition is functorial in the following sense:
    The strongly compact  (resp. the topologically unipotent, resp. the residually semi\-simple) elements
    define functors from $\mathbf{tJZ}_p$ to $\mathbf{Set}$.
    For each morphism $\pHi$ in $\mathbf{tJZ}_p$ we have
    $\pHi((\cdot)_s)=(\pHi(\cdot))_s$ und $\pHi((\cdot)_u)=(\pHi(\cdot))_u$.
\neueZeile
    Especially:
\item[{\bf (5a)}]
    If $H$ is a closed subgroup of $G$ and $g\in H$,
    then also  $g_s$ and $g_u$ are in $H$.
\end{engeListe}
\end{cor}
\begin{cor}\label{topologische Jordanzerlegung}
Let $F$ and $\tilde G$ be as in example \ref{G(F) in tJZ-Kategorie} and assume furthermore that
 $|\pi_0(\tilde G)|$ is prime to $p$.
Each strongly compact element $g\in\tilde G(F)$ has a unique
topological  Jordan decomposition:
\[g = g_{u}\cdot g_{s} = g_{s}\cdot g_{u}\;,\]
where $g_{s}\in\tilde G(F)$ is residually semisimple and $g_{u}\in (\tilde G)^\circ(F)$
topologically unipotent.

\Absatz
The functoriality implies the following statements:
\begin{engeListe}
\item[{\em (1)}] Let $\rho: \tilde G\rightarrow\tilde G'$ be a morphism of
    (not necessarily connected)
    reductive groups, defined over a finite extension of $F$.
    Then we have $\rho(g)_s=\rho(g_s)$ and $\rho(g)_u=\rho(g_u)$.
\item[{\em (2)}] If $g\in\tilde G(\mathcal{O}_F)$, then the image of the topological Jordan decomposition
    under the  reduction map  is the Jordan decomposition in $\tilde G(\F_q)$.
\end{engeListe}
\end{cor}
(Topologically unipotent elements must lie in $(\tilde G)^\circ$ since by assumption
$p$ does not divide $|\pi_0(\tilde G)|$.)
        \Seitenumbruch

\meinchapter{Classification of $\Theta$-conjugacy classes}

\Numerierung If $(G,\Theta)$ is as in the examples  \ref{BspPGL2n+1Sp2n} or \ref{BspGL2nGl1GSpin2n+1},
  the problem of determining the $\Theta$-conjugacy classes of elements $s\Theta\in \tilde G(F)$
  is equivalent to determine the classes of $h=sJ$ under the transformations
  $h\mapsto g\cdot h\cdot{}^t g$. Namely we have the following commutative diagram:
  \begin{align*}\begin{CD}
    \GL_{2n} @>{h\mapsto hJ^{-1}\Theta}>> \widetilde GL_{2n}\\
    @V{h\mapsto g\cdot h\cdot \tran g}VV  @VV{x\mapsto g\cdot x\cdot g^{-1}}V \\
    \GL_{2n} @>>{h\mapsto hJ^{-1}\Theta}> \widetilde GL_{2n}
  \end{CD}\end{align*}
  If we decompose $h=q+p$ in its symplectic part $p$
  and its symmetric part $q$ we thus have to consider the problem of simultaneous normal forms for a
  symplectic and a symmetric bilinear form. To obtain results for orbital integrals we
  have to deal with this problem also over the ring of integers $\OOO_F$. The problem
  can be attacked if we assume $s\Theta$ to be semisimple (resp. residually semisimple
  if we work over $\OOO_F$).

\Absatz
\Numerierung \label{R theta semisimple}
{\bf Notations:} In the following $R$ denotes either a field of characteristic $0$
or the ring of integers $\OOO_F$
of a local $p$-adic field $F$, where $p\ne 2$. We denote by $\mMm$ the maximal ideal
of $R$ (i.e. $\mMm=(0)$ if $R$ is a field) and by $\kappa=R/\mMm$ the residue field in
the case $R=\OOO_F$.
\absatz
Let $M$ denote a free $R$-module of finite rank $r$ with basis $(b_i)_{1\le i\le r}$.
 A bilinear form $q:M\times M\to R$ is called
{\it unimodular} if $\Delta(q):= \det(q(b_i,b_j))\in R^\times$. This definition is obviously
independent of the chosen basis $(b_i)$ since $\Delta(q)$ is an invariant in $R/(R^\times)^2$.
 For $h\in \GL_n(R)$ we have the following
bilinear forms $b_h$ and $b'_h$ on the module $M=R^n$ of column vectors:
 $b_h(m_1,m_2)={}^tm_1\cdot h\cdot m_2$ and $b'_h(m_1,m_2)={}^tm_1\cdot{}^th\cdot m_2$.
\absatz
An element $g\in \GL (M)$ is called $R$-{\it semisimple} iff
\begin{engeListe}
\item[$\bullet$] $g$ is semisimple in the case $R$ is a field,
\item[$\bullet$] $g$ is residually semisimple
(i.e. has finite order prime to $char(\kappa)$) in the case $R=\OOO_F$.
\end{engeListe}
\medskip
For $h\in \GL_n(R)$ we call $N(h)=h\cdot {}^t h^{-1}$ the {\it(right) norm} of $h$ and $N_l(h)=
{}^t h^{-1}\cdot h$ the {\it left norm} of $h$. $N(h)$ and $N_l(h)$ are conjugate by $h$ in
$\GL_n(R)$. Then $h$ is called {\it $R$-$\Theta$-semisimple} if
$N(h)$ (or equivalently $N_l(h)$) is $R$-semisimple.
\absatz
We remark that $h$ is $R$-$\Theta$-semisimple if and only if $h\cdot J^{-1}\cdot \Theta$ is semisimple respectively
 residually semisimple as an element of $\GL_n(R)\rtimes\langle\Theta\rangle$.
\begin{lemma} \label{symplektischeNormalform}
If $p:M\times M\to R$ is a unimodular symplectic form, then there exists a
basis $(e_1,\ldots,e_g,f_g,\ldots,f_1)$ of $M$, such that $p$ has standard form with
respect to this basis, i.e. $p(e_i,e_j)=p(f_i,f_j)=0$ and $p(e_i,f_j)=\delta_{ij}$.
\end{lemma}
Proof: The standard procedure to get a symplectic basis of $M$ applies for unimodular forms.
\qed

\begin{lemma} \label{orthogonaleNormalform}
If $q:M\times M\to R$ is a unimodular symmetric bilinear form and $R=\OOO_F$, then
there exists a basis $(e_i)_{1\le i\le r}$ of $M$ such that $q(e_i,e_j)=\delta_{ij}$ for $(i,j)\ne(r,r)$
and $q(e_r,e_r)$ is some given element in the class of $\Delta(q)$ in $R^\times/(R^\times)^2$.
\end{lemma}
Proof:
Consider the reductions $\kappa=R/\mMm,\quad \bar{M} =M/\mMm M$ and $\bar{q}:\bar{M}\times\bar{M}
\to \kappa.$ Since quadratic forms over finite fields are classified by their discriminants,
the analogous statement for $\bar q$ holds. By lifting a basis from $\bar{M}$ to $M$
 we can therefore assume that
$q(b_i,b_j)\cong \delta_{ij}\mod \mMm$ for $(i,j)\ne (r,r)$. But now we can apply
the Gram-Schmidt-Orthogonalization  procedure (observe that elements
congruent to $1$ modulo $\mMm$ are squares since $p\ne 2$) to obtain the claim.
\qed

\begin{lemma}
  $\ $ \label{Eigenwertzerlegung}
  \begin{itemize}
  \item[(a)]
    If $g\in \GL(M)$ is $R$-semisimple then there exists a finite \'etale galois extension $R'/R$ such
    that $M'=M\otimes_RR'$ decomposes into the direct sum of eigenspaces: $M'=\bigoplus_\lambda M'_\lambda$,
    where $g$ acts on $M'_\lambda$ as the scalar $\lambda$.
  \item[(b)]
    If $g=N_l(h)$ for an $R$-$\Theta$-semisimple $h\in \GL_n(R)$ (see \ref{R theta semisimple}) then
    $b_h(M'_\lambda,M'_\mu)=0=b'_h(M'_\lambda,M'_\mu)$ unless $\lambda\mu=1$.
  \item[(c)] The restrictions of the forms $b_h$ and $b'_h$ to $M'_1$ and $M'_{-1}$ are unimodular.
    For $\lambda\ne \pm 1$ also the restrictions of $b_h,b'_h, b_h+b'_h$ and $b_h-b'_h$ to the modules
    $N'_\lambda =M'_\lambda\oplus M'_{\lambda^{-1}}$ are unimodular.
  \end{itemize}
\end{lemma}
Proof: (a) The minimal polynomial $\chi(X)$ of $g$ decomposes in pairwise different linear factors
  $\chi(X)=\prod_{i=1}^r (X-\lambda_i)$ over some extension ring of $R$. The ring
  $R'=R[\lambda_i]_{1\le i\le r}$ is finite \'etale and galois over $R$, since the $\lambda_i$ are roots of
  unity of order prime to $char(\kappa)$ in the case $R=\OOO_F$. By the same reason we have
  \begin{gather}\label{differenzensindeinheiten}
    \lambda_i-\lambda_j\in (R')^\times\qquad\text{ for }i\ne j
  \end{gather}
  in both cases for $R$. We remark for later use that this statement remains correct if
  we add $\pm 1$ to the set of the $\lambda_i$ (if they are not already among them).  Therefore
  $\chi_i(X)=\prod_{j\ne i} ((X-\lambda_j)\cdot(\lambda_i-\lambda_j)^{-1})\in R'[X]$. We have
  $\sum_{i=1}^r\chi_i(X)=1$ since the left hand side is a polynomial of degree $r-1$ which has the
  value $1$ at $r$ different places. Therefore $M'$ is the sum of the subspaces
  $M'_{\lambda_i}=\chi_i(g)(M)$. Since   $(g-\lambda_i)\cdot\chi_i(g)$ equals
  $\chi(g)\cdot\prod_{j\ne i}(\lambda_i-\lambda_j)^{-1}=0$,
  the spaces $M'_{\lambda_i}$ are eigenspaces for $g$ and the sum
  $M'=\Sigma_{i=1}^r M'_{\lambda_i}$ is direct.
 \absatz
(b) For $m\in M'_\lambda$ and $n\in M'_\mu$ we have $m=\lambda^{-1}\cdot gm$ and $n=\mu\cdot g^{-1}n$.
  The claims follow immediately from the relations
  $ b_h(m,n)=
  \lambda^{-1}\cdot b'_h(m,n)$ and
  \begin{gather}\label{Gleichungmu}
     b_h(m,n)=
     \mu \cdot b'_h(m,n).
  \end{gather}
 \absatz
 (c) In view of the orthogonality relations (b) and the unimodularity of $h$ and $\tran h$ the claims
  for the restrictions of $b_h$ and $b'_h$ follow immediately. By the formula (\ref{Gleichungmu})
  above we have for
  $m\in M_\lambda, n\in M_{\lambda^{-1}}$:
  \begin{align*} (b_h\pm b'_h)(m,n)&=(1\pm\lambda)b_h(m,n) \\
    (b_h\pm b'_h)(n,m)&=(1\pm\lambda^{-1})b_h(n,m).
  \end{align*}
  Since $1\pm\lambda,1\pm\lambda^{-1}\in (R')^\times$ by the remark following (\ref{differenzensindeinheiten})
  above, the claim follows also for the restrictions of $b_h\pm b'_h$.
\qed

\begin{lemma}\label{+-0Zerlegung}
  For an $R$-$\Theta$-semisimple $h\in \GL_n(R)$ with decomposition $h=p+q$, where $p$ is skew-symmetric
  and $q$ symmetric, we have a direct sum decomposition for $M=R^n$
  $$ M = M_+ \oplus M_- \oplus M_0,$$
  where $M_+=\ker p, M_-=\ker q$ and $M_0=(M_+)^{\bot_q}\cap (M_-)^{\bot_p}$ is the intersection of
  the orthogonal complement of $M_+$ with the symplectic orthogonal complement of $M_-$.
  The restrictions
  \begin{align*}
    q_+=q\mid M_+\times &M_+, \qquad p_-=p\mid M_-\times M_-, \\
    q_0=q\mid M_0\times &M_0, \qquad p_0 = p\mid M_0\times M_0
  \end{align*}
  are unimodular.
\end{lemma}
Proof: We identify the matrices $p,q\in \GL_n(R)$ with the forms $b_p,b_q$.
  We take an extension $R'/R$ as in Lemma \ref{Eigenwertzerlegung}(a) and compute
  $$M'_{\pm 1}=\{m\in M'\mid \tran h^{-1}\cdot h\cdot m=\pm m\}=
   \{m\in M'\mid hm=\pm\tran hm\}=\ker(h\mp \tran h).$$
  This means $M'_1=\ker(p\mid M')$ and $M'_{-1}=\ker(q\mid M')$ and implies
  $M'_1=M_+\otimes_RR', M'_ {-1}=M_-\otimes_RR'$. Since unimodularity can be
  checked after the extension $R'/R$ and $b_h$ restricts to $q_+$ resp. $p_-$
  on $M_+$ resp. $M_-$, we conclude from Lemma \ref{Eigenwertzerlegung}(c) that $q_+$ and
  $p_-$ are unimodular. Then it is clear that we have the claimed decomposition in (orthogonal and
  symplectic orthogonal) direct summands. By Lemma \ref{Eigenwertzerlegung}(b) we get
  $M_0\otimes_RR'=\bigoplus_{\lambda\ne\pm 1}M'_\lambda$. By Lemma \ref{Eigenwertzerlegung}(c) again we
  conclude that the restrictions of $b_h+b'_h=2q$ and $b_h-b'_h=2p$ to this module are unimodular.
  So $p_0$ and $q_0$ are unimodular.
\qed

\begin{lemma}[Cayley transformation]\label{Cayleytransformation}
  Let $p\in\GL_n(R)$ be a skew-symmetric matrix. Let
  $Sym_n(R)_{p-ess}$ denote the set of symmetric matrices $q$ such that $q\pm p\in \GL_n(R)$
  and $\Sp(p,R)_{ess}$ the set of symplectic transformations $b$
  such that $b-1\in \GL_n(R)$. Then the following holds:
  \begin{itemize}
  \item[(a)] We have a bijection
    \begin{gather*}
     C: Sym_n(R)_{p-ess} \rightarrow \Sp(p,R)_{ess},\quad q\mapsto (q-p)^{-1}\cdot (q+p)=N_l(p+q).
    \end{gather*}
    The inverse map
    is $C^{-1}: b\mapsto p\cdot (b+1)\cdot(b-1)^{-1}$.
  \item[(b)] $C$ induces a bijection between those elements $q$ of $Sym_n(R)_{p-ess}$,
    for which $p+q$ is $R$-$\Theta$-semisimple, and the $R$-semisimple elements of $\Sp(p,R)_{ess}$.
  \item[(c)] The map $C$ satisfies $C(\tran g\cdot q \cdot g)=g^{-1}\cdot C(q)\cdot g$ for
    $g\in \Sp(p,R)$.
  \end{itemize}
\end{lemma}
Proof: (a) For $q\in Sym_n(R)_{p-ess}$ we put $h=p+q$ and $b=\tran h^{-1}\cdot h$. We have
  \begin{align}
    \notag \tran b\cdot h\cdot b\quad &= \quad\tran h\cdot h^{-1}\cdot h\cdot\tran h^{-1}\cdot h
                                 \quad = \quad\tran h\cdot \tran h^{-1}\cdot h \quad\text{ i.e.}\\
    \label{tranbhbisth} \tran b\cdot h\cdot b\quad &= \quad h\qquad \text{and by transposing}\\
    \label{tranbtranhbisttranh}\tran b\cdot \tran h\cdot b\quad &= \quad \tran h.
  \end{align}
  Subtracting the last two equations we get $\tran b\cdot p \cdot b=p$, i.e. $b\in \SP(p,R)$.
  Furthermore $b-1=(q-p)^{-1}\cdot((p+q)-(q-p))=(q-p)^{-1}\cdot 2p\in \GL_n(R)$
  by the assumptions. The map $C$ is therefore defined.
  \klabsatz
  Conversely we get for $b\in \Sp(p,R)_{ess}$ and
  $q=p\cdot(b+1)\cdot(b-1)^{-1}$ the equivalences:
  \begin{align*}q=\tran q&\Leftrightarrow
  p\cdot(b+1)\cdot(b-1)^{-1}=\tran( b-1)^{-1}\cdot\tran(b+1)\cdot(-p)\\
  &\Leftrightarrow (\tran b-1)p(b+1)=(\tran b+1)p(1-b)\\
  &\Leftrightarrow \tran bpb +\tran bp-pb-p=-\tran b pb+\tran bp-pb+p\\&\Leftrightarrow \tran bpb=p
  \Leftrightarrow b\in \Sp(p,R).
  \end{align*}
  Furthermore $q\pm p= p\cdot\left((b+1)\pm(b-1)\right)\cdot(b-1)^{-1}\in \GL_n(R)$
  since $(b-1)^{-1},2b,2,p\in \GL_n(R)$.
  Therefore the map $C^{-1}$ is also well defined. An easy calculation (as in the case of the usual Cayley transform)
  shows that the maps $C$ and $C^{-1}$ are inverse to another in their domain of definition.
  \absatz
(b) Since $C(q)=N_l(p+q)=(p+q)^{-1}\cdot N(p+q)\cdot (p+q)$ this follows from
  the definition of $R$-$\Theta$-semisimplicity.
  \absatz
(c) We have $C(\tran g\cdot q \cdot g )=(\tran gq g-p)^{-1}(\tran gqg+p)=g^{-1}(q-p)\tran g^{-1}\cdot
  \tran g(q+p)g=g^{-1}\cdot(q-p)^{-1}\cdot(q+p)\cdot g=g^{-1}\cdot C(q)\cdot g$ for $g\in Sp(p,R)$.
  \qed

\begin{lemma}\label{Abspaltung1Sp}
  If $p$ is a unimodular symplectic form on a free $R$-module $N$ and $b\in \Sp(p,R)$ is $R$-semisimple
  then there exists a $b$-invariant and with respect to $p$ orthogonal direct sum decomposition
  $N=N_1\oplus N_*$ such that $b$ acts as identity on $N_1$ and
  $b|N_* \in \Sp(p_*,R)_{ess}$, where $p_*$ is the restriction of $p$ to $N_*$.
\end{lemma}
Proof: We argue as before: By lemma \ref{Eigenwertzerlegung}(a) we have for some
  finite \'{e}tale ring extension $R'/R$  a decomposition
  of $N'=N\otimes_RR'$ into eigenspaces of $b$: $N'=\bigoplus N'_\lambda$, where $b$ acts as
  the scalar $\lambda$ on $N'_\lambda$. As in
  lemma \ref{Eigenwertzerlegung}(b) we can see, that $p(N'_\lambda,N'_\mu)=0$ unless $\lambda\cdot\mu=1$.
  This implies that $p$ is unimodular on $N'_1$ and therefore on $N_1$, thus $N$ is the
  direct sum of $N_1$ and the $p$-orthogonal complement $N_*$ of $N_1$. Since $b$ is a symplectic
  transformation, it leaves $N_*$ invariant. By the orthogonality relations for the $N_\lambda$ we
  have $N_*\otimes_RR'=\bigoplus_{\lambda\ne 1}N'_\lambda$. Since $\lambda-1\in (R')^\times$ for
  $\lambda\ne 1$ the endomorphism $b-1$ of $N_*$ induces an automorphism of $N_*\otimes_RR'$ and
  is therefore itself an automorphism of $N_*$,\qed

\begin{lemma} \label{Zentralisatoren}
  Let $h=p+q\in\GL_n(R)$ be $R$-$\Theta$-semisimple.
  Let $G^{h,\Theta}(R)=\{g\in \GL_n(R)| \tran g\cdot h\cdot g=h\}$. Then the following holds:
  \begin{itemize}
  \item[(a)] With the notations of lemma \ref{+-0Zerlegung} and of lemma \ref{Cayleytransformation} we have
    \begin{align*}G^{h,\Theta}(R) &= \ORTH(q_+,R)\times \Sp(p_-,R)\times \left(\Sp(p_0,R)\cap \ORTH(q_0,R) \right)\\
            &\cong \ORTH(q_+,R)\times \left(\Sp(p_-\oplus p_0,R)\cap \ORTH(q_-\oplus q_0,R) \right)\\
            &\cong \ORTH(q_+,R)\times Cent\left(C(q_-\oplus q_0),\Sp(p_-\oplus p_0,R) \right).\\
    \end{align*}
  \item[(b)] In the situation and with the notations of lemma \ref{Eigenwertzerlegung} we have moreover
    \begin{align*} \quad\left(\Sp(p_0,R')\cap \ORTH(q_0,R') \right)&=
      \left\{\left. (\phi_\lambda)\in\prod_{\lambda\ne \pm 1}GL(M'_\lambda)\right|
      \phi_{\lambda^{-1}}=\tran\phi_\lambda \text{ for all }\lambda\right\} \\
       &\cong\prod_{\lambda\in\LLL}\GL(M'_\lambda)
    \end{align*}
    where $\phi_{\lambda^{-1}}=\tran\phi_\lambda$ means that
    $ b_h(\phi_{\lambda^{-1}}m_{\lambda^{-1}},
    \phi_\lambda m_\lambda)=b_h(m_{\lambda^{-1}},m_\lambda)$ for all
    $ m_{\lambda^{-1}}\in M'_{\lambda^{-1}},m_\lambda\in M'_\lambda$ and where
    $\LLL$ denotes a subset of the set of all $\lambda\ne\pm1$, which
    takes from every pair $\{\lambda,\lambda^{-1}\}$ exactly one member.
  \item[(c)] $ \left(\Sp(p_-\oplus p_0)\cap \ORTH(q_-\oplus q_0) \right)\cong
    Cent\left(C(q_-\oplus q_0),\Sp(p_-\oplus p_0) \right)$
    is a connected reductive smooth group scheme $/R$ with connected special fiber, which becomes split
    over the finite \'{e}tale extension $R'/R$.
  \item[(d)] We have in the situation of \ref{explizitesMatchingfuerresiduellhalbeinfach}
    \begin{align*} Cent(\NNN(h),\Sp_{2n})&\cong \Sp_{2(n-g)}\times
    Cent\left(C(q_-\oplus q_0),\Sp(p_-\oplus p_0)\right)
    \end{align*}
    where $2g$ is the rank of $M_-\oplus M_0$.
  \item[(e)] To obtain the intersections of $G^{h,\Theta}(R)$ with $\SL_n(R)$
    one has only to replace $\ORTH(q_+,R)$ by $\SO(q_+,R)$ on the right hand sides of (a).
  \end{itemize}
\end{lemma}
Proof: (a) Since every $g\in G^{h,\Theta}(R)$ stabilizes the decomposition of
  lemma  \ref{+-0Zerlegung} one immediately gets the first two isomorphisms.
  The last one follows from lemma \ref{Cayleytransformation}(c).
  \klabsatz
(b) Every $g\in G^{h,\Theta}(R)$ centralizes $N_l(h)$ and therefore has to respect the
  decomposition of $M_0\otimes_RR'$ in eigenspaces of $N_l(h)$. The first description of
  $\Sp(p_0,R')\cap \ORTH(q_0,R')$ follows now from \ref{Eigenwertzerlegung}(b).
  Since $b_h$ is unimodular on $M'_{\lambda^{-1}}\oplus M'_\lambda$ it induces an identification
  of $M'_{\lambda^{-1}}$ with the dual space of $M'_\lambda$. This means that $\phi_\lambda$ can
  vary through the whole $\GL(M'_\lambda)$, while $\phi_{\lambda^{-1}}$ is then uniquely determined
  as the inverse of its adjoint. \klabsatz
  We remark that the condition $\phi_{\lambda}=\tran\phi_{\lambda^{-1}}$
  is equivalent to the condition $\phi_{\lambda^{-1}}=\tran\phi_\lambda$ and gives no extra
  restrictions. This is clear since we have $b_h(m_{\lambda},m_{\lambda^{-1}})=
  b'_h(m_{\lambda^{-1}},m_\lambda)=\lambda\cdot b_h(m_{\lambda^{-1}},m_\lambda)$  for
       $ m_{\lambda^{-1}}\in M'_{\lambda^{-1}},m_\lambda\in M'_\lambda$ by (\ref{Gleichungmu}),
       so the two possible identifications of $M'_{\lambda^{-1}}$ with the dual of
        $M'_\lambda$ differ by a scalar and create the same adjoint.
   The last isomorphism follows.
  \klabsatz
(c) This follows from (a) and (b).
  \klabsatz
(d) This follows from the definition of $\NNN$ by the remark, that an element of $Cent(b,\Sp_{2n}(R))$
  has to respect the decomposition of lemma \ref{Abspaltung1Sp}.
  \klabsatz
(e) is clear, since symplectic transformations have determinant $1$.
\qed

%
\begin{lemma}\label{SpkonjbereitsueberR}
  Let $G/R=\OOO_F$ be a connected reductive group with connected special fiber $G\times_{\OOO_F}\kappa$
  and $b\in G(R)$ be $R$-semisimple.\par
   If $b'=h_F^{-1}\cdot b\cdot h_F \in G(R)$
  for some $h_F\in G(\bar{F})$ then there exists $h_R\in G(R)$ with $b'=h_R^{-1}\cdot b\cdot h_R$.
\end{lemma}
Proof: This follows from \cite[Prop. 7.1.]{K-elliptic}.
  \qed

\begin{lemma}\label{GlkonjbereitsueberR}
  Let $R=\OOO_F$ and $h\in \GL_n(R)$ be $R$-$\Theta$-semisimple and
   $h'=\tran g_F\cdot h\cdot g_F \in \GL_n(R)$ for some
   $g_F\in \GL_n(\bar{F})$. Then we have:
\begin{itemize}
\item[(a)] If  additionally  $\det(g_F)\in F^\times$ there exists $g_R\in \GL_n(R)$ with
  $h'=\tran g_R\cdot h\cdot g_R$.
\item[(b)] If we only assume $g_F\in \GL_n(\bar{ F})$ and if $n$ is odd there exist $g_R\in \GL_n(R)$ and
  $\epsilon\in\OOO_F^\times$ such that $h'=\epsilon\cdot\tran g_R\cdot h\cdot g_R$.
\item[(c)]
  We get the statement of (a) if we additionally assume that the discriminants of
  $q_+$ and $q'_+$ coincide in $R^\times/(R^\times)^2$.
\item[(d)] Under the additional conditions $h,h'\in \SL_n(\OOO_F), \enskip g_F\in \SL_n(\bar{F)}$ and
  $n$ odd  we can find $g_R\in \SL_n(R)$ with $h'=\tran g_R\cdot h\cdot g_R$.
\end{itemize}
\end{lemma}
Proof: We use the objects occurring in lemma \ref{+-0Zerlegung} for $h$ and denote the corresponding
  objects for $h'$ by a $'$. We have $rank(M_+)=\dim(M_+\otimes_RF)=\dim(M'_+\otimes_RF)=rank(M'_+)$.
  By transforming $h$ and $h'$ with elements of $\GL_n(R)$ we can therefore assume (using lemma
  \ref{symplektischeNormalform}) that
  \begin{gather}\label{Standardisierung von hh'}
   M_+=M'_+=R^m,\quad M_0\oplus M_-=M'_0\oplus M'_-,\quad p_*:= p_0\oplus p_-=p'_0\oplus p'_-.
  \end{gather}
  The assumption and lemma \ref{Cayleytransformation}(c) (applied in the case $R=\bar{F}$) now imply that
  the elements $C(0\oplus q_0)$ and  $C(0\oplus q'_0)$ of $\Sp(p_*,R)$ are conjugate by an element
  of $\Sp(p_*,\bar{F})$.
  By lemma \ref{SpkonjbereitsueberR} they are conjugate by an element $g_*\in Sp(p_*,R)$,
  hence we get from lemma \ref{Cayleytransformation}(c) the equality $q'_0=\tran g_*\cdot q_0\cdot g_*$
  and therefore $p'_0+p'_-+q'_0=\tran g_*(p_0+p_-+q_0)g_*$ in $M_0\oplus M_-$.
  We have $\det(q'_+)=\det(h')\cdot\det(p'_0+p'_-+q'_0)^{-1}=\det(g_F)^2\cdot\det(h)\cdot
  \det(p_0+p_-+q_0)^{-1}=\det(g_F)^2\cdot\det(q_+)$ (observe $\det g_*=1$).
\absatz
  If case (a) we conclude using
  $R^\times\cap (F^\times)^2=(R^\times)^2$ and lemma \ref{orthogonaleNormalform}, that $q'_+$ and
  $q_+$ are transformed via an element $g_+\in \GL_n(M_+)$, a statement which has been an
  additional assumption in (c) in view of lemma \ref{orthogonaleNormalform}.
    We put $g_*$ and $g_+$ together to
  $g_R\in\GL_n(R)$ which does the required job in cases (a) and (c).
\absatz
  We prove (b) for $\epsilon=\det(q'_+)/\det(q_+)$: We have $h":=\epsilon^{-1}h'=\tran g'_F\cdot h\cdot g'_F$
  for $g'_F=\sqrt{\epsilon^{-1}}\cdot g_F\in \GL_{n}(\bar{F})$. If $2r+1$ is the rank of $M'_+$ we have
  $\det(q"_+)=\det(q'_+)\cdot \epsilon^{2r+1}=\det(q_+)\cdot \epsilon^{2r}$. Thus
   the additional assumption of (c) is fulfilled and we get $g_R\in \GL_n(R)$
   with $h"=\tran g_R \cdot h \cdot g_R$.
\absatz
  To prove (d) observe at first that we can assume the matrices transforming
  $h$ and $h'$ into the standard form (\ref{Standardisierung von hh'}) being in $\SL_n(R)$ since one
  can modify them by elements of $\GL(M_+)$ and since $\rank(M_+)\ge 1$. From $\det(g_F)=1$ we furthermore
  get $\det(q'_+)=\det(q_+)$ and therefore $\det(g_+)=\pm 1$. Since  we can replace
  $g_+$ by $-g_+$ if necessary and $\rank(M_+)$ is odd we can achieve $\det(g_+)=1$ and therefore
  $\det(g_R)=1$.
\qed

  \Seitenumbruch
 \meinchapter{Orbital integrals}

\Numerierung\label{DefOrbint}
For a (not necessarily connected) reductive group $\tilde G/\OOO_F$ with connected component
$G=\tilde G^\circ$ and elements
$\gamma\in \tilde G(F)$, $f\in \CCC^\infty_c(\tilde G(F))$ we define the orbital integral by:
\begin{gather*} O_\gamma(f,\tilde G(F))= \int_{G(F)/G(F)^\gamma} f(x\gamma x^{-1}) dx/dx^\gamma
\end{gather*}
where $G(F)^\gamma$ denotes the centralizer of $\gamma$ in $G(F)$ and where we have chosen
Haar measures $dx$ resp. $ dx^\gamma$ on $G(F)$ resp. $G(F)^\gamma$ such that
\begin{gather*}
vol_{dx}(G(\OOO_F)) =1\quad\text{ and }\quad vol_{dx^\gamma}((G^\gamma)^\circ(\OOO_F))=1.
\end{gather*}
If $1_K$ denotes
the characteristic function of a compact open subset $K\subset G(F)$, we will use the following
abbreviation:
\begin{gather*} O_\gamma(1,\tilde G)\quad =\quad O_\gamma(1_{\tilde G(\OOO_F)}, \tilde G(F))
\end{gather*}
We further introduce stable orbital integrals
\begin{align*} \STO_\gamma(f,\tilde G(F))\quad&=\quad\sum_{\gamma'\sim \gamma} O_{\gamma'}(f,\tilde G(F))
\qquad\text{ respectively}\\
\STO_\gamma(1,\tilde G)\quad&=\quad\sum_{\gamma'\sim \gamma} O_{\gamma'}(1_{\tilde G(\OOO_F)},\tilde G(F))
\end{align*}
where $\gamma'$ runs through a set of representatives for the conjugacy classes inside the
stable conjugacy class of $\gamma$.
\Absatz\Numerierung
  Recall the construction of the {\bf quotient measure} $dg/dh$ on $G/H$ for totally disconnected locally
  compact groups $H\subset G$, where $G$ and $H$ are unimodular
  (e.g. $G$ and $H$ are the sets of $F$-valued points of reductive groups).
  One defines
  \begin{align*}
     vol(K\gamma H/H)\quad=\quad \int_{G/H} 1_{K\gamma H/H}(g) dg/dh \quad := \quad
     \frac{vol_{dg}(K)}{vol_{dh}(\gamma^{-1}K\gamma\cap H)},
  \end{align*}
  where $K\subset G$ is any open compact subgroup, and extends this by linearity to the space of all
  locally constant compactly supported functions on $G/H$.
\absatz
  Of course one has to prove a compatibility condition, if $K'\subset K$ is another
  open compact subgroup: For $\gamma\in G$ let
  \begin{gather*}\label{DoppelnebenklassenKK}
     K\gamma H\quad= \quad\bigcup_{j\in J}^{.}\; K'\cdot \gamma_j\cdot\gamma\cdot H
  \end{gather*}
  be a disjoint double coset decomposition with $\gamma_j\in K$. We have to prove:
  \begin{gather}\label{Maassinvarianz}
    vol(K\gamma H/H)\quad =\quad \sum_{j\in J}\; vol(K'\gamma_j\gamma H/H).
  \end{gather}
  Define $C_j:=\{ K'y\in K'\backslash K\mid K'y\gamma\subset K'\gamma_j\gamma H\}\subset
  K'\backslash K$ for $j\in J$.
  Then we have a disjoint decomposition
  \begin{align*}
    K'\backslash K\quad =\quad \bigcup_{j\in J}^{.} C_j\qquad\qquad&\text{ and isomorphisms}\\
    i_j:\quad (\gamma^{-1}\gamma_j^{-1}K'\gamma_j\gamma\cap H)\backslash(\gamma^{-1}K\gamma\cap H)
      \quad &\xrightarrow{\ \sim \ }\quad C_j \\
    h\qquad &\mapsto \quad K'\cdot\gamma_j\cdot\gamma\cdot h\cdot \gamma^{-1}.
  \end{align*}
  This implies
  \begin{align*}
    \frac{vol_{dg}(K)}{vol_{dg}(K')}\quad=\quad[ K:K']\quad&=\quad
      \sum_{j\in J}  [ (\gamma^{-1}K\gamma\cap H):(\gamma^{-1}\gamma_j^{-1}K'\gamma_j\gamma\cap H) ] \\
    \quad&=\quad \sum_{j\in J}
      \frac{vol_{dh}(\gamma^{-1}K\gamma\cap H)}{vol_{dh}(\gamma^{-1}\gamma_j^{-1}K'\gamma_j\gamma\cap H)}
  \end{align*}
  or equivalently
  \begin{align*}
    \frac{vol_{dg}(K)}{vol_{dh}(\gamma^{-1}K\gamma\cap H)}\quad=\quad
    \sum_{j\in J}\frac{vol_{dg}(K')}{vol_{dh}(\gamma^{-1}\gamma_j^{-1}K'\gamma_j\gamma\cap H)}\, ,
  \end{align*}
  which is the claim (\ref{Maassinvarianz}).
\Absatz
  The crucial statement we need in the following is the following type of a fundamental lemma:
\begin{conj}\label{BCVermutung}
  If the regular algebraically semisimple and topologically unipotent elements
  $u\in \SO_{2n+1}(F)$ and $v\in \Sp_{2n}(F)$ are $BC$-matching (see \ref{BCMatching}) then
  \begin{gather}\tag{$BC_n$}
    \STO_u(1,\SO_{2n+1})=\STO_v(1,\Sp_{2n}).
  \end{gather}
\end{conj}
  The (easy) case ($\text{BC}_{\text 1}$) is proved in \cite[Stable case I in Proof of Theorem]{FPGL3}.
  The case ($\text{BC}_{\text 2}$) is essentially proved in \cite[Part II]{FGL4}, as will be explained
  in \ref{BC2Theorem}.
\Absatz
  {\bf Warning:} While ($\text{BC}_{\text 1}$) is an immediate consequence of the exceptional
  isogeny $i_2:\SP_{2}=\SL_2\twoheadrightarrow \PGL_2=\SO_3$ and the fact, that $\gamma^2$ and
  $i_2(\gamma)$ are $BC$-matching for $\gamma\in\SL_2(F)$, the statement
  ($\text{BC}_{\text 2}$) is much deeper, since the exceptional isogeny $i_4:\SP_4\twoheadrightarrow\SO_5$
  does not satisfy the analogous matching property.
\begin{remark}\label{Existenz von BC matching elements}\em
   It follows immediately from the construction in \ref{BCMatching} that we have a bijection between
   $F$-rational conjugacy classes in $\SO_{2n+1}(\bar F)$ and in $\SP_{2n}(\bar F)$. By the theorem of
   Steinberg each $F$-rational conjugacy class in $\SP_{2n}(\bar F)$ contains a rational element, since
   $\SP_{2n}$ is quasisplit and simply connected. But the same statement holds for $F$-rational
   topologically unipotent conjugacy classes in $\SO_{2n+1}(\bar F)$ as well:
\klabsatz
   If $u\in \SO_{2n+1}(\bar F)$ is topologically unipotent
   and represents an $F$-rational conjugacy class, consider its two preimages $v_1$ and $v_2=-v_1$ in
   $\Spin_{2n+1}$. Since $p\ne 2$ we have $\lim_{n\to\infty}v_2^{p^{n}}=-\lim_{n\to\infty}v_1^{p^{n}}$,
   so that exactly one of the elements $v_1,v_2$ is topologically unipotent, say $v_1$. Since the Galois group respects
   the property to be topologically unipotent, the conjugacy class of $v_1$ is $F$-rational and
   therefore contains an $F$-rational element $v'$ by the theorem of Steinberg. The image of $v'$ in
   $\SO_{2n+1}( F)$ is the desired $F$-rational representative of the conjugacy class of $u$.
\klabsatz
   Thus to every topologically unipotent element in $v\in\SP_{2n}(F)$ is associated at least one
   $BC$-matching $u\in\SO_{2n+1}(F)$ and vice versa.
\end{remark}

\begin{lemma}[Kazhdan-Lemma]\label{Kazhdan-Lemma}
$\ $
\begin{itemize}
\item[(a)] For $\tilde G=G\rtimes \langle\Theta\rangle$ as in \ref{splitmitauto} let us assume
  that the following statement holds:
  \begin{itemize}
    \item[($*$)]  If $s_1\Theta$ and $s_2\Theta$ for $s_1,s_2\in G(\OOO_F)$ are residually semisimple and conjugate
    by an element of $G(F)$ then they are also conjugate by an element of $G(\OOO_F)$.
  \end{itemize}
  If  $\gamma\Theta=u\cdot s\Theta=s\Theta\cdot u$ is a topological Jordan decomposition, where
  $\gamma\in G(\OOO_F)$,  $u$ is topologically unipotent and $s\Theta$ residually semisimple, we have
  \begin{gather*} O_{\gamma\Theta}(1,\tilde{G})\quad =\quad
  \frac{1}{[G^{s\Theta}(\OOO_F):(G^{s\Theta})^\circ(\OOO_F)]}\cdot O_u(1,G^{s\Theta})
  \end{gather*}
\item[(b)] Let $H/\OOO_F$ be connected reductive with connected special fiber.
  For $h\in H(\OOO_F)$  with  topological Jordan decomposition
  $h=v\cdot b=b\cdot v$, where $v$ is topologically unipotent and $b$ residually semisimple, we have
  \begin{gather*} O_h(1,H)\quad =\quad
  \frac{1}{[H^{b}(\OOO_F):(H^{b})^\circ(\OOO_F)]}\cdot O_v(1,H^b)
  \end{gather*}
\end{itemize}
\end{lemma}
Proof: (a) We first prove:
\klabsatz
\begin{itemize}
\item[(**)]
  We have $g\gamma\Theta g^{-1}\in G(\OOO_F)\Theta$ if and only if $g$ is of the form $g=k\cdot x$
  where  $k\in G(\OOO_F)$ and $x\in G^{s\Theta}(F)$ satisfies $xux^{-1}\in G^{s\Theta}(\OOO_F)$.
\end{itemize}
\klabsatz The direction "{}$\Leftarrow $"{} is easy: Under the hypothesis we have
  $g\gamma\Theta g^{-1}=kxus\Theta x^{-1}k^{-1}=
  k(xux^{-1})(s\Theta)k^{-1}\in G(\OOO_F)$.
  For the converse direction "{}$\Rightarrow$"{} let us assume that $g\gamma\Theta g^{-1}\in G(\OOO_F)\Theta$.
  The topological Jordan decomposition is $g\gamma\Theta g^{-1}=(gug^{-1})\cdot (gs\Theta g^{-1})$.
  Since $\langle G(\OOO_F),\Theta\rangle$ is a closed subgroup of $\widetilde G(F)$ we conclude from
  \ref{top Jordan Eigenschaften}(4) that $gs\Theta g^{-1}\in G(\OOO_F)\Theta$ and $gug^{-1}\in G(\OOO_F)$.
  By the first inclusion and assumption ($*$)  we get an element $k\in G(\OOO_F)$ such that
  $g(s\Theta)g^{-1}=k(s\Theta)k^{-1}$, which implies $x=k^{-1}\cdot g\in G^{s\Theta}(F),$
  where $G^{s\Theta}$ is the centralizer
  of $s\Theta$ in $G$. Using $g=kx$ the inclusion $gug^{-1}\in G(\OOO_F)$ is now equivalent to
  $xux^{-1}\in G(\OOO_F)$, which proves $(**)$.
\absatz
  To finish the proof we introduce the double coset decomposition
  \begin{gather*}
   \{g\in G(F)\mid g\gamma\Theta g^{-1}\in G(\OOO_F)\Theta\}\quad=\quad
   \bigcup_{i\in I}^{.}\; G(\OOO_F)\cdot g_i\cdot G^{\gamma\Theta},
  \end{gather*}
  where we can assume $g_i\in G^{s\Theta}(F)$ in view of $(**)$. Again from $(**)$ we
   get the double coset decomposition
  \begin{gather*}
   \{x\in G^{s\Theta}(F)\mid xux^{-1}\in G^{s\Theta}(\OOO_F)\}\quad=\quad
   \bigcup_{i\in I}^{.}\; G^{s\Theta}(\OOO_F)\cdot g_i\cdot G^{\gamma\Theta},
  \end{gather*}
  so that it remains to prove
  \begin{align}\label{VolumengleichungGsTheta}
   \sum_{i\in I}&\frac{vol_{dg}(G(\OOO_F))}{vol_{dh}(g_i^{-1}\cdot G(\OOO_F)\cdot g_i\cap G^{\gamma\Theta}(F))}\\
   \notag &=\quad\frac{1}{[G^{s\Theta}(\OOO_F):(G^{s\Theta})^\circ(\OOO_F)]}\cdot\sum_{i\in I}\;
   \frac{vol_{d\eta}(G^{s\Theta}(\OOO_F))}
   {vol_{dh}(g_i^{-1}\cdot G^{s\Theta}(\OOO_F)\cdot g_i\cap G^{\gamma\Theta}(F))}\; ,
  \end{align}
  where $d\eta$ is a Haar measure on $G^{s\Theta}(F)$ satisfying $vol_{d\eta}((G^{s\Theta})^\circ(\OOO_F))=1$.
  This implies $vol_{d\eta}(G^{s\Theta}(\OOO_F))= [G^{s\Theta}(\OOO_F):(G^{s\Theta})^\circ(\OOO_F)]$.
  On the other hand we claim
  \begin{gather}\label{GGsThetaSchnittGgammaThetagleichung}
   g_i^{-1}\cdot G(\OOO_F)\cdot g_i\cap G^{\gamma\Theta}(F)=
   g_i^{-1}\cdot G^{s\Theta}(\OOO_F)\cdot g_i\cap G^{\gamma\Theta}(F).
  \end{gather}
  The inclusion "{}$\supset$"{} being trivial let us assume that $g=g_i^{-1}\cdot \sigma\cdot g_i$ is
  an element of the left hand side. But $g_i\in G^{s\Theta}(F)$ and
  $g\in G^{\gamma\Theta}(F)\subset G^{s\Theta}(F)$ imply $\sigma\in G^{s\Theta}(F)$.
  Since $G^{s\Theta}(F)\cap G(\OOO_F)=G^{s\Theta}(\OOO_F)$ we get that $g$ lies in the right hand
  side, i.e. (\ref{GGsThetaSchnittGgammaThetagleichung}) is proved. (\ref{VolumengleichungGsTheta})
  now follows immediately.
\Absatz
(b) is now clear: We have $\tilde G=G$ i.e. $\Theta=1$ and the assumption ($*$) is satisfied
by lemma \ref{SpkonjbereitsueberR}.
\qed
\Absatz

 The following lemmas will be useful in later chapters.

\begin{lemma}\label{Homogenitaet}
  If $N\in \NN$ is prime to $p$ then we have for a reductive group $G/\OOO_F$ and $\gamma\in G(F)$
  \begin{gather*} O_{\gamma^N}(1,G)= O_\gamma(1,G).
  \end{gather*}
\end{lemma}
Proof: Notice that $g\cdot\gamma\cdot g^{-1}$ lies in the closure of
  $(g\cdot\gamma^N\cdot g^{-1})^\ZZ$ if $N\in \ZZ_p^\times$.
  This gives the equivalence
  $g\cdot\gamma^N\cdot g^{-1}\in G(\OOO_F)\Longleftrightarrow g\cdot\gamma\cdot g^{-1}\in G(\OOO_F)$,
  which implies the identity of orbital integrals.
\qed

\begin{lemma}\label{EndlichesZentrum}
  If $G/\OOO_F$ is of the form $G=G_1\times Z$ with a reductive group $G_1$ and a finite group
  $Z\simeq Z(\OOO_F)$ then we have for $\gamma\in G_1(F)\subset G(F)$ the following identity of
  orbital integrals:
  \begin{gather*}
    O_\gamma(1,G)\quad =\quad O_\gamma(1,G_1).
  \end{gather*}
\end{lemma}
Proof: We have $G(F)/G(F)^\gamma\simeq G_1(F)/G_1(F)^\gamma$ since $Z\subset G^\gamma$,
and the normalized Haar-measures on $G_1(F)$ and $G_1^\gamma(F)$ are the restrictions of
the normalized Haar-measures on $G(F)$ and $G^\gamma(F)$, since $G^\circ=G_1^\circ$ and
$(G^\gamma)^\circ=(G_1^\gamma)^\circ$. The claim follows.
\qed

\begin{lemma} \label{Redorbintaufadjungiert}
  Let $1\to T\to G\to H\to 1$ be an exact sequence of algebraic groups over $\OOO_F$
  where $T$ is a split torus. Then we have for $\gamma\in G(F)$ with image $\eta\in H(F)$:
  \begin{gather*}
    \STO_\gamma(1,G)=\STO_\eta(1,H).
  \end{gather*}
\end{lemma}
Proof: We use the fact that the image of $(G^\gamma)^\circ$ in $H$ is $(H^\eta)^\circ$.
By Hilbert 90 we get exact sequences
$1\to T(F)\to G(F)\to H(F)\to 1$ and $1\to T(F)\to (G^\gamma)^\circ(F)\to (H^\eta)^\circ(F)\to 1$,
so that we have an isomorphism $G(F)/(G^\gamma)^\circ(F)\simeq H(F)/(H^\eta)^\circ(F)$.
Since $(H^\eta)^\circ$ has finite index in $H^\eta$ we can compute $O_\eta(1,H)$ as
$\int_{H(F)/(H^\eta)^\circ(F)} 1_{H(\OOO_F)}(h\eta h^{-1}) dh/dh^\eta$.
Similarly $$O_\gamma(1,G)=\int_{G(F)/(G^\gamma)^\circ(F)} 1_{G(\OOO_F)}(g\gamma g^{-1}) dg/dg^\gamma.$$
Now the quotient measures on $G(F)/(G^\gamma)^\circ(F)$ and $H(F)/(H^\eta)^\circ(F)$ coincide since
$G(\OOO_F)\twoheadrightarrow H(\OOO_F)$ and
$(G^\gamma)^\circ(\OOO_F)\twoheadrightarrow (H^\eta)^\circ(\OOO_F)$, and
 we conclude $O_\gamma(1,G)=O_\eta(1,H)$.
\absatz
It remains to check that the
set $St_\gamma$ of conjugacy classes inside the stable conjugacy class of $\gamma$ maps
bijectively to the corresponding set $St_\eta$ associated to $\eta$. But in  the
following commutative diagram of abelian groups with exact rows and columns the map $\iota$ must
be an isomorphism:
\[\begin{CD}
&&&&  H^1(F,T)  @= H^1(F,T)\\
&&&&    @VVV @VVV\\
1 @>>> St_\gamma @>>> H^1_{ab}(F,G^\gamma) @>>> H^1_{ab}(F,G)\\
&& @V\iota VV  @VVV       @VVV\\
1 @>>> St_\eta@>>> H^1_{ab}(F,H^\gamma) @>>>  H^1_{ab}(F,H)\\
& &&& @VVV @VVV \\
&&&& H^2(F,T)  @= H^2(F,T).
\end{CD}\]
Here $H^1_{ab}(F,.)$ denotes the abelianized cohomology of \cite{Borovoi} which coincides for
nonarchimedean $F$ as a pointed set with the usual cohomology.
\qed
 
   \Seitenumbruch
   
 \meinchapter{Comparison between $\PGL_{2n+1}$ and $\SP_{2n}$ }
  \Absatz
  Recall (see \ref{R theta semisimple}) that $R$ is either a field of characteristic $0$
  or the  ring of integers $\OOO_F$ of
   a local non archimedean field $F$ with residue characteristic $\ne 2$.
  \absatz

\Numerierung \label{explizitesMatchingfuerresiduellhalbeinfach} {\bf The explicit norm map $\NNN$.}
  Our final goal being the comparison of $\Theta$-twisted stable orbital integrals on $\PGL_{2n+1}$ with
  stable orbital integrals on $\SP_{2n}$, we will represent elements of $\PGL_{2n+1}$ by elements of
   the groups $\GL_{2n+1}$ resp. $\SL_{2n+1}$.
  Let $\GL_n(R)_{R\Theta ss}/traf$ resp. $\SL_n(R)_{R\Theta ss}/traf$
  be the set of transformation classes of $R$-$\Theta$-semisimple (see \ref{R theta semisimple})
  elements of $h\in\GL_n(R)$ resp. $h\in\SL_n(R)$ under the transformations $h\mapsto \tran ghg$
  for $g\in GL_n(R)$ resp $g\in SL_n(R)$. Similarly let $Sp_{2g}(R)_{Rss}/conj$ be the set of
  conjugacy classes of $R$-semisimple elements in $Sp_{2g}(R)$. We define a norm map
  $$\NNN: \GL_{2n+1}(R)_{R\Theta ss}/traf \longrightarrow \Sp_{2n}(R)_{Rss}/conj$$
  as follows: If $h=p+q\in \GL_{2n+1}(R)$ represents a class of the left hand side, we decompose
  $M=R^{2n+1}=M_+\oplus M_-\oplus M_0$ as in lemma \ref{+-0Zerlegung}. We consider $M_+$ as the degenerate part
  of $M$ with respect to $p$ and denote the non degenerate part by $M_*:=M_-\oplus M_0$. Since $p_*=p_-\oplus p_0$ is
  a unimodular form on $M_*$ we can find a basis $(e_1,\ldots,e_g,f_g,\ldots,f_1)$ of
  $M_*$ such that $p_*$ has standard form with respect to this basis
  by lemma \ref{symplektischeNormalform}. Let $P_*$ resp. $Q_*$ be the matrix describing the
  {(skew-)} symmetric bilinear
  form $p_*$ resp. $q_-\oplus q_0$ with respect to this basis ($q_-$ is the zero form). Thus $P_*=J_{2g}$ and
   $\Sp(P_*)=\Sp_{2g}$. Now $\NNN(h)$ or more precisely the image of the class of $h$ under the norm map $\NNN$
   is defined to be the $\Sp_{2n}(R)$-conjugacy  class of
  $1_{2(n-g)}\times C(Q_*)\in \Sp_{2(n-g)}(R)\times \Sp_{2g}(R)\subset \Sp_{2n}(R)$, where we use the
  Cayley-transform-map $C$ from lemma \ref{Cayleytransformation}.

\begin{remark}\label{andereBeschreibung Norm}\em
  In the situation where the decomposition $M=R^{2n+1}=M_+\oplus M_*$ is of the form
  $M=R^{2(n-g)+1}\oplus R^{2g}$ the matrix $h$ splits into the blocks $h_+\in \GL_{2(n-g)+1}(R)$ and
  $h_*\in \GL_{2g}(R)$ so that $N_l(h_*)= \tran h_*^{-1}\cdot h_*$ is a symplectic transformation with respect
  to the alternating part $p_*$ of $h_*$. Then $C(Q_*)\in \SP_{2g}(R)$ is the conjugate of $N_l(h_*)$
  by a matrix, which transforms $p_*$ into the standard form $J_{2g}$.
\end{remark}

\begin{prop}\label{NNN prop PGLSP}
Let $R$ be as above. Then the following statements hold:
\begin{itemize}
\item[(a)] The map $\NNN: \GL_{2n+1}(R)_{R\Theta ss}/traf \longrightarrow \Sp_{2n}(R)_{Rss}/conj$
  is well defined and surjective. In the case  $R=\OOO_F $
  its fibers are of order $2=\#(R^\times/(R^\times)^2)$ and describe the two different classes of
  unimodular quadratic forms on $M_+$.
\item[(b)] The restriction $\NNN_{\SL}$ of $\NNN$ to $\SL_{2n+1}(R)_{R\Theta ss}/traf$ is surjective as well.
  It is bijective if $R$ is an algebraically closed field or if $R=\OOO_F$.
\item[(c)] If $h$ represents a class in $\GL_{2n+1}(R)_{R\Theta ss}/traf$ then the image of
  $h\cdot J^{-1}\Theta$ in $ \PGL_{2n+1}(R)\rtimes \langle\Theta\rangle$
  matches with $\NNN(h)$ in the sense of $\Theta$-endoscopy.
\end{itemize}
\end{prop}
Proof: (a) and (b) The choices made in constructing $\NNN(h)$ only allow $Q_*$ to be replaced
  by some $\tran g\cdot Q_*\cdot g$ for $g\in \Sp(P_*,R)$. By lemma \ref{Cayleytransformation}(c) this
  does not change the conjugacy class of $\NNN(h)$. Therefore the map $\NNN$ is well defined. To prove
  surjectivity first observe that each class in $\Sp_{2n}(R)_{Rss}/conj$ has a representative of the
  form $(1_{2(n-g)},b)$ with $b\in \Sp_{2g}(R)_{ess}$ by lemma \ref{Abspaltung1Sp} with a unique
  $g\le n$. The
  $\Sp_{2g}(R)$-conjugacy-class of $b$ is unique. The bijectivity of the Cayley-transform map
  and property \ref{Cayleytransformation}(c) then imply that there is a $Q_*\in Sym_{2g}(R)$,
  which is unique up to transformations with elements of $\Sp_{2g}(R)=\Sp(P_*,R)$, such that $b=C(Q_*)$.
  Now we consider the unimodular bilinear form $h_*=P_*+Q_*$ on $R^{2g}$ and some unimodular symmetric
  bilinear form $q_+$ on $R^{2(n-g)+1}$. The form $q_+\oplus h_*$ on $R^{2n+1}$ is then unimodular
  and $R$-$\Theta$-semisimple. Since we can choose $q_+$ in such a way that $\det(q_+\oplus h_*')=1$
  we get the surjectivity statements of (a) and (b). Since the transformation class of $h_*'$ is unique
  by the considerations above and since $h=q_+\oplus h_*$
  we conclude that the fibers of $\NNN$ correspond to the transformation classes
  of unimodular quadratic forms on $M_+$. The remaining statements of (a) and (b) now follow from
  lemma \ref{orthogonaleNormalform}.
  \absatz
(c) By the definition of matching (\ref{DefMatching}) we can work over $R=\bar F$ and therefore may assume that
$\gamma=h\cdot J_{2n+1}^{-1}$
has diagonal form $\gamma=diag(t_1,\ldots,t_{2n+1})$. After applying a permutation in
$W_{\SO_{2n+1}}$ we may assume
\begin{gather}\label{Annahme an h}
 t_i\ne t_{2n+2-i} \text{ for }i\le g \text{ and }t_i=t_{2n+2-i} \text{ for }
g+1\le i\le 2n+1-g.
\end{gather}
 We have:
\begin{align*} h\quad&=\quad antidiag(t_1,-t_2,t_3,\ldots,t_{2n+1})\\
              h\pm\tran h\quad&=\quad antidiag(t_1\pm t_{2n+1},-t_2\mp t_{2n},t_3\pm t_{2n-1},
              \ldots,t_{2n+1}\pm t_1)\\
              \tran h^{-1}\cdot h\quad&=\quad
              diag(t_{2n+1}/t_1,t_{2n}/t_2,\ldots,t_{n+2}/t_n,1,t_n/t_{n+2},\ldots,t_1/t_{2n+1})
\end{align*}
This means that $M_+\simeq R^{2(n-g)+1}$ is spanned by the standard basis elements
 $e_{g+1},\ldots, e_{2n+1-g}$ of $R^{2n+1}$, and $M_*=M_-\oplus M_0$ by $e_1,\ldots,e_g,
 e_{2n+2-g},\ldots,e_{2n+1}$. Since  $h-\tran h$ is an antidiagonal matrix, its non degenerate part
  can be transformed by a diagonal matrix $d$ into the standard form $J_{2g}$.
  Now we use remark \ref{andereBeschreibung Norm} to get the following representative for
  $\NNN(h)$, observing that conjugation by $d$ does not change a diagonal matrix:
\begin{gather*}
diag(t_{2n+1}/t_1,t_{2n}/t_2,\ldots,t_{n+2}/t_n,t_n/t_{n+2},\ldots,t_1/t_{2n+1}),
 \end{gather*}
 which may be conjugated by an element of the Weylgroup into the form
 $$ diag(t_1/t_{2n+1},t_2/t_{2n},\ldots,t_n/t_{n+2},t_{n+2}/t_n,\ldots,t_{2n+1}/t_1).$$
 The claim now follows from example \ref{BspMatchingPGL2n+1Sp2n}.
 \qed

\begin{cor}\label{Existenz matching rational PGLSp}
  For every semisimple $\bar\gamma\Theta\in \widetilde{\PGL}_{2n+1}(F)$ there exists
  a semisimple $\eta\in \Sp_{2n}(F)$ matching with $\eta$ in the sense of \ref{DefMatching} and vice versa.
\end{cor}
Proof: If $\gamma\in \GL_{2n+1}(F)$ represents a given $\bar\gamma$ one applies part (c) of the proposition
  to $h=\gamma\cdot J_{2n+1}$. If $\eta$ is given one applies (b) and (c).
\qed

\begin{prop}\label{OrbintredvonPGLaufGlSl}
  Let $Z=Cent(\GL_{2n+1})\simeq \GG_m$ denote the center of $\GL_{2n+1}$.
  Let $\bar{\gamma}\in \PGL_{2n+1}(F)$ be represented by $\gamma\in \GL_{2n+1}(F)$. Since $2n+1$ is odd
  we can achieve that   $\det(\gamma)$ has even valuation.   Then
  \begin{align}\label{OrbintredvonPGLaufGL}
    O_{\bar\gamma\Theta}(1,\widetilde{\PGL}_{2n+1})\quad=
      \quad 2\cdot O_{\gamma\Theta}(1,\widetilde{\GL}_{2n+1}).
  \end{align}
  If moreover $\gamma\Theta$ is strongly compact with topological Jordan decomposition
  $\gamma\Theta=u\cdot(s\Theta)=(s\Theta)\cdot u$ we have $u\in \SL_{2n+1}(F)$ and get
  \begin{gather}\label{OrbintredvonPGLaufSL}
    O_{\bar\gamma\Theta}(1,\widetilde{\PGL}_{2n+1})\quad=\quad O_u(1,\SL_{2n+1}^{s\Theta})
  \end{gather}
\end{prop}
Proof:  The relation $\bar{g}\cdot\bar{\gamma}\Theta\cdot \bar{g}^{-1}\in \widetilde{\PGL_{2n+1}}(\OOO_F)$
 means $g\cdot \gamma\cdot\Theta(g)^{-1}=\zeta\cdot k$ with $\zeta \in Z(F)\simeq F^*$
 and $k\in \GL_{2n+1}(\OOO_F)$. Since $\det(\Theta(g))=\det(g)^{-1}$ the relation implies taking
 determinants
 \begin{gather}\label{detgleichung}
 \det(g)^2\cdot\det(\gamma)\cdot\zeta^{-2n-1}\in \OOO_F^*.
 \end{gather}
  This implies that $\zeta$  has even valuation $2m$ for $m\in\Z$, since the valuation of
  $\det(\gamma)$ was assumed to be even.
 If we replace $g$ by $g'=\zeta_\OOO\cdot\varpi^{-m}\cdot g$ for $\zeta_\OOO\in \OOO_F^*$
 we get $g'\cdot\gamma\cdot\Theta(g')^{-1}\in\GL_{2n+1}(\OOO_F)$. Conversely
 the equation (\ref{detgleichung}) implies that every $g'\in g\cdot Z(F)$ with this property
 must be of the stated form.
\klabsatz
 Next observe that the condition $\bar{g}\in \PGL_{2n+1}^{\gamma\Theta}(F)$ means that we
 have for some representative $g\in \GL_{2n+1}(F)$ of $\bar g$ and some $\zeta\in Z(F)\simeq F^*$ the
 relation $g\gamma\Theta(g)^{-1}=\zeta\gamma$. This implies the determinant equation:
 $\det(g)^2=\zeta^{2n+1}$. Putting $\rho=\det(g)/\zeta^n\in Z(F)$ this implies $\zeta=\rho^2$ and
 $\det(g)=\rho^{2n+1}$. If we replace $g$ by $\rho^{-1}\cdot g$ we get
 $g\gamma\Theta(g)^{-1}=\gamma$ and $\det(g)=1$ . The only other element in $g\cdot Z(F)$ having
 the first property is
 $-g$, but $\det(-g)=-1$. This means that we have isomorphisms
 \begin{align*}\GL_{2n+1}(F)^{\gamma\Theta}\quad&\xrightarrow{\ \sim\ }
 \quad\SL_{2n+1}(F)^{\gamma\Theta}\times\{\pm 1\} \qquad\text{ and }\\
 \SL_{2n+1}(F)^{\gamma\Theta}\quad&\xrightarrow{\ \sim\ }\quad\PGL_{2n+1}^{\bar{\gamma}\Theta}(F).
 \end{align*}
 Since the normalized Haar measure on $\PGL_{2n+1}(F)$ is the
 quotient of the normalized Haar measure on $\GL_{2n+1}(F)$ by the norma\-lized Haar measure on
 $Z(F)$ (i.e. $vol(Z(\OOO_F))=1$) and since the normalized measure on
 $\GL_{2n+1}(F)^{\gamma\Theta}$ restricts to the normalized Haar measure
  on $\PGL_{2n+1}^{\bar{\gamma}\Theta}(F)\simeq\SL_{2n+1}^{\gamma\Theta}(F) $,
  the above considerations imply the relation (\ref{OrbintredvonPGLaufGL}).
\klabsatz
 If $\gamma\Theta$ is strongly compact we can assume that $\gamma\in\GL_{2n+1}(\OOO_F)$ and apply
 lemma \ref{Kazhdan-Lemma} to get
 \begin{gather*}
    O_{\bar\gamma\Theta}(1,\widetilde{\PGL}_{2n+1})\quad=\quad O_u(1,\GL_{2n+1}^{s\Theta})
 \end{gather*}
 observing that $[\GL_{2n+1}^{s\Theta}(\OOO_F):(\GL_{2n+1}^{s\Theta})^\circ(\OOO_F)]=2$.
 But since  $\GL_{2n+1}(F)^{s\Theta}\simeq\SL_{2n+1}(F)^{s\Theta}\times\{\pm 1\}$
  we can apply lemma \ref{EndlichesZentrum} to conclude (\ref{OrbintredvonPGLaufSL}).
 \qed
\begin{lemma}\label{Matchinglemma}
  Let $h=sJ\in \GL_{2n+1}(\OOO_F)$ be $R$-$\Theta$-semisimple and
  $b=(1_{2(n-g)},b_*)\in \Sp_{2n}(\OOO_F)$
  a representing element of $\NNN(h)$ with $b_*\in \SP_{2g}(\OOO_F)_{ess}$.
  Since $M_+$ is of odd rank $2(n-g)+1$ we can identify $(M_+,q_+)$
  with $(\OOO_F^{2(n-g)+1}, \epsilon q_{sp})$ for some $\epsilon\in \OOO_F^\times$ and the standard
  splitform $q_{sp}$.
  Assume that we have $BC$-matching algebraically semisimple and topologically unipotent elements
  $$u_+\in \SO_{2(n-g)+1}(F)\simeq \SO(q_+ ) \quad\text{ and }\quad
  v_+\in \Sp_{2(n-g)}(F)\simeq \ker(b-1)(F)\cap \Sp_{2n}(F)$$
  and an additional algebraically semisimple and topologically unipotent element
  $$u_*\in \SO(q_*)(F)\cap \Sp(p_*)(F)\simeq Cent(b_*,\Sp_{2g}(F)).$$
  \neueZeile
  Then the elements
  $\gamma\Theta= s\Theta\cdot (u_+,u_*)=(u_+,u_*)\cdot s\Theta\in \PGL_{2n+1}(F)\Theta$ and
  $\eta:=(v_+^2,u_*^2)\cdot b=b\cdot (v_+^2,u_*^2)\in \Sp_{2n}(F)$ match.
\end{lemma}
 Proof: As in the proof of lemma \ref{explizitesMatchingfuerresiduellhalbeinfach}(c) we work
 in the case $F=\bar F$ and assume that $\gamma$ resp. $\eta$ lie in the diagonal tori. The same
 holds for the residually semisimple parts $s$ resp. $b$ and the topologically unipotent parts
 $u=(u_+,u_*)$ and $v=(v_+^2,u_*^2)$. As the matching of $s\Theta$ and $b$ is already proved
 in  \ref{explizitesMatchingfuerresiduellhalbeinfach}(c) we only have to examine the topologically
 unipotent elements. We can make the assumption (\ref{Annahme an h}) and write
 \begin{align*}
  u_+\quad &=\quad diag(w_{g+1},\ldots,w_n,1,w_n^{-1},\ldots,w_{g+1}^{-1})\in \SO_{2(n-g)+1}(\bar F)\\
 u_*\quad &=\quad diag(w_1,\ldots,w_g,w_g^{-1},\ldots,w_1^{-1})\in Cent(b_*,\Sp_{2g}(\bar F))
  \end{align*}
 By the definition of $BC$-matching we can assume
 \begin{gather*}
   v_+\quad =\quad diag(w_{g+1},\ldots,w_n,w_n^{-1},\ldots,w_{g+1}^{-1})\in \Sp_{2(n-g)}(\bar F)
 \end{gather*}
 Taking everything together we get from the description of $M_+$ and $M_*$ in the proof of
 lemma \ref{explizitesMatchingfuerresiduellhalbeinfach}(c):
 \begin{align*}
    u\quad &=\quad(w_1,\ldots,w_n,1,w_n^{-1},\ldots,w_1^{-1})\\
    v\quad &=\quad(w_1^2,\ldots,w_n^2,w_n^{-2},\ldots,w_1^{-2})
\end{align*}
and the claim follows again from example \ref{BspMatchingPGL2n+1Sp2n}.
 \qed
\Absatz
  The statement of the following theorem is the {\bf fundamental lemma} for semisimple elements in
  the stable endoscopic situation $(\Sp_{2n}, \widetilde{\PGL}_{2n+1})$.
   Recall that the fundamental lemma also predicts the vanishing
  of orbital integrals for those rational elements, which match with no rational elements on the other side.
  But in view of corollary \ref{Existenz matching rational PGLSp} this case does not occur.
\absatz
\begin{thm}\label{FulemmaSp2nPGL2n+1}
  If the semisimple elements $\bar\gamma\Theta\in \widetilde{\PGL}_{2n+1}(F)$ and $\eta\in \Sp_{2n}(F)$
  match in the sense of \ref{DefMatching} and if
  conjecture ($\text{BC}_{\text{m}}$) is true for all $m\le n$  then we have
  \begin{gather}\label{Fule2n2n+1Gleichung}
     \STO_{\bar\gamma\Theta}(1,\widetilde{\PGL}_{2n+1})=\STO_\eta(1, \Sp_{2n}).
  \end{gather}
\end{thm}
  Proof: \underline{Step 1} {\bf (Reductions):} In the first step we will prove that the nonvanishing of
  one side of (\ref{Fule2n2n+1Gleichung}) implies the nonvanishing of the other side and that we can
  reduce to the following situation:
  \begin{itemize}
    \item $\gamma\in \GL_{2n+1}(\OOO_F)$
    \item $\eta\in \Sp_{2n}(\OOO_F)$

  \item the topological Jordan decompositions are of the form
  $$\gamma\Theta=(u_+,u_*)\cdot s\Theta\quad\text{ and } \eta=(v_+^2,v_*^2)\cdot b\qquad\text{ such that}$$
    \begin{itemize}
      \item $b$ lies in $\NNN(h)$ where $h=s\cdot J_{2n+1}$,
      \item $u_+$ and $v_+$ are $BC$-matching,
      \item $u_*$ can be identified with $v_*$ under an isomorphism
      $Cent(b_*,\Sp_{2g}(\OOO_F))\simeq Aut(h_*)$
    \end{itemize}
  \end{itemize}
\absatz
  So let us assume that the right hand side of (\ref{Fule2n2n+1Gleichung})
  does not vanish. Then there exists $\eta'\in \Sp_{2n}(F)$ stably conjugate to $\eta$ which has a
  nonvanishing orbital integral, i.e. can be conjugated
  into $\Sp_{2n}(\OOO_F)$. We can assume that $\eta'\in \Sp_{2n}(\OOO_F)$ and that its
  topological Jordan decomposition satisfies $\eta'=b'\cdot v'=v'\cdot b$  with residually semisimple
  $b'=(1_{2(n-g)},b_*)\in \Sp_{2(n-g)}(\OOO_F)\times \Sp_{2g}(\OOO_F)$ and topologically unipotent
  $v'$. We write $v'$ in the form $v' =((v'_+)^2,(u'_*)^2)$ with
  $v'_+\in \Sp_{2(n-g)}(\OOO_F)$ and $u'_*\in Cent(b_*,\Sp_{2g}(\OOO_F))$ using \ref{top Jordan Eigenschaften}(4)
  and the general assumption $p\ne 2$. Thus we have nonvanishing $O_{h'}(1,\Sp_{2n})$ and get from
  the Kazhdan-lemma \ref{Kazhdan-Lemma} and lemma \ref{Homogenitaet}:
  \begin{align}\label{SPOrbitalreduktion}
   O_{\eta'}(1,\Sp_{2n})\quad &=\quad O_{v'}(1,Cent(b',\Sp_{2n}))\\
   \notag &=\quad O_{(v'_+)^2}(1,\Sp_{2(n-g)})\cdot O_{(u'_*)^2}(1,Cent(b_*,\Sp_{2g}))\\
   \notag &=\quad O_{v'_+}(1,\Sp_{2(n-g)})\cdot O_{u'_*}(1,Cent(b_*,\Sp_{2g})).
  \end{align}
  Hence the stable orbital integral $\STO_{v'_+}(1,\Sp_{2(n-g)})$
  (being the sum of integrals of nonnegative functions) is strictly
  positive.
\absatz
  By remark \ref{Existenz von BC matching elements} there exists a $BC$-matching between $v'_+$ and
   some $u'_+\!\in\! \SO_{2(n-g)+1}(F)$. Then the equation ($\text{BC}_{n-g}$)
  implies that there exists $u_+\in \SO_{2(n-g)+1}(F)$  with strictly
  positive orbital integral and $BC$-matching with  $v'_+$, i.e. we can assume $u_+\in \SO_{2(n-g)+1}(\OOO_F)$.
\absatz
  Let $h=sJ\in \GL_{2n+1}(\OOO_F)_{R\Theta ss}$
  be a residually semisimple element with $\NNN(h)=b'$ and define the element
  $\gamma'\Theta=(u_+,u'_*)\cdot s\Theta= s\Theta\cdot (u_+,u'_*)\in\widetilde{\GL_{2n+1}}(\OOO_F)$.
  Here we identify the $Cent(s\Theta,{\GL_{2n+1}}\simeq G^{h,\Theta}\simeq
  \ORTH(q_+,R)\times Cent\left(C(q_*),\Sp(p_*) \right)\simeq \SO_{2(n-g)+1}\times Cent(b_*,\Sp_{2g}$,
   so that $(u_+,u'_*)$ can be viewed as  an element of the left hand side.
  The element $\overline{\gamma'}\Theta\in \widetilde{\PGL_{2n+1}}(\OOO_F)$
  matches with $\eta'$ (and therefore also with $\eta$) by lemma \ref{Matchinglemma}
  and therefore lies in the stable conjugacy class of $\gamma\Theta$.
\absatz
  If the left hand side of
  (\ref{Fule2n2n+1Gleichung}) does not vanish, it is immediate that there exists $\gamma'\Theta\in
  \widetilde{\GL_{2n+1}}(\OOO_F)$ in the stable conjugacy class of $\gamma\Theta$. By reversing the
  above arguments we see that there exists $\eta'\in \Sp_{2n}(\OOO_F)$ in the stable class of $\eta$.
  So excluding the tautological case that (\ref{Fule2n2n+1Gleichung}) means $0=0$ we may assume without
  loss of generality that $\gamma\in \GL_{2n+1}(\OOO_F)$ and $\eta\in \Sp_{2n}(\OOO_F)$.
  We may furthermore assume that $\gamma\Theta=(u_+,u_*)\cdot s\Theta$ and $\eta=(v_+^2,u_*^2)\cdot b$
  are the topological Jordan decompositions with $BC$-matching $u_+$ and $v_+$ and matching
  residually semisimple $s\Theta$ and $b$.
 \absatz

\underline{Step 2} {\bf (Calculation of the symplectic orbital integral):}
  If $\eta'\in \Sp_{2n}(F)$ is stable conjugate to $\eta$ then the residually semisimple
  parts $b'$ and $b$ are stable conjugate as well. If $\eta'$ has nonvanishing orbital integral then
  $\eta'$ and therefore also $b'$ can be conjugated into $\Sp_{2n}(\OOO_F)$, i.e. we can assume
  $b'\in  \Sp_{2n}(\OOO_F)$. By the Kottwitz lemma \ref{SpkonjbereitsueberR}
  $b'$ and $b$ are conjugate over $\Sp_{2n}(\OOO_F)$ i.e. we can assume $b'=b$.
   This means that
  we obtain all relevant conjugacy classes  in the stable conjugacy class of $\eta$ if we let
  $v'_+$ vary through a set of representatives for the conjugacy classes inside the stable conjugacy class of $v_+$ in
  $\Sp_{2(n-g)}(F)$ and $u'_*$ through a set of representatives for the conjugacy classes
  inside the stable conjugacy class of $u_*$
  in $Cent(b_*,\Sp_{2g})$. Then the corresponding $\eta'$ are of the form
   $$\eta'= b\cdot ((v'_+)^2,(u'_*)^2).$$
  We get using (\ref{SPOrbitalreduktion}) and lemma \ref{Homogenitaet}:
\begin{align}\label{Sporbital}
  \STO_{\eta}(1,\Sp_{2n})\quad
          &=\quad \sum_{v'_+\sim v_+} O_{(v'_+)^2}(1,\Sp_{2(n-g)})
             \cdot \sum_{u'_*\sim u_*} O_{(u'_*)^2}(1,Cent(b_*,\Sp_{2g}))\\
  \notag  &=\quad \sum_{v'_+\sim v_+} O_{v'_+}(1,\Sp_{2(n-g)})
             \cdot \sum_{u'_*\sim u_*} O_{u'_*}(1,Cent(b_*,\Sp_{2g})).
\end{align}
\underline{Step 3} {\bf (Calculation of the $\Theta$-twisted orbital integral):}
  We can repeat this argument in the $\Theta$-twisted situation,
  since by lemma \ref{GlkonjbereitsueberR}(b)
  the class of the residually semisimple part $\bar{s}\Theta$ of $\bar{\gamma}\Theta$   is the only
   $\PGL_{2n+1}(F)$-conjugacy class inside the stable class of $\bar{s}\Theta$,
  which meets $\PGL_{2n+1}(\OOO_F)$. If we denote by $u'_+$ a set of representatives for the
  $\SO_{2(n-g)+1}(F)$-conjugacy classes in the stable class of $u_+\in \SO_{2(n-g)+1}(\OOO_F)$ we
  therefore get using proposition \ref{OrbintredvonPGLaufGlSl}
 \begin{align}\label{Glorbital}
 \STO_{\bar{\gamma}\Theta}(1, \widetilde{\PGL}_{2n+1})&\quad = \quad
  \sum_{(u'_+,u'_*)\sim (u_+,u_*)} O_{(u'_+,u'_*)}(1,\SL_{2n+1}^{s\Theta})\\
  \notag=\quad &\sum_{u'_+\sim u_+} O_{u'_+}(1,\SO_{2(n-g)+1})\cdot
  \sum_{u'_*\sim u_*} O_{u'_*}(1,Cent(b_*,\Sp_{2g})).
 \end{align}
\underline{Step 4} {\bf (End of the proof):}
  Since $v_+$ and $u_+$ are $BC$-matching it only remains to apply ($\text{BC}_{n-g}$)
  in order to identify
  $$\sum_{v'_+\sim v_+} O_{v'_+}(1,\Sp_{2(n-g)})\quad\text{ with }\quad
  \sum_{u'_+\sim u_+} O_{u'_+}(1,\SO_{2(n-g)+1}).$$
  Thus the right hand sides of  (\ref{Sporbital}) and (\ref{Glorbital}) coincide, and the proof of the theorem
  is finished.

\qed
  
  \Seitenumbruch
\meinchapter{Comparison between $\GL_{2n}\times \GL_1$ and $\GSPIN_{2n+1}$}

\begin{lemma}[Cayley transformation again]\label{Cayleytransformationorthogonal}
  For a symmetric matrix $q\in\GL_n(R)$ the following holds:
  \begin{itemize}
\item[(a)] We have a bijection
  \begin{gather*}
   \tilde C: Alt_n(R)_{q-ess} \rightarrow \ORTH (q,R)_{ess},\quad p\mapsto (p-q)^{-1}\cdot (q+p)=-N_l(p+q)
  \end{gather*}
  between the set $Alt_n(R)_{q-ess}$ of skew-symmetric matrices $p$ such that $p\pm q\in \GL_n(R)$
  and the set
  $\ORTH(q,R)_{ess}$ of orthogonal transformations $b$
  such that $b-1\in \GL_n(R)$. The inverse map
  is $\tilde C^{-1}: b\mapsto q\cdot (b+1)\cdot(b-1)^{-1}$.
\item[(b)] $\tilde C$ induces a bijection between those elements $q$ of $Alt_n(R)_{q-ess}$,
  for which $p+q$ is $R$-$\Theta$-semisimple, and the $R$-semisimple elements of $\ORTH (q,R)_{ess}$.
\item[(c)] The map $\tilde C$ satisfies $\tilde C(\tran g\cdot p \cdot g)=g^{-1}\cdot \tilde C(p)\cdot g$ for
  $g\in \ORTH(q,R)$.
\item[(d)] We have $\det(b)=(-1)^n$ for $b\in \ORTH(q,R)_{ess}$.
\end{itemize}\end{lemma}
Proof: (a) For $p\in Alt_n(R)_{q-ess}$ we put $h=p+q$ and $b=(p-q)^{-1}\cdot (q+p)=-\tran h^{-1}\cdot h$
  Adding the formulas (\ref{tranbhbisth}) and (\ref{tranbtranhbisttranh}) in the proof of lemma
  \ref{Cayleytransformation} we get $(-\tran b)\cdot q\cdot(- b)=q$, i.e. $b\in\ORTH(q,R)$.
\klabsatz
  Furthermore $b-1=(p-q)^{-1}\cdot((p+q)-(p-q))=(p-q)^{-1}\cdot 2q\in \GL_n(R)$
  by the assumptions. The map $\tilde C$ is therefore defined.
\klabsatz
  Conversely we get for $b\in \ORTH(q,R)_{ess}$ and
  $p=q\cdot(b+1)\cdot(b-1)^{-1}$ the equivalences:
  \begin{align*}
  p=-\tran p&\Leftrightarrow
  q\cdot(b+1)\cdot(b-1)^{-1}=\tran( b-1)^{-1}\cdot\tran(b+1)\cdot(-q)\\
  &\Leftrightarrow (\tran b-1)q(b+1)=(\tran b+1)q(1-b)\\
  &\Leftrightarrow \tran bqb +\tran bq-qb-q=-\tran b qb+\tran bq-qb+q\\&\Leftrightarrow \tran bqb=q
  \Leftrightarrow b\in \ORTH(q,R).
  \end{align*}
  Furthermore $p\pm q= q\cdot\left((b+1)\pm(b-1)\right)\cdot(b-1)^{-1}\in \GL_n(R)$
  since $(b-1)^{-1},2b,2,q\in \GL_n(R)$.
  Therefore the map $\tilde C^{-1}$ is also well defined. An easy calculation using the relation
  $(b+1)\cdot(b-1)^{-1}=(b-1)^{-1}\cdot(b+1)$
  shows that the maps $\tilde C$ and $\tilde C^{-1}$ are inverse to another in their domain of definition.
\absatz
  (b) and (c) follow as in the proof of lemma \ref{Cayleytransformation}.
\absatz
  (d) is clear since every $b\in \ORTH_n(R)$ with $\det(b)=(-1)^{n-1}$ has $1$ as an eigenvalue.
    (Alternatively we can use (a) and the computation $\det(-\tran h^{-1}\cdot h)=(-1)^n$.)
  \qed

\begin{lemma}\label{Abspaltung1O}
  If $q$ is a unimodular symmetric bilinear form on a free $R$-module $N$ and
   $b\in \ORTH(q,R)$ is $R$-semisimple
  then there exists a $b$-invariant $q$-orthogonal direct sum decomposition
  $N=N_1\oplus N_*$ such that $b$ acts as identity on $N_1$ and
  $b|N_* \in \ORTH(q_*,R)_{ess}$, where $q_*$ is the restriction of $q$ to $N_*$.
\end{lemma}
Proof: The proof of lemma \ref{Abspaltung1Sp} can be adapted with obvious modifications.
\qed
\Absatz
\Numerierung \label{explizitesMatchingorthogonal}{\bf The explicit norm map $\NNN$.}
Let $(\GL_{2n}(R)\times R^\times)_{R\Theta ss}/traf$
  be the set of transformation classes of $R$-$\Theta$-semisimple
  elements $(h,a)\in\GL_{2n}(R)\times R^\times$  under the transformations
  $(h,a)\mapsto (\tran ghg,\det g^{-1}\cdot a)$
  for $g\in GL_{2n}(R), a\in R^\times$. Similarly let $\SO_{2n+1}(R)_{Rss}/conj$ be the set of
  conjugacy classes of $R$-semisimple elements in $\SO_{2n+1}(R)$. We define a norm map
  $$\NNN: (\GL_{2n}(R)\times R^\times)_{R\Theta ss}/traf \longrightarrow \SO_{2n+1}(R)_{Rss}/conj$$
  as follows: If $(h,a)\in \GL_{2n}(R)\times R^\times$  represents a class of the left hand side
  and if $h=p+q$ is the decomposition in the symmetric part $q$ and the skew-symmetric part $p$,
  we decompose
  $M=R^{2n}=M_+\oplus M_-\oplus M_0$ as in lemma \ref{+-0Zerlegung}. The form $q_*=q_+\oplus q_0$ on
  $M_*=M_+\oplus M_0$ is unimodular. Since the ranks of $M$ and $M_-$ are even we have
   $M_*\simeq R^{2r}$ for some $r\in\NN_0$. Let $p'_*$ and $q'_*$ be the $2r\times 2r$-matrices which
  describe $p_*$ and $q_*$ with respect to the standard basis ob $R^{2r}$.
  Let $\tilde q_-$ be a symmetric bilinear form
  on $\tilde M_-:=R^{2(n-r)+1}$ such that $\Delta(q'_*)\cdot \Delta(\tilde q_-)\in (R^\times)^2$.
  By lemma \ref{orthogonaleNormalform} we have an isomorphism of quadratic spaces
  \begin{align*}
    i: (M_*,q_*)\oplus (\tilde M_-,\tilde q_-)\qquad\tilde\longrightarrow \qquad (R^{2n+1},J_{2n+1})
  \end{align*}
  (observe $\det(J_{2n+1})=1$) which induces an injection
  \begin{align*}
    j: \ORTH(M_*,q_*)\times \ORTH(\tilde M_-,\tilde q_-) \quad\hookrightarrow\quad
    \ORTH \left( M_*\oplus \tilde M_-,\; q_*\oplus \tilde q_-\right)
     \quad\tilde\rightarrow \quad\ORTH_{2n+1}.
  \end{align*}
  This injection is canonical (i.e. independent of the chosen isomorphism $i$) on the set of
  conjugacy classes.
\klabsatz
  Now $\NNN(h)$, the image of the class of $h$ under $\NNN$, is defined to be the
  $\ORTH_{2n+1}(R)$-conjugacy  class of
  $j(\tilde C(p'_*),1_{2(n-r)+1})\in \ORTH_{2n+1}(R)$, where we use the
  Cayley-trans\-form-map $\tilde C$ with respect to $q'_*$ from lemma \ref{Cayleytransformationorthogonal}.
  We observe that $\det(\tilde C(p'_*))=1$ by lemma \ref{Cayleytransformationorthogonal}(d) and therefore
  $\NNN(h)$ lies in $\SO_{2n+1}(R)$. Since the centralizer of $j(\tilde C(p'_*),1_{2(n-r)+1})$ in $\ORTH_{2n+1}(R)$
  contains $\{1_{2r}\}\times\ORTH_{2(n-r)+1}(R)$ i.e. elements of determinant $-1$,
  the $\ORTH_{2n+1}(R)$-conjugacy class is in fact a $\SO_{2n+1}(R)$-conjugacy class.

\begin{lemma}\label{Spinnormformel}
  In the notations of \ref{explizitesMatchingorthogonal} the spinor norm of $\NNN(h)$ is the class of
  $\det(h)\mod (R^\times)^2$.
\end{lemma}
Proof: It is sufficient to consider the case $R=F$, since we have an injection
  $\OOO_F^\times/(\OOO_F^\times)^2\hookrightarrow F^\times/(F^\times)^2$. If $\sigma$ denotes the
  spinor norm of $\NNN(h)$ we have by a theorem of Zassenhaus (comp. \cite{Zas}) in the version of \cite{Mas}
  \begin{align*}
    \sigma\quad\equiv\quad \det\left(id-\tilde C(p'_*)\right)\cdot
    \Delta( q'_*)\mod (R^\times)^2.
  \end{align*}
  But $id-\tilde C(p'_*)=(q'_*-p'_*)^{-1}\cdot 2\cdot q'_*$ so that we get
  $\sigma\equiv\det(q'_*-p'_*)^{-1}\cdot 2^{2r}\equiv\det(q'_*-p'_*)\mod (R^\times)^2$.
  Furthermore $\det(q'_*-p'_*)=\det\tran(q'_*-p'_*)=
  \det(p'_*+q'_*)$. Since the discriminant of $p_-$ is a square
  we finally get $\sigma\equiv\det(p'_*+q'_*)\cdot \det(p_-)\equiv\det(h)\mod (R^\times)^2$.
\qed

\begin{prop}\label{explizitesMatchingSpin}
$\ $
\begin{itemize}
\item[(a)] The map $\NNN: (\GL_{2n}(R)\times R^\times)_{R\Theta ss}/traf \longrightarrow \SO_{2n+1}(R)_{Rss}/conj$
  is well defined and surjective.
  Two classes lie in the same fiber iff they have representatives of the form $(h,a_1)$ and $(h,a_2)$.
\item[(b)] If $(h,a)$ represents a class in $(\GL_{2n}(R)\times R^\times)_{R\Theta ss}/traf$ then
  $(h\cdot J^{-1},a)\Theta\in (\GL_{2n}(R)\times R^\times)\rtimes \langle\Theta\rangle$
  matches in the sense of $\Theta$-endoscopy with some ele\-ment $\eta\in\GSPIN_{2n+1}(R)$,
  which maps to $\NNN(h)$ under the projection $pr_{ad}:\GSPIN_{2n+1}\to \SO_{2n+1}$.
\end{itemize}
\end{prop}
Proof: (a) If we replace  $p'_*$ by some $\tran g\cdot p'_*\cdot g$ for $g\in \ORTH(q'_*,R)$,  this
  does not change the conjugacy class of $\NNN(h)$ by lemma
  \ref{Cayleytransformationorthogonal}(c). Since the effect of the other choices has already been
  considered, the map $\NNN$ is well defined.
  \klabsatz To prove
  surjectivity first observe that each class $b \in \SO_{2n+1}(R)_{Rss}/conj$ can be represented after some
  transformation of $J_{2n+1}$ in the
  form $(b',1_{2(n-r)+1})$ with $b'\in \SO(q'_*,R)_{ess}$ by lemma \ref{Abspaltung1O} with a unique
  $r\le n$ and some symmetric $q'_*\in \GL_{2r}(R)$.
  One should think of  $(b',1_{2(n-r)+1})$ as a block-matrix
  $$\begin{pmatrix} B_{11} & 0 & B_{12}\\ 0& 1_{2(n-r)+1} &0\\ B_{21}&0& B_{22}
  \end{pmatrix}
  \quad \text{ with }\quad b'=
  \begin{pmatrix}B_{11}&B_{22}\\ B_{21}& B_{22}
  \end{pmatrix}$$
  Since the class of $\Delta(q'_*)$ in
  $R^\times/(R^\times)^2$ is the inverse of the class of $\Delta(J_{2n+1}|\ker(b-1))$,
  the transformation class of $q'_*$ is unique by lemma \ref{orthogonaleNormalform}.
  Up to this the $\SO(q'_*,R)$-conjugacy-class of $b'$ is unique. The bijectivity of the Cayley-transform map
  and property \ref{Cayleytransformationorthogonal}(c) then imply that there is a $p'\in Alt_{2r}(R)$,
  which is unique up to transformations with elements of $\SO(q'_*,R)$, such that $b=\tilde C(p')$.
  Now we consider the unimodular bilinear form $h'=p'+q'$ on $R^{2r}$, which is unique up to transformations
  with elements of $\GL_{2r}(R)$, and some unimodular skew symmetric
   form $p_-$ on $R^{2(n-r)}$. The form $p_-\oplus h'$ on $R^{2n}$ is then unimodular
  and $R$-$\Theta$-semisimple i.e. corresponds to a $R$-$\Theta$-semisimple transformation class $h$.
  For every $R$-semisimple $a\in R^\times$ we get $\NNN(h,a)=b$.  Since the transformation class of $h'$
  is unique by the above considerations and by lemma \ref{symplektischeNormalform}
  we conclude that the fibers of $\NNN$ correspond to the different choices for
  the $R$-semisimple element $a\in R^\times$.
  \absatz
(b) At first we consider the case that $R=\bar F$ is an algebraically closed field, so that we may assume that
  $\gamma=h\cdot J_{2n+1}^{-1}$
  has diagonal form $\gamma=diag(t_1,\ldots,t_{2n})$. After applying a permutation in
  $W_{\Sp_{2n}}$ we may assume
  \begin{gather}\label{Annahme an h zwei}
    t_i\ne t_{2n+1-i} \text{ for }i\le r \text{ and }t_i=t_{2n+1-i} \text{ for }
    r+1\le i\le 2n-r.
  \end{gather}
  We have:
  \begin{align*} h\quad&=\quad antidiag(t_1,-t_2,t_3,\ldots,-t_{2n})\\
              h\pm\tran h\quad&=\quad antidiag(t_1\mp t_{2n},-(t_2\mp t_{2n-1}),t_3\mp t_{2n-2},
              \ldots,-(t_{2n}\mp t_1))\\
              -\tran h^{-1}\cdot h\quad&=\quad diag(t_{2n}/t_1,t_{2n-1}/t_2,\ldots,t_1/t_{2n})
  \end{align*}
  Thus $M_-\simeq R^{2(n-r)}$ is generated by the standard basis elements
  $e_{r+1},\ldots, e_{2n-r}$ of $R^{2n}$ and $M_*=M_-\oplus M_0$ has the basis $e_1,e_2,,\ldots,e_r,
  e_{2n+1-r},\ldots,e_{2n}$. The matrix of the symmetric bilinear form given by
  $\tran d\cdot (h+\tran h)\cdot d$
  with respect to this latter basis has the standard form $J'_{2r}$, if we take
  $d=diag((t_1-t_{2n})^{-1},(t_2-t_{2n-1})^{-1},\ldots,(t_r-t_{2n+1-r})^{-1},1,\ldots,1)$.
  Since $d$ and $\tran h^{-1}\cdot h$ commute we get from the definition of $\NNN(h)$ that
  $\NNN(h)$ is represented by the diagonal matrix
  \begin{gather*}
    diag(t_{2n}/t_1,t_{2n}/t_2,\ldots,t_{n+1}/t_n,1,t_n/t_{n+1},\ldots,t_1/t_{2n}),
  \end{gather*}
  which can be conjugated into the form of example \ref{BspMatchingGL2nGSpin2n+1}. This proves
  the claim in the case that $R$ is an algebraically closed field.
\Absatz
  In the case that  $R$ is arbitrary we consider the commutative diagram with
  exact rows and columns and a connecting homomorphism marked with $\ldots$ (snake lemma):
  \begin{align*}\begin{CD}
   &&&&&&   1&&\\
   &&&&&&@VVV&&\\
   &&&&1 @>>>  \{\pm 1\}@>>> \ldots\\
   &&&& @VVV  @VVV &&\\
   && 1 @>>> \Spin_{2n+1}(R) @= \Spin_{2n+1}(R) @>>> 1\\
   && @VVV  @VVV  @VVV &&\\
   1 @>>> R^\times @>>> \GSPIN_{2n+1}(R) @>pr_{ad}>> \SO_{2n+1}(R) @>>> 1\\
   && @VVV @VV\mu V @VV\text{Spinnorm}V &&\\
   \ldots @>>> R^\times @>{r\mapsto r^2}>> R^\times @>>> R^\times/(R^\times)^2 @>>> 1\\
   && @VVV @VVV @VVV &&\\
   && 1 && 1 && 1 &&
  \end{CD}\end{align*}
  It follows from this diagram and lemma \ref{Spinnormformel} that a matrix $\eta_0$ in the class $\NNN(h)$
  has a preimage $\eta\in \GSPIN_{2n+1}(R)$ such that $\mu(\eta)= \det(h)\cdot a^{2}$ and that the set
  $\{x\in \GSPIN_{2n+1}(\bar F)| pr_{ad}(x)=\eta_0,\; \mu(x)= \det(h)\cdot a^{2}\}$ just consists of $\pm \eta$.
  On the other hand by example \ref{BspMatchingmurelation}
   an element $\eta'\in \GSPIN_{2n+1}(\bar F)$
  matching with $(h,a)$ satisfies $\mu(\eta')= \det(h)\cdot a^{2}$. From the validity of the proposition
  over $\bar F$ now follows that either $\eta$ or $-\eta$ matches with $(h,a)$. This element has all desired
  properties.
\qed\absatz
  {\it Remark}: To get (b) it is even in the case $R=F$ not sufficient just to apply Steinbergs
  theorem on rational elements to the rational conjugacy class inside $\GSPIN_{2n+1}(\bar F)$,
   which matches with $(h\cdot J^{-1},a)\Theta$:
  If $\eta'\in \GSPIN_{2n+1}(F)$ denotes such an element, then we only know from the case
  $R=\bar F$ that $pr_{ad}(\eta')$ and $\NNN(h)$ are stably conjugate elements in $\SO_{2n+1}(F)$.
  But the Spinor norm is not invariant under stable conjugation. Thus it is not clear without the
  use of lemma \ref{Spinnormformel} that $\NNN(h)$ can be lifted to a class in $\GSPIN_{2n+1}(F)$, on which
   the multiplier $\mu$ takes the correct value.

\begin{cor}\label{SurjMatchingGLSPIN}
  For each semisimple $\eta\in \GSPIN_{2n+1}(F)$ there exists an $F$-$\Theta$-semisimple
  $(h\cdot J^{-1},a)\Theta\in \left(\GL_{2n}(F)\times F^\times\right)\rtimes \langle\Theta\rangle$ matching with $\eta$.
\end{cor}
Proof:  By \ref{explizitesMatchingSpin}(a) for $R=F$ there exists $(h,a_1)\in \GL_{2n}(F)\times F^\times$
  with $pr_{ad}(\eta)\in \NNN(h)$ and by (b) there exists $\eta_1\in \GSPIN_{2n+1}(F)$ matching with
  $(h\cdot J^{-1},a_1)\Theta$ such that $pr_{ad}(\eta_1)=pr_{ad}(\eta)$. It follows $\eta=\eta_1\cdot b$ for some
  $b\in F^\times\simeq Center(\GSPIN_{2n+1}(F))$. Then 
  $(h\cdot J^{-1},a_1\cdot b)\Theta$ 
   matches with $\eta$.
\qed
\begin{lemma}\label{GlGmkonjbereitsueberR}
  For $G=\GL_{2n}\times \GG_m$ let $\gamma_1,\gamma_2\in G(\OOO_F), g_F\in G(\bar F) $ be such that
  $\gamma_2\Theta=g_F\cdot \gamma_1\Theta\cdot g_F^{-1}$ with $\Theta$ as in
  example \ref{BspGL2nGl1GSpin2n+1}. Then there exists $g_R\in G(\OOO_F)$ with
  $\gamma_2\Theta=g_R\cdot \gamma_1\Theta\cdot g_R^{-1}$.
\end{lemma}
Proof: Write $\gamma_i=(h_i\cdot J_{2n}^{-1},a_i),\; g_F=(h_F,\tilde a)$. Then the assumption means:
  $h_2=h_F\cdot h_1\cdot \tran h_F$ and $a_2=\tilde a\cdot a_1 \cdot \det(h_F)^{-1}\cdot \tilde a^{-1}$
  which implies $\det(h_F)\in\OOO_F^\times$. By lemma \ref{GlkonjbereitsueberR}(a) there exists
  $h_R\in \GL_{2n}(\OOO_F)$ with $h_2=h_R\cdot h_1\cdot \tran h_R$. This implies $\det(h_R)^2=\det(h_F)^2$.
  If $\det(h_R)=-\det(h_F)$ then $h_F^{-1}\cdot h_R\in\ORTH(h_1)(\bar F)$ where
  $\ORTH(h_1)=\{h\in\GL_{2n}\mid \tran h\cdot h_1\cdot h=h_1\}=\ORTH(q_{+,1})\times
  \left(\SP(p_{-,1}\oplus p_{0,1})\cap  \ORTH(q_{-,1}\oplus q_{0,1})\right)$ has determinant $-1$.
   This implies $M_{+,1}\ne 0$ so that we get an
  element $h_{\epsilon}\in \ORTH(h_1)(\OOO_F)$ of determinant $-1$. Replacing $h_R$ by $h_R\cdot h_{\epsilon}$
   we can now assume $\det(h_R)=\det(h_F)$.
  With $g_R=(h_R,1)$ we now have $\gamma_2\Theta=g_R\cdot \gamma_1\Theta\cdot g_R^{-1}$.
  \qed
\begin{lemma}\label{MatchinglemmaGspin}
  Let $(h,a)=(sJ,a)\in \GL_{2n}(\OOO_F)\times \OOO_F^\times$ be $\OOO_F$-$\Theta$-semisimple and
  $b=(1_{2(n-r)+1},b_*)\in \SO_{2n+1}(\OOO_F)$
  a representing element of $\NNN(h)$. With $p_*,q_*,p_-$ as in \ref{explizitesMatchingorthogonal}
  assume that we have matching topologically unipotent elements
  $u_-\in \Sp_{2(n-r)}(F)\simeq \Sp(p_- )$ and
  $v_-\in \SO_{2(n-r)+1}(F)\simeq (\ker(b-1)(F)\cap \SO_{2n+1}(F))$
  and an additional topologically unipotent element
  $u_*\in \SO(q_*)(F)\cap \Sp(p_*)(F)\simeq Cent(b_*,\SO(q_*)(F))$.
  \neueZeile
  Then the element
  $\gamma\Theta= (s,a)\Theta\cdot (u_-,u_*)=(u_-,u_*)\cdot (s,a)\Theta\in
  \left(GL_{2n}(F)\times F^\times\right)\rtimes \langle\Theta\rangle$
  matches with some element $\eta\in \GSPIN_{2n+1}(F)$, which projects to
  $\beta:=(v_-^2,u_*^2)\cdot b=b\cdot (v_-^2,u_*^2)\in \SO_{2n+1}(F)$.
\end{lemma}
Proof: We first prove the existence of $\eta\in \GSPIN_{2n+1}(\bar F)$ with the desired
  properties and thus work over the algebraically closed
  field $F=\bar F$ as in the proof of lemma \ref{NNN prop PGLSP}(c). Thus we may assume that
  $\gamma$ and $\beta$ lie in the corresponding diagonal tori. The same
  holds for the residually semisimple parts $s$ resp. $b$ and the topologically unipotent parts
  $u=(u_+,u_*)$ and $v=(v_-^2,u_*^2)$. As $(s,a)\Theta$ matches with some $\eta_s$ which projects to $b$
  by  \ref{explizitesMatchingSpin}(b) we only have to examine the topologically
  unipotent elements. We can make the assumption (\ref{Annahme an h zwei}) and write
  \begin{align*}
    u_-\quad &=\quad diag(w_{r+1},\ldots,w_n,w_n^{-1},\ldots,w_{r+1}^{-1})\quad\in\quad \Sp_{2(n-r)}(\bar F)\\
    u_*\quad &=\quad diag(w_1,\ldots,w_r,w_r^{-1},\ldots,w_1^{-1})\quad\in\quad \SO_{2r}(\bar F)^{b_*}.
    \intertext{By the definition of $BC$-matching we can assume}
    v_- \quad &=\quad  diag(w_{r+1},\ldots,w_n,1,w_n^{-1},\ldots,w_{r+1}^{-1})\quad \in\quad \SO_{2(n-r)+1}(\bar F).
  \end{align*}
  Taking everything together we get by the description of $M_-$ and $M_*$ in the proof of
  lemma \ref{explizitesMatchingSpin}(b):
  \begin{align*}
    u\quad &=\quad(w_1,\ldots,w_n,w_n^{-1},\ldots,w_1^{-1})\\
    v\quad &=\quad(w_1^2,\ldots,w_n^2,1,w_n^{-2},\ldots,w_1^{-2})
  \end{align*}
  and the matching between $\gamma\Theta$ and some $\eta\in \GSPIN_{2n+1}(\bar F)$ which projects to $\beta$
  follows from example \ref{BspMatchingGL2nGSpin2n+1}.
\absatz
  To get $\eta$ as an element of $\GSPIN_{2n+1}(F)$ we observe that the determinant of $\gamma J_{2n}^{-1}$
  equals the spinor norm of $\beta$ as an element of $F^\times/(F^\times)^2$: This is already clear by
  \ref{Spinnormformel} for the residually semisimple parts, but both topologically unipotent parts lead to
  the neutral element in $F^\times/(F^\times)^2$, since $2\ne p$ by assumption. Now one argues as in the
  proof of \ref{explizitesMatchingSpin}(b) to get $\eta$ as an $F$-rational element.
\qed

\begin{thm} \label{BC2Theorem}
 ($\text{BC}_{\text 2}$) is true.
\end{thm}
 Proof:
 We observe that every pair of  $BC$-matching (topologically unipotent) elements
 $\bar\gamma\in\SO_5(F)$ and $\eta_1\in \Sp_4(F)$ can be obtained from
 a pair of (topologically unipotent) elements $\gamma\in \GSp_4(F)\simeq \GSPIN_5(F)$
 and $\eta\Theta=\Theta\eta\in (\widetilde{\GL_4}\times \GL_1)(F)$ such that $\bar\gamma=pr_{ad}(\gamma)$ and
 $\eta=(\eta_1,a)\in (\GL_4\times\GL_1)^\Theta(F)\simeq (\Sp_4\times \GL_1)(F)$
 and such that $\gamma^2$ matches with $\eta\Theta$ in the sense of
 \ref{DefMatching}. This follows immediately from the definition of $BC$-matching
 \ref{BCMatching} and example \ref{BspMatchingGL2nGSpin2n+1}.
\absatz
 If we apply lemma \ref{Redorbintaufadjungiert} in the case $G=\GSPIN_5\simeq \GSp_4,\;
 T=\G_m, H=\SO_5$ and lemma \ref{Homogenitaet} we get
 $$\STO_{\gamma^2}(1,\GSp_4)=\STO_{\bar \gamma}(1,\SO_5)$$
 Since we have $\STO_\eta(1,\Sp_4\times\GL_1)=\STO_{\eta_1}(1,\Sp_4)$ by lemma \ref{Redorbintaufadjungiert}
 the statement of ($\text{BC}_{\text 2}$) is equivalent to the identity
 $$\STO_\gamma(1,\GSp_4)=\STO_\eta(1,\Sp_4\times\GL_1)$$ for matching topologically unipotent
 $\gamma\in \GSp_4(F)$
 and $\eta\Theta=\Theta\eta\in (\widetilde{\GL_4}\times \GL_1)(F)$.
 In the case that $\eta$ is strongly $\Theta$-regular
 this has been proved in \cite[ch. II]{FGL4}. The general case follows by the germ
 expansion principle as in \cite{HalesSp4SO5}, \cite{Rogawski}.
\qed
\begin{cor}[Fundamental lemma for $\Sp_4\leftrightarrow\widetilde{\PGL}_5$]\label{BC2TheoremCor}
  If $\gamma\Theta\in \widetilde{\PGL}_5(F)$ and $h\in \Sp_{4}(F)$ are matching semisimple elements  then we have
  \begin{gather*}
    \STO_{\gamma\Theta}(1,\widetilde{\PGL}_5)=\STO_h(1, \Sp_{4}).
  \end{gather*}
\end{cor}
Proof: This follows from theorem \ref{BC2Theorem}, ($\text{BC}_{\text 1}$) (compare \ref{BCVermutung}) and theorem
  \ref{FulemmaSp2nPGL2n+1}.
\qed
\begin{thm}\label{FulemmaGSpin2n+1GL2n} Let $G=\GL_{2n}\times\GG_m$.
If $\gamma\Theta\in \tilde G(F)$ and $\eta\in \GSPIN_{2n+1}(F)$ are matching semisimple elements and if
conjecture ($\text{BC}_{\text{m}}$) is true for all $m\le n$  then we have
\begin{gather}\label{Fule2n+12nGleichung}
\STO_{\gamma\Theta}(1,\tilde G)=\STO_\eta(1,\GSPIN_{2n+1}).
\end{gather}
\end{thm}
Proof: Let $h=pr_{ad}(\eta)\in \SO_{2n+1}(F)$. In view of lemma \ref{Redorbintaufadjungiert} we have to prove
  \begin{align} \label{FuleO2n+12nGleichung}
    \STO_{\gamma\Theta}(1,\tilde G)=\STO_h(1,\SO_{2n+1}).
  \end{align}
  The proof is now similar to the proof of Theorem \ref{FulemmaSp2nPGL2n+1}.
\absatz
\underline{Step 1}: Let us assume that the right hand side of (\ref{FuleO2n+12nGleichung})
does not vanish. Then there exists $h'\in \SO_{2n+1}(F)$ stably conjugate to $h$ which has a
nonvanishing orbital integral, i.e. can be conjugated
into $\SO_{2n+1}(\OOO_F)$. We can assume that $h'\in \SO_{2n+1}(\OOO_F)$. Since $h'=pr_{ad}(\eta')$ for some
$\eta'\in\GSPIN_{2n+1}(F)$ in the stable conjugacy class of $\eta$ we can assume without loss
of generality that $\eta'=\eta$ and thus $h'=h$. Furthermore we can assume that the topological
Jordan decomposition  is of the form $h=b\cdot v=v\cdot b$  with residually semisimple
$b=(1_{2(n-r)+1},b_*)\in \SO_{2(n-r)+1}(\OOO_F)\times \SO(q_*,\OOO_F)_{ess}$ and topologically unipotent
$v$, where $q_*$ denotes the restriction of $J_{2n+1}$ to the orthogonal complement of
$\ker(b-1)\simeq \OOO_F^{2(n-r)+1}$. Here we observe that the restriction of $J_{2n+1}$ to $\ker(b-1)$
can be assumed to be a multiple of the standard form $J_{2(n-r)+1}$ and that $b_*$ has determinant $1$ as
$1_{2(n-r)+1}$ has determinant $1$.
We can write $v$ in the form $v =((v_-)^2,(u_*)^2)$ with
$v_-\in \SO_{2(n-r)+1}(\OOO_F)$ and $u_*\in Cent\left(b_*,\SO(q_*,\OOO_F)\right)$ since $2\in \Z_p^\times$.
\absatz
We remark that condition ($*$) in the Kazhdan-lemma \ref{Kazhdan-Lemma} is satisfied by lemma
\ref{GlGmkonjbereitsueberR}, so that we get  nonvanishing
 \begin{align*}
   O_{h}(1,\SO_{2n+1})\quad &=\quad O_{v}(1,Cent(b,\SO_{2n+1}))
 \end{align*}
Observing that we have an isomorphism
 \begin{align*}
   Cent(b,\SO_{2n+1})\quad\simeq\quad \SO_{2(n-r)+1}\times Cent(b_*,\SO(q_*))\times \{\pm 1\}
 \end{align*}
we can decompose the orbital integral on the right hand side using lemma \ref{EndlichesZentrum}:
 \begin{align*}
   O_{h}(1,\SO_{2n+1})\quad &=\quad
   O_{(v_-)^2}(1,\SO_{2(n-r)+1})\cdot O_{(u_*)^2}(1,Cent(b_*,\SO(q_*)))\\
   \notag &=\quad O_{v_-}(1,\SO_{2(n-r)+1})\cdot O_{u_*}(1,Cent(b_*,\SO(q_*))),
 \end{align*}
 (use lemma \ref{Homogenitaet} in the last step)
 i.e. $\STO_{v_-}(1,\SO_{2(n-r)+1})$ (being the sum of integrals of nonnegative functions) is strictly
 positive. Since $v_-$ is $BC$-matching with some $u'_-\in \Sp_{2(n-r)}(F)$
 the equation ($\text{BC}_{n-r}$)
 implies that there exists $u_-\in \SP_{2(n-r)}(F)$ matching with  $v_-$ and with strictly
 positive orbital integral, i.e. we can assume $u_-\in \SP_{2(n-r)}(\OOO_F)$. Let $s\in
 \GL_{2n}(\OOO_F)$ be a residually semisimple element with $\NNN(s\cdot J_{2n}^{-1},a)=b$ for some
  $a\in \OOO_F^\times$, chosen such that we can identify
 the corresponding $q_*$ on $M_*$ with the above obtained $q_*$. By modifying $a$ we can assume that
 $(s,a)\Theta$ matches with the residually semisimple part $\eta_s$ of $\eta$.  We define the element
 $\gamma'\Theta=(u_-,u_*)\cdot s\Theta= s\Theta\cdot (u_-,u_*)\in\tilde{G}(\OOO_F)$.
 The element $\gamma'\Theta\in \tilde G(\OOO_F)$ matches with $h$ by lemma \ref{MatchinglemmaGspin}
 and therefore lies in the stable conjugacy class of $\gamma\Theta$.
 \absatz
  If the left hand side of
 (\ref{FuleO2n+12nGleichung}) does not vanish, it is immediate that there exists $\gamma'\Theta\in
 \tilde G(\OOO_F)$ in the stable conjugacy class of $\gamma_\Theta$. By reversing the
 above arguments we see that there exists $h'\in \SO_{2n+1}(\OOO_F)$ in the stable class of $h$.
 So excluding the tautological case that (\ref{Fule2n+12nGleichung}) means $0=0$ we may assume without
 loss of generality that $\gamma\in G(\OOO_F)$ and $h\in \SO_{2n+1}(\OOO_F)$.
 We may furthermore assume that $\gamma\Theta=(u_+,u_*)\cdot s\Theta$ and $h=(v_+^2,u_*^2)\cdot b$
 are the topological Jordan decompositions with $BC$-matching $u_+$ and $v_+$ and matching
 residually semisimple $s\Theta$ and $\eta_s$.
 \absatz

\underline{Step 2}:  As in the proof of theorem \ref{FulemmaSp2nPGL2n+1} we get from
 lemma \ref{SpkonjbereitsueberR} the fact that we
 obtain all relevant conjugacy classes  in the stable conjugacy class of $h$ if we let
  $v'_-$ vary through a set of representatives of the stable conjugacy class of $v_-$ in
  $\SO_{2(n-r)+1}(F)$ and $u'_*$ through a set of representatives of the stable conjugacy class of $u_*$
  in $Cent(b_*,\SO(q_*))$ and then consider all $h'= b\cdot ((v'_+)^2,(u'_*)^2)$.  i.e.
\begin{align}\label{SOorbital}
  \STO_{h}(1,\SO_{2n+1})\;
          &=\; \sum_{v'_+\sim v_+} O_{(v'_+)^2}(1,\SO_{2(n-r)+1})
             \cdot\!\!\! \sum_{u'_*\sim u_*} O_{(u'_*)^2}(1,Cent(b_*,\SO(q_*)))\\
  \notag  &=\quad \sum_{v'_+\sim v_+} O_{v'_+}(1,\SO_{2(n-r)+1})
             \cdot \sum_{u'_*\sim u_*} O_{u'_*}(1,Cent(b_*,\SO(q_*))).
\end{align}
\underline{Step 3}: We can repeat this argument in the $\Theta$-twisted situation,
  since by lemma \ref{GlGmkonjbereitsueberR}
  the class of the residually semisimple part $(s,a)\Theta$ of $\gamma\Theta$   is the only
   $G(F)$-conjugacy class inside the stable class of $(s,a)\Theta$,
  which meets $G(\OOO_F)$ and since the Kazhdan-Lemma \ref{Kazhdan-Lemma} holds for $\tilde G$ by the same
  lemma. We remark that $G^{(s,a)\Theta} \simeq \SP_{2(n-r)}\times \Cent(b_*,SO(q_*))\times \GG_m$
  by the definition of $\Theta$ and lemma \ref{Zentralisatoren}(e), so $G^{(s,a)\Theta}$ is connected
   and we can use lemma
  \ref{Redorbintaufadjungiert}
   to get rid of the $\GG_m$ factors in the following orbital integrals.
  If we denote by $u'_-$ a set of representatives for the
  $\SP_{2(n-r)}(F)$-conjugacy classes in the stable class of $u_-\in \SP_{2(n-r)}(\OOO_F)$ we get
 \begin{align}\label{GlGmorbital}
 \STO_{\bar{\gamma}\Theta}(1, \tilde G)\quad &= \quad
  \sum_{(u'_-,u'_*)\sim (u_-,u_*)} O_{(u'_-,u'_*)}(1,G^{s\Theta})\\
  \notag&=\quad \sum_{u'_-+\sim u_-} O_{u'_-}(1,\Sp_{2(n-r)})\cdot
  \sum_{u'_*\sim u_*} O_{u'_*}(1,Cent(b_*,\SO(q_*))).
 \end{align}
\underline{Step 4}: Since $v_-$ and $u_-$ are $BC$-matching we can apply ($\text{BC}_{n-r}$) to get
that the right hand sides of  (\ref{SOorbital}) and (\ref{GlGmorbital}) coincide, which proves the theorem.

\qed

 \Seitenumbruch
 \meinchapter{Comparison between $\SO_{2n+2}$ and $\Sp_{2n}$}
  Let $R$ be as in \ref{R theta semisimple}.
\begin{lemma}\label{Abspaltungplusminus1}
  Let $N$ be a free $R$-module.
 \begin{itemize}
 \item[(a)]
  If $p$ is a unimodular symplectic form on $N$ and if $\beta\in \Sp(p,R)$ is $R$-semisimple
  then there exists a $\beta$-invariant orthogonal (with respect to $p$) direct sum decomposition
  $N=N_+\oplus N_-\oplus N_*$ such that $\beta$ acts as identity on $N_+$, as $-id$ on $N_-$ and
  $\beta_*=\beta|N_*\in \Sp(p_*) $ satisfies $\beta_*-\beta_*^{-1}\in\GL(N_*)$,
  where $p_*$ is the restriction of $p$ to $N_*$.
 \item[(b)]
  If $q$ is a unimodular symmetric bilinear form on $N$ and $b\in \ORTH(q,R)$ is $R$-semisimple
  then there exists a $b$-invariant orthogonal (with respect to $q$) direct sum decomposition
  $N=N_+\oplus N_-\oplus N_*$ such that $b$ acts as identity on $N_+$, as $-id$ on $N_-$ and
  $b_*=b|N_*\in \ORTH(q_*) $ satisfies $b_*-b_*^{-1}\in\GL(N_*)$, where $q_*$ is the restriction of $q$ to $N_*$.
 \end{itemize}
\end{lemma}
Proof: The proof of lemma \ref{Abspaltung1Sp} can be adapted with obvious modifications: We have
  $b-b^{-1}=b^{-1}\cdot (b-1)\cdot (b+1)$, so that $b-b^{-1}\in\GL(N_*)$ is equivalent to $b-1,b+1\in
  \GL(N_*)$.
\qed
\begin{lemma}\label{SOSPschieben}
 Let  $b\in \GL_n(R)$ satisfy $b-b^{-1}\in \GL_n(R)$. Then the following holds:
 \begin{itemize}
  \item[(a)] If $q\in \GL_n(R)$ is symmetric and $b\in \ORTH(q,R)$ then the matrix $p=q\cdot(b-b^{-1})$
    is unimodular skew-symmetric and we have $b\in \SP(p,R)$.
  \item[(b)]  If $p\in \GL_n(R)$ is skew-symmetric and $b\in \SP(p,R)$ then the matrix
    $q=p\cdot(b-b^{-1})^{-1}$ is unimodular symmetric and we have $b\in \SO(q,R)$.
  \item[(c)]  Under the conditions of (a) and (b) we have:
    \begin{gather*}
      Cent(b,\ORTH(q))=Cent(b,\SP(p))=Cent(b,\SO(q)).
     \end{gather*}
  \item[(d)] The above statements and formulas are invariant under the substitutions
  $b\mapsto g^{-1}bg, q\mapsto \tran gqg, p\mapsto \tran gpg$ for $g\in \GL_n(R)$.
 \end{itemize}
\end{lemma}
Proof: (a) We have $p=q\cdot b-(\tran b\cdot q\cdot b)\cdot b^{-1}=qb-\tran(qb)$ and
   $\tran bpb=\tran b(qb)b-\tran bqb^{-1}\cdot b=qb-\tran bqb\cdot b^{-1}=qb-qb^{-1}=p$.
   Unimodularity of $p$ follows from $q,(b-b^{-1})\in \GL_n(R)$.
\klabsatz (b) We have
   $\tran q=q\Leftrightarrow -\tran(b-b^{-1})^{-1}\cdot p=p\cdot(b-b^{-1})^{-1}\Leftrightarrow
   p(b-b^{-1})=(-\tran b+\tran b^{-1})p\Leftrightarrow
   \tran b^{-1}\cdot(\tran bpb -p)=(p-\tran bpb)\cdot b^{-1} \Leftarrow b\in \SP(p,R)$ and
   $\tran bqb=\tran bp(b-b^{-1})^{-1}\cdot b=\tran bpb\cdot (b-b^{-1})^{-1}=p(b-b^{-1})^{-1}=q$.
   As an element of a symplectic group $b$ has determinant $1$.
 \klabsatz (c) For $x\in Cent(b,\GL_n(R))$ we have $xb=bx$ and $b^{-1}x=xb^{-1}$ which imply
   $x(b-b^{-1})=(b-b^{-1})x$, so that we get
   $\tran x qx=q \Leftrightarrow \tran x qx(b-b^{-1})=q(b-b^{-1})\Leftrightarrow
   \tran xq(b-b^{-1})x=q(b-b^{-1})\Leftrightarrow \tran xpx=p$. This proves the first "=".
   The second follows immediately since elements of $\SP(p)$ have determinant $1$.
 \klabsatz (d) follows by almost trivial computations.
\qed
  \Absatz
\Numerierung \label{explizitesMatchingSOSP}
  If $s\in \ORTH_{2n+2}$ with $\det(s)=-1$ denotes a reflection, we can identify the semidirect product
   $\SO_{2n+2}\rtimes\langle\Theta\rangle$ where $\Theta=int(s)$  with the orthogonal group $\ORTH_{2n+2}$.
 \klabsatz
  Let $\ORTH_{2n+2}(R)^-_{Rss}/conj$
  be the set of $\SO_{2n+2}(R)$-conjugacy classes of $R$-semisimple (=$R$-$\Theta$-semisimple)
  elements of $h\in\ORTH_{2n+2}(R)$ with $\det(h)=-1$. Recall that $Sp_{2n}(R)_{Rss}/conj$ is the set of
  conjugacy classes of $R$-semisimple elements in $Sp_{2n}(R)$. We define a norm map
  $$\NNN: \ORTH_{2n+2}(R)^-_{Rss}/conj \longrightarrow \Sp_{2n}(R)_{Rss}/conj$$
  as follows: If $b\in \ORTH_{2n+2}(R)$ represents a class of the left hand side, we decompose
  $N=R^{2n+2}=N_+\oplus N_-\oplus N_*$ as in lemma \ref{Abspaltungplusminus1}(b). Let
  $b_+=id_{N_+},b_-=-id_{N_-}$ and $b_*=b|N_*$. Let $q_*$ be the restriction of the form
  $J_{2n+2}$ to $N_*$. We may think of $q_*$ as a symmetric matrix after introducing a basis of $N_*$.
   Since $b_*\in \SP(p_*)$ for $p_*=q_*\cdot(b_*-b_*^{-1})$ by
  lemma \ref{SOSPschieben}(a) we have $\det(b_*)=1$. Therefore $-1=\det(b)=\det(b_+)\cdot \det(b_-)
  \cdot\det(b_*)=1\cdot(-1)^{rank N_-}\cdot 1$, i.e. $rank(N_-)$ is odd $=1+2r_-$. Since $rank(N_*)$
  is even by lemma \ref{SOSPschieben} we have $rank(N_+)=1+2r_+$ for some $r_+\in\NN_0$.
  Now we equip the $R$-module $M=M_+\oplus M_-\oplus N_*$ where $M_+\simeq R^{2r_+}, M_-\simeq R^{2r_-}$
  with the alternating form $p=J_{2r_+}\oplus J_{2r_-}\oplus p_*$ and the linear automorphism
  $\beta=id_{M_+}\times -id_{M_-}\times b_*\in \SP(p)$. Identifying the symplectic space
  $(M,p)\simeq (R^{2n},J_{2n})$ we can think of $\beta$ as an element of $\SP_{2n}(R)$.
  The conjugacy class of $\beta$ in $\SP_{2n}(R)$ does not depend on the choices we made (apply
  lemma \ref{SOSPschieben}(d) ) and is the desired $\NNN(b)$. It is clear that $\NNN(b)$ is
  $R$-semisimple.

\begin{prop}\label{propoSOSP}
Let $R$ be as in \ref{R theta semisimple}.
\begin{itemize}
\item[(a)] The map $\NNN: \ORTH_{2n+2}(R)^-_{Rss}/conj \longrightarrow \Sp_{2n}(R)_{Rss}/conj$
  is well defined. Each $b\in\ORTH_{2n+2}(R)^-_{Rss}/conj\;$ matches with $\NNN(b)$ in the sense of $\Theta$-endoscopy
  (compare examples \ref{BspSO2n+2Sp2n}, \ref{BspMatchingSO2n+2Sp2n}).
\item[(b)] The map $\NNN$ is surjective, if $R=\OOO_F$.
  Its fibers are of order $2=\#(R^\times/(R^\times)^2)$ and describe the two different pairs
  $(q_+,q_-)$ of classes of unimodular quadratic forms on $(M_+,M_-)$ such that $\Delta(q_+)\cdot
  \Delta(q_-)\equiv \det(q_*)^{-1} \mod (R^\times)^2$.
\end{itemize}
\end{prop}
Proof: (a)That $\NNN$ is well defined is already clear.
  By the definition of matching we can work over $R=\bar F$, so that we may assume that
  $\gamma=b\cdot s^{-1}\in \SO_{2n+2}(R)$
   has diagonal form $\gamma=diag(t_1,\ldots,t_{n+1},t_{n+1}^{-1},\ldots,t_1^{-1})$, where $s$
   is the reflection defined in \ref{BspSO2n+2Sp2n}. We have:
\begin{align*}b\quad&=\quad
 \begin{pmatrix}diag(t_1,\ldots,t_n)&&&\\
                &0&t_{n+1}&\\ &t_{n+1}^{-1}&0&\\
                &&& diag(t_n^{-1},\ldots,t_1^{-1})
 \end{pmatrix}\\
 &\hbox{\vbox{\vskip10pt}}\\
     b-b^{-1}\quad&=\quad diag(t_1-t_1^{-1},\ldots,t_n-t_n^{-1},0,0,t_n^{-1}-t_n,\ldots,t_1^{-1}-t_1)
\end{align*}
 With the standard basis $(e_i)_{1\le i\le 2n+2}$ of $R^{2n+2}$ we get:
 \begin{align*}
   M_+\quad &=\quad\langle e_i,e_{2n+3-i}\mid t_i=1,\; 1\le i\le n\rangle
                    \oplus \langle t_{n+1}\cdot e_{n+1}+ e_{n+2}\rangle\\
   M_-\quad &=\quad\langle e_i,e_{2n+3-i}\mid t_i=-1, 1\le i\le n\rangle
                    \oplus \langle t_{n+1}\cdot e_{n+1}- e_{n+2}\rangle\\
   M_*\quad &=\quad\langle e_i,e_{2n+3-i}\mid t_i\ne \pm 1, 1\le i\le n\rangle
 \end{align*}
 The corresponding description of $N=N_+\oplus N_-\oplus N_*$ can be arranged such that:
 \begin{align*}
   N_+\quad &=\quad\langle e'_i,e'_{2n+3-i}\mid t_i=1,\; 1\le i\le n\rangle\\
   N_-\quad &=\quad\langle e'_i,e'_{2n+3-i}\mid t_i=-1, 1\le i\le n\rangle\\
   M_*\quad &=\quad\langle e'_i, e'_{2n+3-i}\mid t_i\ne \pm 1, 1\le i\le n\rangle
 \end{align*}
 where $e'_i=(-1)^{i}\cdot(t_i-t_i^{-1})^{-1}e_i$ if $t_i\ne \pm 1$ and $1\le i\le n$ and  $e'_j=e_j$ else.
 With respect to this new basis
 of $M_*$ the symplectic form given by $p_*=q_*\cdot(b_*-b_*^{-1})$ has standard form $J_{2g}$, so that
 the symplectic form $p$ on $R^{2n}$ can be assumed to be of standard form $J_{2n}$ with
 respect to the basis $e'_1,\ldots,e'_n,e'_{n+3},\ldots,e'_{2n+2}$. The symplectic transformation
 $\beta=id_{N_+}\times (-id_{N_-})\times b_*$ in $\NNN(b)$ has the diagonal form
 $diag(t_1,\ldots,t_n,t_n^{-1},\ldots,t_1^{-1})$ with respect to this basis. The claim now follows
 from example \ref{BspMatchingSO2n+2Sp2n}.
 \Absatz

(b) Let $\beta\in\SP_{2n}(R)_{Rss}$. We decompose
  $N=R^{2n}=N_+\oplus N_-\oplus N_*$ as in lemma \ref{Abspaltungplusminus1}(a). Since this decomposition is
  $J_{2n}$-orthogonal the restrictions $p_+,p_-,p_*$ of the symplectic form $J_{2n}$ to $N_+,N_-$ and
  $N_*$ are unimodular, so these spaces have even rank:
  $N_+\simeq R^{2r_+},\, N_-\simeq R^{2r_-},\, N_*\simeq R^{2g}$. If we view $p_*$ as skew symmetric matrix
  and $\beta_*\in \SP(p_*)\subset \SL_{2g}$ we can form the symmetric matrix (bilinear form)
  $q_*=p_*\cdot(\beta_*-\beta_*^{-1})^{-1}$ and get $\beta_*\in \SO(q_*)$. For $\epsilon_{\pm}\in
  R^\times/(R^\times)^2$ we consider the symmetric bilinear forms $q_+=\epsilon_+\cdot J_{1+2r_+}$
  on $M_+=R^{1+2r_+}$ and $q_-=\epsilon_-\cdot J_{1+2r_-}$ on $M_-=R^{1+2r_-}$. By lemma
  \ref{orthogonaleNormalform} there are two different choices of pairs $(\epsilon_+,\epsilon_-)$
  such that the quadratic space $(M,q)=(M_+,q_+)\oplus (M_-,q_-)\oplus (N_*,q_*)$ is isomorphic
  to the standard space $(R^{2n+2},J'_{2n+2})$. For these two choices the element
  $b=id_{M_+}\times (-id_{M_-})\times \beta_*\in \ORTH(q,R)^-$ can be viewed as an element of
  $\ORTH_{2n+2}(R)^-_{Rss}$, which maps to $\beta$ under $\NNN$. It is clear from the constructions
  that the two classes just obtained are all $\SO_{2n+2}(R)$-conjugacy classes in
  $\ORTH_{2n+2}(R)^-_{Rss}$ mapping to $\beta$ under $\NNN$.

\qed

\begin{lemma}\label{O2n+2konjbereitsueberR}
 Let $\gamma_1\in \ORTH_{2n+2}(\OOO_F)^-$  be
 $R$-$\Theta$-semisimple.
 \begin{itemize}
  \item[(a)] If $\gamma_2:=g_F^{-1}\cdot \gamma_1\cdot g_F\in \ORTH_{2n+2}(\OOO_F)^-$
    for $g_F\in \SO_{2n+2}(F)$ there exists $g_R\in \SO_{2n+2}(\OOO_F)$ with
    $\gamma_2=g_R^{-1}\cdot \gamma_1\cdot g_R$.
 \item[(b)] There is a unique $\SO_{2n+2}(\OOO_F)$-conjugacy class $\{\gamma'_1\}$ in
    $\SO_{2n+2}(\OOO_F)$ different from $\{\gamma_1\}$ such that for every $g_F\in \SO_{2n+2}(\bar F)$ with
    $\gamma_2:=g_F^{-1}\cdot \gamma_1\cdot g_F\in \ORTH_{2n+2}(\OOO_F)^-$
    there either exists $g_R\in \SO_{2n+2}(\OOO_F)$ with
    $\gamma_2=g_R^{-1}\cdot \gamma_1\cdot g_R$ or $g'_R\in \SO_{2n+2}(\OOO_F)$ with
    $\gamma_2=(g'_R)^{-1}\cdot \gamma'_1\cdot g'_R$
 \end{itemize}
\end{lemma}
Proof: (a) Let $M=\OOO_F^{2n+2}=M_{+,i}\oplus M_{-,i}\oplus M_{*,i}$ for $i=1,2$ be the orthogonal
  decompositions with respect to $\gamma_i$ as in lemma \ref{Abspaltungplusminus1}(b) and let
  $q_{\pm,i},q_{*,i}$ denote the restrictions of the standard form $J'_{2n+2}$ to these subspaces.
  Since $M_{\pm,i}$ are eigenspaces of $\gamma_i$ we have
  \begin{align*}
    g_F(M_{\pm,2}\otimes_{\OOO_F} F)\quad &= \quad M_{\pm,1}\otimes_{\OOO_F} F \qquad
  \end{align*}
  Since $g_F\in\SO_{2n+2}(F)$ we get for the orthogonal complement:
  \begin{align*}
    g_F(M_{*,2}\otimes_{\OOO_F} F)\quad &= \quad M_{*,1}\otimes_{\OOO_F} F.
  \end{align*}
  In fact $g_F$ induces isomorphisms of quadratic spaces:
  \begin{align*} 
    \left(M_{\pm,2}\otimes_{\OOO_F}F, q_{\pm,2}\right)\quad &\tilde\to\quad
     \left(M_{\pm,1}\otimes_{\OOO_F}F,q_{\pm,1}\right)\quad \text{ and }\quad \\
    g_{F,*}:\quad\left(M_{*,2}\otimes_{\OOO_F}F, q_{*,2}\right)\quad &\tilde\to\quad
     \left(M_{*,1}\otimes_{\OOO_F}F,q_{*,1}\right).
  \end{align*}
  Since the quadratic spaces are defined over $\OOO_F$ and become isomorphic over $F$ and since
  the forms are unimodular, the spaces are isomorphic over $\OOO_F$ by lemma \ref{orthogonaleNormalform},
  i.e. there exists  $g'_R\in \SO_{2n+2}(\OOO_F)$ inducing isomorphisms
  \begin{align*}
    \left(M_{\pm,2}, q_{\pm,2}\right)\quad \tilde\to\quad\left(M_{\pm,1},q_{\pm,1}\right)\quad
    \text{ and }\quad g'_{R,*}:\quad\left(M_{*,2}, q_{*,2}\right)\quad \tilde\to\quad \left(M_{*,1},q_{*,1}\right).
  \end{align*}
 \klabsatz
  If   $\gamma_{*,i}$ denotes the restriction of $\gamma_i$ to $M_{*,i}$ we get
  \begin{align}
    \label{gamma2}\gamma_{*,2}\quad&=\quad g_{F,*}^{-1}\cdot \gamma_{*,1}\cdot g_{F,*}\qquad\text{ and }\\
    \label{gamma3} \gamma_{*,3}\quad:&=\quad g'_{R,*}\cdot\gamma_{*,2}\cdot(g'_{R,*})^{-1} \in \SO(q_{*,1},\OOO_F).
  \end{align}
  We have $g_{F,*}\cdot (g'_{R,*})^{-1}\in \SO(q_{*,1},F)$.
  Now it follows from (\ref{gamma2}), (\ref{gamma3}) and lemma \ref{SpkonjbereitsueberR} that there
  exists $g_*\in \SO(M_{*,1},q_{*,1})$ satisfying $g_*\cdot \gamma_{*,3}\cdot g_*^{-1}=\gamma_{*,1}$. Then
  $g_R:=(id_{M_{+,1}}\times id_{M_{-,1}}\times g_*)\cdot g'_R\in \SO_{2n+2}(\OOO_F)$ satisfies
  $g_R\cdot \gamma_2\cdot g_R^{-1} =\gamma_1$.
\absatz
(b) Let us assume that $g_F\in \SO_{2n+2}(\bar F)$ satisfies
  $\gamma_2:=g_F^{-1}\cdot \gamma_1\cdot g_F\in \ORTH_{2n+2}(\OOO_F)^-$.
  We only know that the quadratic spaces become isomorphic over $\bar F$, but we have the
  additional discriminant conditions $\Delta(q_{+,1})\cdot \Delta(q_{-,1})\cdot \Delta(q_{*,1})
  =\Delta(q_{+,2})\cdot \Delta(q_{-,2})\cdot \Delta(q_{*,2})$ and $\Delta(q_{*,1})=\Delta(p_{*,1})
  \cdot \det(\gamma_{*,1}-\gamma_{*,1}^{-1})=1\cdot \det(g_{F,*}(\gamma_{*,2}-\gamma_{*,2}^{-1})g_{F,*}^{-1})
  = \Delta(p_{*,2})\cdot \det(\gamma_{*,2}-\gamma_{*,2}^{-1})=\Delta(q_{*,2})$ in
  $\OOO_F^\times/(\OOO_F^\times)^2$, where we use the fact that the
  $p_{*,i}:=q_{*,i}\cdot(\gamma_{*,i}-\gamma_{*,i}^{-1})$
  are unimodular skew symmetric by lemma \ref{SOSPschieben} and thus have square determinants.
  The isomorphy-typ of the quadratic spaces $\left(M_{\pm ,1},q_{\pm ,1}\right),
  \left(M_{*,1},q_{*,1}\right)$ being fixed this means that there are two choices for the equivalence class of
  $q_{+,2}$ but the isomorphy-typ of the other quadratic spaces $(M_{-,2},q_{-,2})$ and $(M_{*,2},q_{*,2})$
  are then uniquely determined. To construct $\gamma'_1$ we change the quadratic
  forms $q_{+,1}$ on $M_{+,1}$ and $q_{-,1}$ on $M_{-,1}$ to the other isomorphy-typ but make no change for $M_{*,1}$,
  consider an isomorphism  of quadratic spaces   $\iota:R^{2n+2}\quad\tilde\to\quad M_{+,1}\oplus M_{-,1}\oplus M_{*,1}$
  with respect to these modified forms on $M_{\pm,1},M_{*,1}$ and the standard form $J'_{2n+2}$ on $R^{2n+2}$,
  and put finally   $\gamma'_1=\iota^{-1}\circ \gamma_1\circ\iota$. The statement of (b) now follows
  as in part (a).
\qed

\begin{lemma}
\label{MatchinglemmaSOSP}
  In the notations of \ref{explizitesMatchingSOSP} let $b\in \ORTH_{2n+2}(\OOO_F)^-$ be residually
  semi\-simple and $\beta=1_{2r_+}\times(- 1_{2r_-})\times b_*\in \Sp_{2n}(\OOO_F)$
  a representing element of $\NNN(b)$ with $b_*-b_*^{-1}\in \GL_{2g}(\OOO_F)$.\klabsatz
  Assume we have $BC$-matching topologically unipotent elements
  $u_+\in \SO_{2r_++1}(F)$ and  $v_+\in  \Sp_{2r_+}(F)$ resp.
  $u_-\in \SO_{2r_-+1}(F)$ and  $v_-\in  \Sp_{2r_-}(F)$
  and an additional topologically unipotent element
  $u_*\in Cent(b_*,\SO(q_*,F))\, \simeq\, Cent(b_*,\Sp(p_*,F))$.
   We form the topologically unipotent elements $u=u_+\times u_-\times u_*\in Cent(b,\SO_{2n+2}(F))$
  and $v=v_+\times v_-\times u_*\in Cent(\beta,\SP_{2n}(F))$.\klabsatz
 Then the elements
  $g:= bu=ub\in \ORTH_{2n+2}(F)^-$ and $\gamma:=\beta v=v\beta\in\SP_{2n}(F)$ match.
\end{lemma}
 Proof: As in the proof of lemma \ref{Matchinglemma}  we work
 in the case $F=\bar F$ and assume that $g$ resp. $\gamma$ lie in the diagonal tori. The same
 holds for the residually semisimple parts $b$ resp. $\beta$ and the topologically unipotent parts
 $u$ and $v$. As the matching of $b$ and $\beta$ is already proved
 in  \ref{propoSOSP} we only have to examine the topologically
 unipotent elements. We can arrange the diagonal matrices $u_\pm\in \SO(q_\pm,F)$
 such that their middle entries $1$ correspond to the eigenvectors $t_{n+1}\cdot e_{n+1}\pm e_{n+2}\in M_\pm$,
 which get lost by the construction of $N_\pm$. Then the claim follows immediately from the
 definition of $BC$-matching \ref{BCMatching}, example \ref{BspMatchingSO2n+2Sp2n} and the constructions
 in the proof of proposition \ref{propoSOSP}.
  \qed

\begin{remark} \em  The surjectivity statement of  Proposition \ref{propoSOSP}(b) is not true if $R$ is a field,
  for example a $p$-adic field $F$: Let $\Delta\in F^*$ denote a non square and
  \begin{align*} \beta\quad&=\quad
  \begin{pmatrix} a_1 & & & b_1\\ & a_2&b_2&\\  &b_2\Delta &a_2&\\ b_1\Delta& & &a_1\end{pmatrix}
  \quad\in \quad \Sp_4(F)\qquad\text{ where } \\
    a_i\quad&=\quad \frac{\lambda_i^2+\Delta}{\lambda_i^2-\Delta},\qquad
    b_i\quad=\quad \frac{2\cdot\lambda_i}{\lambda_i^2-\Delta}\qquad \text{ for }\lambda_i\in F^*,\quad i=1,2.
  \end{align*}
  Then we have $N_*=N$ and $\beta_*=\beta$ for $N=F^4$ and can compute
  \begin{align*}
    q_*\quad:&=\quad p_*\cdot(\beta_*-\beta_*^{-1})^{-1}\quad=\quad
    J_4\cdot antidiag(2b_1,2b_2,2\Delta b_2,2\Delta b_1)^{-1}\\&=\quad
    diag\left( \frac{1}{2b_1}, \frac{-1}{2b_2},\frac{1}{2\Delta b_2},\frac{-1}{2\Delta b_1}\right)
    \quad=\quad\frac{-1}{2\Delta b_1}\cdot diag\left(-\Delta,\Delta\cdot\frac{b_1}{b_2},-\frac{b_1}{b_2},1\right).
  \end{align*}
  Thus the quadratic form $q_*$ on $N$ is anisotropic if $b_1\cdot b_2^{-1}$ is not a norm of the
  extension $F\sqrt\Delta/F$. In this case $(N,q_*)$ cannot be obtained as direct summand of the six dimensional
  quadratic split space $(F^6,J'_6)$. The considerations of \ref{explizitesMatchingSOSP} and \ref{propoSOSP}(b)
  then show that the conjugacy class of $\beta$ is not in the image of $\NNN$.
  \klabsatz

  The following theorem is again the fundamental lemma
  for a stable endoscopic lift modulo the $BC$-conjecture. But the non surjectivity of $\NNN$ in the case
  of local fields forces us to include the vanishing statement
  for orbital integrals of elements, that do not match.
\end{remark}

\begin{thm}\label{FulemmaSO2n+2SP2n}
  Assume that conjecture ($\text{BC}_{\text{m}}$) is true for all $m\le n$.
 \begin{itemize}
 \item[(a)]
  If   $g\in\widetilde{\SO}_{2n+2}(F)=\ORTH_{2n+2}(F)$ with $\det(g)=-1$ and $\gamma\in \Sp_{2n}(F)$ are
  matching semisimple elements  then we have
 \begin{gather}\label{Fule2n+22nGleichung}
  \STO_{g}(1,\widetilde{\SO}_{2n+2})\quad=\quad\STO_\gamma(1, \Sp_{2n}).
 \end{gather}
 \item[(b)] If the semisimple $\gamma\in \Sp_{2n}(F)$ matches with no element of $\widetilde{\SO}_{2n+2}(F)$, then we have
 $\STO_\gamma(1, \Sp_{2n})=0$.
 \end{itemize}
\end{thm}
Proof: Since the proof of (a) is similar to the proofs of theorems \ref{FulemmaSp2nPGL2n+1} and
\ref{FulemmaGSpin2n+1GL2n} we will be sketchy in some steps. We remark that (b) is an immediate
corollary of the considerations in Step 1: If the right hand side of (\ref{Fule2n+22nGleichung}) does not vanish
one can construct an element $g\in\widetilde{\SO}_{2n+2}(F)$ matching with $\gamma$ using proposition \ref{propoSOSP}
and lemma \ref{MatchinglemmaSOSP}.
\absatz

\underline{Step 1} As in the cited proofs we may assume without
 loss of generality that $g\in \ORTH_{2n+2}(\OOO_F)$ and $\gamma\in \Sp_{2n}(\OOO_F)$.
 We may furthermore assume that $g=(u_+,u_-,u_*)\cdot s$ and $\gamma=(v_+,v_-,u_*)\cdot \sigma$
 are the topological Jordan decompositions with $BC$-matching topologically unipotent
 $u_+\in \SO_{2r_++1}(\OOO_F) $ and $v_+\in \SP_{2r_+}(\OOO_F)$ respectively
 $u_-\in \SO_{2r_-+1}(\OOO_F)$ and $v_-\in \SP_{2r_-}(\OOO_F)$, matching
 residually semisimple $s\in\ORTH_{2n+2}(\OOO_F)$ and $\sigma\in \SP_{2n}(\OOO_F)$ and
 topologically unipotent $u_*\in Cent(\sigma_*,\SP_{2g})(\OOO_F)$.
 \absatz

\underline{Step 2} As in \ref{FulemmaSp2nPGL2n+1}
  we obtain all relevant conjugacy classes  in the stable conjugacy class of $\gamma$ if we let
  $v'_+$ resp. $v'_-$ vary through a set of representatives for the conjugacy classes in the stable
  conjugacy class of $v_+$  resp. $v_-$ in  $\Sp_{2r_+}(F)$ resp. $\Sp_{2r_-}(F)$  and $u'_*$
  through a set of representatives for the conjugacy classes in the stable conjugacy class of $u_*$
  in $Cent(\sigma_*,\Sp_{2g})(F)$ and then consider all $\gamma'= \sigma\cdot (v'_+,v'_-,u'_*)$ i.e.
\begin{align}\label{Sporbital2}
  &\STO_{\gamma}(1,\Sp_{2n})\quad \\
  \notag  =&\quad \sum_{v'_+\sim v_+} O_{v'_+}(1,\Sp_{2r_+})\cdot\sum_{v'_-\sim v_-} O_{v'_-}(1,\Sp_{2r_-})
             \cdot \sum_{u'_*\sim u_*} O_{u'_*}(1,Cent(\sigma_*,\Sp_{2g})).
\end{align}
\underline{Step 3} In the $\Theta$-twisted situation $\ORTH_{2n+2}$ all relevant $\Theta$-conjugacy
 classes are of the form $g'=s'\cdot(u'_+,u'_-,u'_*)$ where $u'_*$ is as above, $u'_\pm$ vary
 through a set of representatives for the conjugacy classes in the stable class of
 $u_\pm\in \SO_{2r_\pm +1}(F)$
 and $s'$ is either $s$ or a representative for the corresponding other conjugacy class $s''$ as in
 lemma \ref{O2n+2konjbereitsueberR}(b). Observe the centralizers $\SO_{2n+2}^s$ and
 $\SO_{2n+2}^{s''}$ can be identified, since the two equivalence classes of symmetric
 unimodular forms on a free $\OOO_F$-module of odd rank have representatives which are
 scalar multiples of each other. Therefore we can use the same collections of $u'_\pm$ and
 $u'_*$ for $s$ as for $s''$. The appearance of $s''$ thus introduces just an additional
 factor $2$ in the computation.  But since the centralizers  $\SO_{2n+2}^s$ and
 $\SO_{2n+2}^{s''}$  have two connected components, there appears an additional factor
 $\frac{1}{2}$ when we apply the Kazhdan-lemma \ref{Kazhdan-Lemma}. Thus we get:
 \begin{align}\label{Oorbital}
  &\STO_{g}(1, \ORTH_{2n+2})\quad = \quad
   \sum_{(u'_+,u'_-,u'_*)\sim (u_+,u_-,u_*)} O_{(u'_+,u'_-,u'_*)}(1,\SO_{2n+2}^{s})\\
  \notag=\; \sum_{u'_+\sim u_+} &O_{u'_+}(1,\SO_{2r_++1})\cdot
   \sum_{u'_-\sim u_-} O_{u'_-}(1,\SO_{2r_-+1})\cdot
   \sum_{u'_*\sim u_*} O_{u'_*}(1,Cent(\sigma_*,\Sp_{2g})).
 \end{align}
 In the last step we applied lemma \ref{EndlichesZentrum} in the situation
 $G=\SO_{2n+2}^s=G_1\times \{\pm 1_{2n+2}\}$ where
 $G_1=\SO_{2r_++1}\times\SO_{2r_-+1}\times Cent(\sigma_*,\SP_{2g})$.
 \absatz
\underline{Step 4} Since $v_\pm$ and $u_\pm$ are $BC$-matching we can apply ($\text{BC}_{r_\pm}$) to get
that the right hand sides of  (\ref{Sporbital2}) and (\ref{Oorbital}) coincide, which proves the theorem.

\qed

   \Seitenumbruch

\appendix
\renewcommand{\Numerierung}{\refstepcounter{meinzaehler}{\bf(\Alph{chapter}.\arabic{meinzaehler}) }}
\renewcommand{\meinchapter}[1]{\refstepcounter{chapter}\section*{\Alph{chapter} $\ $ #1}
    \addcontentsline{toc}{chapter}{\Alph{chapter}$\ $ #1}
    \markboth{}{\scshape #1}
    }

{ 
\small
\setlength{\itemsep}{0pt}
\setlength{\parsep}{0pt}
\setlength{\topsep}{0pt}
\setlength{\partopsep}{0pt}
\setlength{\baselineskip}{0pt}
\setlength{\smallskipamount}{0pt}
\setlength{\medskipamount}{0pt}
\setlength{\bigskipamount}{0pt}

}
\Absatz\Absatz

Mathematisches Institut \par
Im Neuenheimer Feld 288\klabsatz
D-- 69120 Heidelberg
\Absatz

weissauer@mathi.uni-heidelberg.de\par
weselman@mathi.uni-heidelberg.de\par
  

\end{document}